\documentclass[a4paper,12pt]{amsart}
\usepackage{amsmath}
\usepackage{amsfonts}
\usepackage{amssymb}
\usepackage{amsbsy}

\usepackage{amscd}
\usepackage{verbatim}
\usepackage[pdftex]{graphicx}
\def\cal{\mathcal}

\def\frak{\mathfrak}

\newenvironment{pf*}[1]{\proof[#1]}{\endproof}
\newtheorem{Theorem}[equation]{Theorem}
\newtheorem{Corollary}[equation]{Corollary}
\newtheorem{Lemma}[equation]{Lemma}
\newtheorem{Proposition}[equation]{Proposition}

\theoremstyle{definition}
\newtheorem{Definition}[equation]{Definition}

\numberwithin{equation}{section}

\newcommand{\meji}[1]{}
\newcommand{\submeji}[1]{}

\newcommand{\mejila}[1]{}

\begin{document}
\baselineskip=16pt
\pagestyle{headings}
\begin{center}
\title[
Simple and explicit constructions of semi-discrete surfaces] 
{Simple and explicit constructions of semi-discrete surfaces and discrete surfaces}

\vspace{3mm} 

\end{center}
\author{Kenji KAJIWARA, Shota SHIGETOMI and Seiichi UDAGAWA}

\address{Kenji KAJIWARA : Institute of Mathematics for Industry, Kyushu University, 744 Motooka, Fukuoka 819-0395, Japan}
\email{kaji@imi.kyushu-u.ac.jp}
\thanks{The first named author is partly supported by 
JST CREST Grant Number JPMJCR1911, JSPS Grant-in-Aid for 
Scientific Research (C) No.~21K03329.}

\address{Shota SHIGETOMI : Institute of Mathematics for Industry, Kyushu University, 744 Motooka, Fukuoka 819-0395, Japan}
\email{s-shigetomi@imi.kyushu-u.ac.jp}

\address{Seiichi UDAGAWA : Nihon University, School of Medicine,
Department of Natural Sciences, Division of Mathematics, Itabashi,
Tokyo 173-0032, Japan}
\email{udagawa.seiichi@nihon-u.ac.jp}
\thanks{The third named author is partly supported by JSPS Grant-in-Aid for Scientific Research (C) No.~21K03235.}

\maketitle
\par\noindent 

\vspace{10pt}
\begin{abstract}
We give a simple and explicit constructions of various
semi-discrete surfaces and discrete $K$-surfaces in terms of 
the Jacobi elliptic functions using $\tau$-functions. 
Their periodicities are also determined.
\end{abstract}
\vspace{15pt}

\centerline{\bf Introduction}\label{Section0}

The discrete sine-Gordon equation was introduced by R. Hirota(\cite{hirota}).
After that, as an application of discrete sine-Gordon equation
to the discrete differential geometry,
significant progress was done by Bobenko-Pinkall in \cite{bobpink}.
They gave the solution of discrete sine-Gordon equation in terms of
the Riemann theta function and gave a formulae for the immersion of 
the discrete $K$-surface.
On one hand, the industrial mathematics researchers has been studying
the discrete space curve and the isoperimetric deformation of the 
discrete space curve
(for example, see \cite{ikmo},\cite{kkp}),
of which the integrability is given by the semi-discrete 
sine-Gordon equation. 
The theory of discrete space curves have 
recently attracted attention as equations of motion 
that describe the Kaleidocycle(see \cite{kkp}, \cite{kks}).
The Kaleidocycle is demonstrated by a deformation of a discrete space curve
which is periodic with respect to the integer parameter.
Very recently, some examples of the Kaleidocycle are given in
\cite{kks} and they are described in terms of the Jacobi theta functions.
On the other hand, in \cite{udagawa},
we gave solutions in terms of the Jacobi elliptic functions for
both the discrete and semi-discrete sine-Gordon equations.
They are said to be \lq\lq {\rm dn}-solution\rq\rq because
the cosine of the solution function may be described in terms of the
Jacobi {\rm dn}-function. In fact, although
the solution in terms of Riemann
theta function has been given in \cite{bobpink},
it is not so clear to deduce the elliptic function solution from that.
It is because the Riemann period matrix is differeent from the ordinary form
and the expressions using Riemann theta function 
for the Jacobi elliptic functions are different from the ordinary one.
For this problem, in \cite{udagawa}, we investigate the spectral curve
for the sine-Godron equation and obtain a Weierstrass ${\frak p}$-function
and determined the two periods of it, which teachs one how the form of
the Riemann period matrix is.  This knowledge made it possible
to clarify the relationship between the elliptic solution and the Riemann
theta function solution (see \cite{udagawa}). Based on these observations,
we may find $\tau$-functions for constructing semi-discrete surface
(for the definition, see the preliminaries) which correspond to
the Baker-Akhiezer functions in the theory of 
Bobenko-Pinkall(\cite{bobpink}).

In this paper, first of all, we shall give several definitions of
the discrete objects which we are interested in.
Next, we give \lq\lq {\rm cn}-solution\rq\rq
of the semi-discrete and discrete sine-Gordon equations, that is,
the cosine of the solution function may be described in terms of 
the Jacobi {\rm cn}-function, which are done in \S1 and \S2.
In \S3, we shall give a general method of construction of semi-discrete
surface in terms of $\tau$-functions.
In \S4, we shall give the $\tau$-functions which are described 
in terms of the ceratin Jacobi theta functions and give interesting
examples of the semi-discrete surfaces which corresponds to the
${\rm dn}$-solution of the semi-discrete SG-equation.
In \S5, choosing another $\tau$-functions which are described in terms of 
the certain Jacobi theta functions, we give another example of the 
semi-discrete surface which corresponds to the
${\rm cn}$-solution of the semi-discrete SG-equation.
Some periodic solutions among the examples constructed in \S4 and \S5
yield the model of the Kaleidocycle. 
In \S6, we shall give interesting examples of discrete $K$-surfaces 
in the sense of \cite{bobpink} using the results in \S4 and \S5. 
\S7 is Appendix, which is devoted to demonstrate
the well-known properties of Jacobi theta functions
and the properties of the elliptic functions needed
in \S1$\sim$\S5.
\vspace{10pt}

\centerline{\bf Preliminaries}

First of all, we shall give the several difinition of discrete objects
which we are interested in.  
\begin{Definition}
\item{\rm (1)} A map $\Gamma : \mathbb{Z}\ni m \longrightarrow \Gamma_{m}\in 
\mathbb{R}^{3}$ is said to be a {\it discrete space curve} if
$\Gamma_{m+1}\not=\Gamma_{m}$ for any $m$ and 
any consecutive three points $\Gamma_{m+1}, \Gamma_{m}$ and
$\Gamma_{m-1}$ are not collinear.
\item{\rm (2)} A map $\Gamma : \mathbb{Z}\times\mathbb{R}\ni (m, t) 
\longrightarrow
\Gamma(m,t)\in \mathbb{R}^{3}$ is said to be
a {\it semi-discrete surface} if $\Gamma_{m}(t):=\Gamma(m,t)$ 
is a deformation of the discrete space curve 
$m \longrightarrow \Gamma_{m}$
with the deformation parameter $t$.
\item{\rm (3)} A map $\Gamma : \mathbb{Z}\times\mathbb{Z}\ni (m, n) 
\longrightarrow \Gamma(m,n)=:\Gamma_{m,n}\in
\mathbb{R}^{3}$ is said to be
{\it discrete K-surface} if 
(i) every point $\Gamma_{m,n}$ and its neighbours $\Gamma_{m+1,n},
\Gamma_{m-1,n},\Gamma_{m,n+1},\Gamma_{m,n-1}$ belong to
one plane ${\cal P}_{m,n}$, (ii) the lengths of the opposite
edges of an elementary quadrilateral are equal
\begin{eqnarray*}
\begin{array}{ll}
|\Gamma_{m+1,n}-\Gamma_{m,n}|
&=|\Gamma_{m+1,n+1}-\Gamma_{m,n+1}|=A_{m}\not=0 ,\\
|\Gamma_{m,n+1}-\Gamma_{m,n}|
&=|\Gamma_{m+1,n+1}-\Gamma_{m+1,n}|=B_{n}\not=0 ,\\
\end{array}
\end{eqnarray*}
where $A_{m}$ does not depend on $n$ and $B_{n}$ not on $m$.
\end{Definition}
The definition of the discrete $K$-surface is due to Bobenko-Pinkall
(\cite{bobpink}).
For a while from now, we are constrate on the semi-discrete surface
with the deformation parameter $t$.
We may construct a Frenet frame $(T_{m}\hskip .1cm 
N_{m}\hskip .1cm B_{m})\in SO(3)$ associated to $\Gamma_{m}(t)$ 
as follows. Define $T_{m}$ by $T_{m}=\dfrac{\Gamma_{m+1}-\Gamma_{m}}{
\varepsilon_{m}}$, where $\varepsilon_{m}=|\Gamma_{m+1}-\Gamma_{m}|$.
Choose $B_{m}$ so that $\Gamma_{m+1}-\Gamma_{m}=a_{m} B_{m+1}\times B_{m}$
for some constant $a_{m}$ for each $m$.
We assume that $\Gamma_{m}$ has constant speed. We then see that  
$\varepsilon_{m}$ is constant and independent of $m$, which we denote by
just $\varepsilon$. Finally, we define $N_{m}$ by
$N_{m}=B_{m}\times T_{m}$. The signed curvature ${\cal K}_{m}$ and
the torsion angle $\nu_{m}$ is defined by
\begin{equation*}
\left<T_{m},T_{m-1}\right>=\cos{\cal K}_{m},\quad
\left<B_{m},B_{m-1}\right>=\cos\nu_{m},\quad
\left<B_{m}, N_{m-1}\right>=\sin\nu_{m} .
\end{equation*}
If we set $\widetilde{\Phi}_{m}=(T_{m}\hskip .1cm 
N_{m}\hskip .1cm B_{m})$, we obtain $\widetilde{\Phi}_{m+1}=
\widetilde{\Phi}_{m}\widetilde{L}_{m}$, where
\begin{equation}\label{trans}
\widetilde{L}_{m}=\left(\begin{matrix}
\cos{\cal K}_{m+1}&-\sin{\cal K}_{m+1}&0\\
\cos\nu_{m+1}\sin{\cal K}_{m+1}&\cos\nu_{m+1}\cos{\cal K}_{m+1}&
\sin\nu_{m+1}\\
-\sin\nu_{m+1}\sin{\cal K}_{m+1}&-\sin\nu_{m+1}\cos{\cal K}_{m+1}&
\cos\nu_{m+1}\\
\end{matrix}\right).
\end{equation}
(see \cite{ikmo}, \cite{kkp}). 
Moreover, we assume that $\nu_{m}$ is constant and independent of
$m$, which we denote by just $\nu$.
Under this situation, the deformation is said to be {\it isoperimetric}
if there is some constant $\rho$,
which is independent of $m$, and some function $w_{m}(t)$ 
such that the following holds.
\begin{equation}
\dfrac{d\Gamma_{m}}{dt}=\rho\left(\cos w_{m}T_{m}+\sin w_{m}N_{m}\right).
\end{equation}
Then, we obtain ${\cal K}_{m}=-w_{m}\pm w_{m-1}$ (see \cite{kkp}).
In the case of ${\cal K}_{m}=-w_{m}+w_{m-1}$, ${\cal K}_{m}$ and
$w_{m}$ are controlled by the solution of the semi-discrete SG-equation.
On the other hand, 
in the case of ${\cal K}_{m}=-w_{m}-w_{m-1}$, ${\cal K}_{m}$ and
$w_{m}$ are controlled by the solution of the semi-discrete 
potential mKdV-equation (see Remark 2 in \cite{kkp}).
However, since our solution below to the semi-discrete SG-equation 
satisfies also the potential mKdV-equation (see \S1),
our solution covers both cases of the relations between ${\cal K}_{m}$
and $w_{m}$.

To obtain the Weierstrass ${\frak p}$-function stated in the introduction,
we consider the spectral curve for the smooth SG-equation.
For a function $\hat{w}=\hat{w}(x,t)$ of the parameters $x$ and $t$,  
the smooth SG-equation is given by $\hat{w}_{xt}=-4\sin(\hat{w})$,
where $\hat{w}_{xt}=\partial_{t}\partial_{x}\hat{w}$. 
Set $u=\frac{x+t}{2}, v=\frac{x-t}{2}$ and assume that $\hat{w}=\hat{w}(u)$,
a function which depends on only parameter $u$.
We then rewrite the SG-equation as $\hat{w}_{uu}=-16\sin\hat{w}
=8i\left(\exp(i\hat{w})-\exp(-i\hat{w})\right)$.
Multiplying it by $\hat{w}_{u}$ and integrating it, we obtain
\begin{equation*}
\left(\hat{w}_{u}\right)^{2}=16\exp(i\hat{w})+16\exp(-i\hat{w})+C,
\end{equation*}
where $C$ is an integration constant and determined by an initial condition.
We give the initial condition by $\hat{w}(0)=0, \hat{w}_{u}(0)=8k$,
where $k$ is the modulus of the Jacobi {\rm dn}-function because
$\hat{w}=2\arccos({\rm dn}(4u))$ is a ${\rm dn}$-solution of the smooth
SG-equation. By the initial condition, we see that $C=32(2k^{2}-1)$.
If we set $Y=\exp(i\hat{w})$, we then obtain 
$\left(\dfrac{dY}{du}\right)^{2}=-16Y(Y^{2}+2(2k^{2}-1)Y+1)$.
If we set $2iB=\dfrac{dY}{du}, \eta=Y+\frac{2}{3}(2k^{2}-1)$,
we then have
\begin{equation*}
B^{2}=4(\eta-e_{1})(\eta-e_{2})(\eta-e_{3}),
\end{equation*}
where $e_{1}=\frac{2}{3}(2k^{2}-1), e_{2}=-\frac{1}{3}(2k^{2}-1)-2kk^{\prime}i,
e_{3}=-\frac{1}{3}(2k^{2}-1)+2kk^{\prime}i$, and
$k^{\prime}=\sqrt{1-k^{2}}$. This is a standard form of the elliptic curve
and parametrized by $\eta={\frak p}(z), B={\frak p}^{\prime}(z)$,
where ${\frak p}(z)$ is a Weierstrass ${\frak p}$-function.
We may express it in terms of the Jacobi elliptic functions as follows.
\begin{Lemma}{\rm (\cite{udagawa})}
\begin{equation}\label{WP}
{\frak p}(z)=\left({\rm dn}(2iz+iK^{\prime})-ik{\rm sn}(2iz+iK^{\prime})
\right)^{2}+e_{1}.
\end{equation}
\end{Lemma}
We find by~(\ref{WP}) that the two periods $\left\{2\omega, 2\omega^{\prime}
\right\}$ of ${\frak p}(z)$ are given by
\begin{equation}
2\omega=2K^{\prime} \qquad 2\omega^{\prime}=iK+K^{\prime} .
\end{equation}
The holomorphic 1-form $\hat{\omega}$ and 
the normalized abelian differential of second kind $\Omega$ on the spectral elliptic curve are given by
\begin{equation}
\hat{\omega}=\dfrac{\pi i}{\omega}\dfrac{d\eta}{B} ,\quad
\Omega=-\left({\frak p}(z)+\dfrac{\zeta_{W}(\omega)}{\omega}\right)
\dfrac{d\eta}{B},
\end{equation}
where $\zeta_{W}$ is the Weierstrass $\zeta$-function
(see \cite{udagawa}).

\vspace{20pt}

\section{Solutions of semi-discrete sine-Gordon equation in terms of the Jacobi elliptic function}

In \cite{udagawa}, it is proved that $\hat{w}_{m}=2\arcsin(k{\rm sn}
(4K\xi_{m}))$ is a solution of the semi-discrete sine-Gordon equation
\begin{equation*}
\frac{d}{dt}\hat{w}_{m+1}-\frac{d}{dt}\hat{w}_{m}
=\widetilde{\alpha}\sin\left(\frac{1}{2}\hat{w}_{m+1}+
\frac{1}{2}\hat{w}_{m}\right) ,
\end{equation*}
where $\xi_{m}=m\Omega+\xi_{0}+At$ 
with some constants $\Omega, \xi_{0}, A$,
and $\widetilde{\alpha}$ is given by
\begin{equation*}
\widetilde{\alpha}=-8KA\dfrac{{\rm sn}(2K\Omega){\rm dn}
(2K\Omega)}{{\rm cn}(2K\Omega)}.
\end{equation*}
Moreover, $\hat{w}_{m}=2\arcsin(k{\rm sn}
(4K\xi_{m}))$ also satisfies the mKdv equation
\begin{equation*}
\frac{d}{dt}\hat{w}_{m+1}+\frac{d}{dt}\hat{w}_{m}
=\widetilde{\beta}\sin\left(\frac{1}{2}\hat{w}_{m+1}-
\frac{1}{2}\hat{w}_{m}\right) ,
\end{equation*}
where
\begin{equation*}
\widetilde{\beta}=8KA\dfrac{{\rm cn}(2K\Omega)}
{{\rm sn}(2K\Omega){\rm dn}(2K\Omega)}.
\end{equation*}
Since we may write $\hat{w}_{m}$ as
$\hat{w}_{m}=2\arccos({\rm dn}(4K\xi_{m}))$,
we call this solution \lq\lq {\rm dn}-solution\rq\rq.
We here prove that there is also \lq\lq {\rm cn}-solution
\rq\rq for the sine-Gordon equation.
For this, we take $f_{m}=\theta_{3}(\xi_{m})\theta_{1}(
\xi_{m})$ and $g_{m}=\theta_{0}(\xi_{m})\theta_{2}(\xi_{m})$ and 
consider $\tan(\frac{\hat{w}_{m}}{4})
=\dfrac{g_{m}}{f_{m}}$. 
We then calculate $\sin(\frac{\hat{w}_{m}}{2})$ using (i) of Lemma~\ref{App1}
as follows.
\begin{eqnarray*}
\begin{array}{ll}
\sin(\frac{\hat{w}_{m}}{2})&=
\dfrac{2\tan(\frac{\hat{w}_{m}}{4})}{
1+\tan^{2}(\frac{\hat{w}_{m}}{4})}
=\dfrac{2f_{m}g_{m}}{f_{m}^{2}+g_{m}^{2}}
=\dfrac{2\theta_{0}(\xi_{m})\theta_{1}(\xi_{m})
\theta_{2}(\xi_{m})\theta_{3}(\xi_{m})}{
\theta_{3}^{2}(\xi_{m})\theta_{1}^{2}(\xi_{m})+
\theta_{0}^{2}(\xi_{m})\theta_{2}^{2}(\xi_{m})}\\
&\\
&=\dfrac{\theta_{1}(2\xi_{m})\theta_{3}(0)}{
\theta_{0}(2\xi_{m})\theta_{2}(0)}
={\rm sn}(4K\xi_{m}) .\\
\end{array}
\end{eqnarray*}
Therefore, we have $\hat{w}_{m}=2\arccos({\rm cn}(4K\xi_{m}))$. 
Then, by the similar calculation as in \cite{udagawa}
we see that $\hat{w}_{m}=2\arccos({\rm cn}(4K\xi_{m}))$
satisfies the semi-discrete sine-Gordon equation
for $\widetilde{\alpha}$ replaced by
$\widetilde{\gamma}=-8k^{2}KA\dfrac{{\rm sn}(2K\Omega)
{\rm cn}(2K\Omega)}{{\rm dn}(2K\Omega)}$.
It also satisfies the mKdV equation for $\widetilde{\beta}$ replaced by $\widetilde{\delta}
=8k^{2}KA\dfrac{{\rm dn}(2K\Omega)}{{\rm sn}(2K\Omega)
{\rm cn}(2K\Omega)}$.

\section{Solutions of discrete sine-Gordon equation in terms of the Jacobi elliptic functions}

Set $\xi_{m,n}=m\Omega+n{\it P}+\xi_{0}$.
In \cite{udagawa}, it is proved that $\hat{w}_{m,n}=2\arcsin(
k{\rm sn}(4K\xi_{m,n}))$ is a solution of the discrete sine-Gordon equation
\begin{eqnarray}\label{sineGordon}
\begin{array}{ll}
&\sin\left(\frac{1}{4}(\hat{w}_{m+1,n+1}-\hat{w}_{m+1,n})-
\frac{1}{4}(\hat{w}_{m,n+1}-\hat{w}_{m,n})\right)\\
=&\hat{\gamma}
\sin\left(\frac{1}{4}(\hat{w}_{m+1,n+1}+\hat{w}_{m+1,n})+
\frac{1}{4}(\hat{w}_{m,n+1}+\hat{w}_{m,n})\right),\\
\end{array}
\end{eqnarray}
where
\begin{equation*}
\hat{\gamma}=-\dfrac{{\rm sn}(2K\Omega){\rm dn}(2K\Omega)}
{{\rm cn}(2K\Omega)}\cdot
\dfrac{{\rm sn}(2K{\it P}){\rm dn}(2K{\it P})}
{{\rm cn}(2K{\it P})}
\end{equation*}
for a solution $\hat{w}_{m,n}=2\arcsin(
k{\rm sn}(4K\xi_{m,n}))$. 
Since we may write $\hat{w}_{m,n}=2\arccos({\rm dn}
(4K\xi_{m,n}))$, this solution is called \lq\lq
${\rm dn}$-solution. Next, we search for the
\lq\lq ${\rm cn}$-solution. For this,
set
\begin{equation*}
f_{m,n}=\theta_{3}(\xi_{m,n})\theta_{1}(\xi_{m,n}) ,\quad
g_{m,n}=\theta_{0}(\xi_{m,n})\theta_{2}(\xi_{m,n}).
\end{equation*}
We consider $\tan(\frac{\hat{w}_{m,n}}{4})=\dfrac{g_{m,n}}{f_{m,n}}$. 
Setting $A=\hat{w}_{m+1,n+1}, B=\hat{w}_{m,n},
C=\hat{w}_{m+1,n}, D=\hat{w}_{m,n+1}$, we calculate
$\tan(\frac{A+B}{4})$ using (i) of Lemma~\ref{App1}, 
we obtain the following.
\begin{eqnarray*}
\begin{array}{ll}
&\tan(\frac{A+B}{4})=\dfrac{f_{m,n}g_{m+1,n+1}+f_{m+1,n+1}
g_{m,n}}{f_{m+1,n+1}f_{m,n}-g_{m+1,n+1}g_{m,n}} \\
&=-\dfrac{\theta_{1}(2\xi_{m,n}+\Omega+{\it P})
\theta_{3}(\Omega+{\it P})}{
\theta_{2}(2\xi_{m,n}+\Omega+{\it P})
\theta_{0}(\Omega+{\it P})}
=-{\rm dn}(2K(\Omega+{\it P}))\cdot
\dfrac{{\rm sn}(4K\xi_{m,n}+2K(\Omega+{\it P}))}{
{\rm cn}(4K\xi_{m,n}+2K(\Omega+{\it P}))},
\end{array}
\end{eqnarray*}
and
\begin{eqnarray*}
\begin{array}{ll}
&\tan(\frac{C+D}{4})=\dfrac{f_{m,n+1}g_{m+1,n}+f_{m+1,n}
g_{m,n+1}}{f_{m+1,n}f_{m,n+1}-g_{m+1,n}g_{m,n+1}} \\
&=-\dfrac{\theta_{1}(2\xi_{m,n}+\Omega+{\it P})
\theta_{3}(\Omega-{\it P})}{
\theta_{2}(2\xi_{m,n}+\Omega+{\it P})
\theta_{0}(\Omega-{\it P})}
=-{\rm dn}(2K(\Omega-{\it P}))\cdot
\dfrac{{\rm sn}(4K\xi_{m,n}+2K(\Omega+{\it P}))}{
{\rm cn}(4K\xi_{m,n}+2K(\Omega+{\it P}))}.
\end{array}
\end{eqnarray*}
Now, we see that
\begin{eqnarray*}
\begin{array}{ll}
\hat{\gamma}&=\dfrac{\sin\left(\frac{A+B}{4}-\frac{C+D}{4}
\right)}{
\sin\left(\frac{A+B}{4}+\frac{C+D}{4}
\right)}
=\dfrac{{\rm dn}(2K(\Omega+{\it P}))-
{\rm dn}(2K(\Omega-{\it P}))}{
{\rm dn}(2K(\Omega+{\it P}))+
{\rm dn}(2K(\Omega-{\it P}))}
\\
&=-k^{2}\dfrac{{\rm sn}(2K\Omega){\rm cn}(2K\Omega)}{
{\rm dn}(2K\Omega)}\cdot
\dfrac{{\rm sn}(2K{\it P}){\rm cn}(2K{\it P})}{
{\rm dn}(2K{\it P})}.
\end{array}
\end{eqnarray*}
Since we have
\begin{equation*}
\tan(\frac{\hat{w}_{m,n}}{4})=
\dfrac{\theta_{0}(\xi_{m,n})\theta_{2}(\xi_{m,n})}{
\theta_{3}(\xi_{m,n})\theta_{1}(\xi_{m,n})}
=\dfrac{{\rm cn}(2K\xi_{m,n})}{
{\rm dn}(2K\xi_{m,n}){\rm sn}(2K\xi_{m,n})},
\end{equation*}
we obtain the following.
\begin{eqnarray*}
\begin{array}{ll}
\sin(\frac{\hat{w}_{m,n}}{2})&=2\sin(\frac{\hat{w}_{m,n}}{4})\cos(\frac{\hat{w}_{m,n}}{4})=
\dfrac{2\tan(\frac{\hat{w}_{m,n}}{4})}{1+\tan^{2}(
\frac{\hat{w}_{m,n}}{4})}\\
&=\dfrac{2{\rm sn}(2K\xi_{m,n}){\rm cn}(2K\xi_{m,n})
{\rm dn}(2K\xi_{m,n})}{1-k^{2}{\rm sn}^{4}(2K\xi_{m,n})}
={\rm sn}(4K\xi_{m,n}) .\\
\end{array}
\end{eqnarray*}
Therefore, we see that $\cos(\frac{\hat{w}_{m,n}}{2})=
{\rm cn}(4K\xi_{m,n})$, which is a \lq\lq ${\rm cn}$-
solution\rq\rq.

\section{Construction of semi-discrete surface}

Given a solution $\hat{w}_{m}$ of the semi-discrete
sine-Gordon equation,we write it as
$\exp(i\dfrac{\hat{w}_{m}}{2})
=\left.-i\dfrac{g_{m}}{f_{m}^{*}}
\right|_{\lambda=\lambda_{0},z=i\lambda_{0}}$.
After choosing such $f_{m}, g_{m}$,
we consider the unitary frame $\Phi_{m}$ defined by
\begin{equation*}
\Phi_{m}=\dfrac{1}{\sqrt{F_{m}}}
\left(\begin{matrix}\alpha_{m}&0\\ 0&\alpha^{-1}_{m}\\
\end{matrix}\right)
\left(\begin{matrix} f_{m}&-g^{*}_{m}\\ g_{m}&f^{*}_{m}\\
\end{matrix}\right) ,\quad F_{m}:=f_{m}f^{*}_{m}+g_{m}g^{*}_{m} ,
\end{equation*}
where $(\alpha_{m})^{*}=\alpha^{-1}_{m}$
and 
all the functions in the entries are functions of some parameters $\lambda, z$ and $t$. 
Setting $U_{m}=\left.\dfrac{\partial \Phi_{m}}{\partial \lambda}
\cdot\Phi^{-1}_{m}\right|_{\lambda=\lambda_{0}}$,
we write the component of $U_{m}$ as $U_{m}=
\left(\begin{matrix} U_{11}& U_{12}\\ U_{21}& U_{22}\end{matrix}
\right)$. 
We use an isomorphism $\varphi : {\mathbb R}^{3} \longrightarrow
{\frak s}{\frak u}(2)$ defined by 
\begin{equation*}
\varphi\left(\begin{matrix}x_{1}\\x_{2}\\x_{3}\\\end{matrix}
\right)=x_{1}E_{1}+x_{2}E_{2}+x_{3}E_{3}=
\left(\begin{matrix}-ix_{3}&-ix_{1}-x_{2}\\
-ix_{1}+x_{2}&ix_{3}\\\end{matrix}\right),
\end{equation*}
where $E_{1}=\left(\begin{matrix}0&-i\\ -i&0\\\end{matrix}\right),
E_{2}=\left(\begin{matrix}0&-1\\ 1&0\\\end{matrix}\right),
E_{3}=\left(\begin{matrix}-i&0\\ 0&i\\\end{matrix}\right)$.
We choose some $\lambda_{0}$ so that 
$\left.\alpha_{m}\right|_{\lambda=\lambda_{0}}=1$.
Set 
\begin{equation*}
R_{m}=\left.\alpha_{m}^{-1}
\dfrac{\partial\alpha_{m}}{
\partial \lambda}\right|_{\lambda=\lambda_{0}} .
\end{equation*}
Hereafter, we write 
$\dfrac{\partial f_{m}}{
\partial\lambda}$ as $(f_{m})_{\lambda}$ and so on.

We then express a vector $\Gamma_{m}=\left(\begin{matrix}
\Gamma_{m}^{1}\\ \Gamma^{2}_{m}\\ \Gamma^{3}_{m}\\
\end{matrix}\right)\in\mathbb{R}^{3}$ as
\begin{equation}
\Gamma_{m}=\left(\begin{matrix}
\frac{i}{2}\left(U_{21}-U^{*}_{21}\right)\\
\frac{1}{2}\left(U_{21}+U^{*}_{21}\right)\\
iU_{11}\\ \end{matrix}\right) ,
\end{equation}
where
\begin{eqnarray}
\hskip 1cm\left\{
\begin{array}{ll}
U_{21}&=\left.\dfrac{1}{F_{m}}\left(f^{*}_{m}(g_{m})_{\lambda} - (f^{*}_{m})_{\lambda}g_{m}\right)
\right|_{\lambda=\lambda_{0}},\\
U_{11}&=\left.\dfrac{1}{2F_{m}}\left(f^{*}_{m}(f_{m})_{\lambda}
-f_{m}(f^{*}_{m})_{\lambda}+g_{m}(g^{*}_{m})_{\lambda}
-g^{*}_{m}(g_{m})_{\lambda}\right)\right|_{\lambda=\lambda_{0}}
+R_{m} .\\
\end{array}
\right.
\end{eqnarray}

For a nonnegative integer $m$, we introduce functions $F_{m}, H_{m}, 
\Psi^{FH}_{m}, \Psi^{FR}_{m}$ satisfying following relations.
\begin{eqnarray}
\hskip 1cm \left\{
\begin{array}{ll}\label{fgFH}
&F_{m}=\left.\left(f_{m}f^{*}_{m}+g_{m}g^{*}_{m}\right)\right|_{\lambda=\lambda_{0}} ,\\
&\\
&H_{m}=\left.\dfrac{i}{2}D_{\lambda}g_{m}\cdot f^{*}_{m}\right|_{\lambda=\lambda_{0}},\\
&F_{m}H_{m+1}-H_{m}F_{m+1}=\left.
\dfrac{\varepsilon}{i}\Psi^{FH}_{m}\right|_{\lambda=\lambda_{0}},\\
&\Psi^{FH}_{m}=f^{*}_{m}f^{*}_{m+1}\left(f_{m}g_{m+1}-f_{m+1}g_{m}\right)+g_{m}g_{m+1}\left(
f^{*}_{m}g^{*}_{m+1}-f^{*}_{m+1}g^{*}_{m}\right) ,\\
&\\
&\left.\dfrac{1}{2}D_{z}F_{m}\cdot F_{m+1}\right|_{z=0}
+i(R_{m+1}-R_{m})F_{m+1}F_{m}
=\left.\dfrac{2\varepsilon}{i}\psi^{FR}_{m}
\right|_{\lambda=\lambda_{0}},\\
&\quad \Psi^{FR}_{m}=
f_{m+1}f^{*}_{m}g_{m}g^{*}_{m+1}-
f^{*}_{m+1}f_{m}g^{*}_{m}g_{m+1} ,
\end{array}
\right.
\end{eqnarray}
where $i$ the imaginary unit, $\varepsilon=\pm 1$ and $D_{z}$
and $D_{\lambda}$ is the Hirota differetiation operator with a real parameter $z$ and $\lambda$, respectively.
Moreover, we assume that the following relations hold.
\begin{eqnarray}
\left\{
\begin{array}{ll}
(f_{m})_{\lambda}&=i(f_{m})_{z},\\
(f^{*}_{m})_{\lambda}&=-i(f^{*}_{m})_{z} ,
\end{array}
\right.
\left\{
\begin{array}{ll}
(g_{m})_{\lambda}&=-i(g_{m})_{z} ,\\
(g^{*}_{m})_{\lambda}&=i(g^{*}_{m})_{z} .\\
\end{array}
\right.
\end{eqnarray}

With these solutions we obtain semi-discrete surface  
$\Gamma_{m}$. We may define a Frenet frame 
$\widetilde{\Phi}_{m}=(T_{m}\hskip .1cm 
N_{m}\hskip .1cm B_{m})\in SO(3)$ associated to $\Gamma_{m}$. 
Set $\lambda=\lambda_{0}$ and $z=0$. 
When ${\cal K}_{m+1}=\frac{1}{2}\hat{w}_{m+2}-\frac{1}{2}\hat{w}_{m}$, we define the frame as follows.
\begin{eqnarray}\label{frenetframe1}
\left\{
\begin{array}{ll}
\varphi(T_{m})&=\pm\Phi_{m}\left(-\sin(\frac{\hat{w}_{m+1}}{2})E_{1}
+\cos(\frac{\hat{w}_{m+1}}{2})E_{2}\right)\Phi_{m}^{-1} ,\\
\varphi(N_{m})&=\pm\Phi_{m}\left(-\cos(\frac{\hat{w}_{m+1}}{2})E_{1}
-\sin(\frac{\hat{w}_{m+1}}{2})E_{2}\right)\Phi_{m}^{-1} ,\\
\varphi(B_{m})&=\Phi_{m}E_{3}\Phi_{m}^{-1}.\\
\end{array}
\right.
\end{eqnarray}
We write $\Phi_{m+1}=\Phi_{m}L_{m}$, where $L_{m}$ corresponds
to $\widetilde{L}_{m}$ in~(\ref{trans}). We find $L_{m}$ as follows.

\begin{Lemma}\label{frame}
For ${\cal K}_{m+1}=\frac{1}{2}\hat{w}_{m+2}-\frac{1}{2}\hat{w}_{m}$,
we have the following.
\item{\rm (1)} 
\begin{equation*}
L^{\pm}_{m}=
\left(\begin{matrix}\cos(\frac{\nu}{2})e^{-\frac{i}{4}(\hat{w}_{m+1}
-\hat{w}_{m})}&
\pm\sin(\frac{\nu}{2})e^{-\frac{i}{4}(\hat{w}_{m+1}+\hat{w}_{m})}\\
\mp\sin(\frac{\nu}{2})e^{\frac{i}{4}(\hat{w}_{m+1}+\hat{w}_{m})}&
\cos(\frac{\nu}{2})e^{\frac{i}{4}(\hat{w}_{m+1}-\hat{w}_{m})}\\
\end{matrix}\right),
\end{equation*}
where the compound sign is in the same order as the compound sign 
in~{\rm (\ref{frenetframe1})}.
\item{\rm (2)} $B_{m+1}\times B_{m}=\sin\nu T_{m}$.
\end{Lemma}

{\it Proof}. (1) We show 
$\left(\varphi(T_{m+1})\hskip .1cm \varphi(N_{m+1})
\hskip .1cm \varphi(B_{m+1})\right)=
\left(\varphi(T_{m})\hskip .1cm \varphi(N_{m})
\hskip .1cm \varphi(B_{m})\right)\widetilde{L}_{m}$. 
For this purpose, we calculate $L_{m}E_{1}L_{m}^{-1}, L_{m}E_{2}L_{m}^{-1}$
and $L_{m}E_{3}L_{m}^{-1}$ and obtain the following.
\begin{eqnarray*}
\begin{array}{ll}
L^{\pm}_{m}E_{1}(L^{\pm}_{m})^{-1}&=\left(\cos\nu\cos(\frac{\hat{w}_{m}}{2})\cos(\frac{\hat{w}_{m+1}}{2})+\sin(\frac{\hat{w}_{m}}{2})
\sin(\frac{\hat{w}_{m+1}}{2})\right) E_{1} \\
&+\left(\cos\nu\cos(\frac{\hat{w}_{m}}{2})\sin(\frac{\hat{w}_{m+1}}{2})
-\sin(\frac{\hat{w}_{m}}{2})\cos(\frac{\hat{w}_{m+1}}{2})\right) E_{2}
\pm\sin\nu\cos(\frac{\hat{w}_{m}}{2}) E_{3} ,\\
&\\
L^{\pm}_{m}E_{2}(L^{\pm}_{m})^{-1}
&=\left(\cos\nu\sin(\frac{\hat{w}_{m}}{2})\cos(\frac{\hat{w}_{m+1}}{2})-\cos(\frac{\hat{w}_{m}}{2})
\sin(\frac{\hat{w}_{m+1}}{2})\right) E_{1} \\
&+\left(\cos\nu\sin(\frac{\hat{w}_{m}}{2})\sin(\frac{\hat{w}_{m+1}}{2})
+\cos(\frac{\hat{w}_{m}}{2})\cos(\frac{\hat{w}_{m+1}}{2})\right) E_{2}
\pm\sin\nu\sin(\frac{\hat{w}_{m}}{2}) E_{3},\\
&\\
L^{\pm}_{m}E_{3}(L^{\pm}_{m})^{-1}
&=\mp\sin\nu\cos(\frac{\hat{w}_{m+1}}{2}) E_{1}
\mp\sin\nu\sin(\frac{\hat{w}_{m+1}}{2}) E_{2} +\cos\nu E_{3} .\\
\end{array}
\end{eqnarray*}
Using ${\cal K}_{m+1}=\frac{1}{2}\hat{w}_{m+2}-\frac{1}{2}\hat{w}_{m}$,
we then have the following.
\begin{eqnarray*}
\begin{array}{ll}
\varphi(B_{m+1})&=\Phi_{m+1}E_{3}\Phi_{m+1}^{-1}
=\Phi_{m}(L_{m}E_{3}L_{m}^{-1})\Phi_{m}^{-1} \\
&=\Phi_{m}\left(\pm\sin\nu\left(-\cos(\frac{\hat{w}_{m+1}}{2})E_{1}
-\sin(\frac{\hat{w}_{m+1}}{2})E_{2}\right) 
+\cos\nu E_{3}\right)\Phi_{m}^{-1}\\
&=\varphi(\sin\nu N_{m}+\cos\nu B_{m}) ,\\
&\\
\varphi(T_{m+1})&=\pm\Phi_{m+1}\left(-\sin(\frac{\hat{w}_{m+2}}{2})E_{1}
+\cos(\frac{\hat{w}_{m+2}}{2})E_{2}\right)\Phi_{m+1}^{-1} \\
&=\pm\Phi_{m}\left(-\sin(\frac{\hat{w}_{m+2}}{2})L_{m}E_{1}L_{m}^{-1}
+\cos(\frac{\hat{w}_{m+2}}{2})L_{m}E_{2}L_{m}^{-1}\right)\Phi_{m}^{-1}\\
&=\varphi(\cos\nu\sin{\cal K}_{m+1} N_{m} +\cos{\cal K}_{m+1} T_{m}
-\sin\nu\sin{\cal K}_{m+1} B_{m}) ,\\
&\\
\varphi(N_{m+1})&=\pm\Phi_{m+1}\left(-\cos(\frac{\hat{w}_{m+2}}{2})E_{1}
-\sin(\frac{\hat{w}_{m+2}}{2})E_{2}\right)\Phi_{m+1}^{-1} \\
&=\pm\Phi_{m}\left(-\cos(\frac{\hat{w}_{m+2}}{2})L_{m}E_{1}L_{m}^{-1}
-\sin(\frac{\hat{w}_{m+2}}{2})L_{m}E_{2}L_{m}^{-1}\right)\Phi_{m}^{-1} \\
&=\varphi(-\sin{\cal K}_{m+1}T_{m}+\cos\nu\cos{\cal K}_{m+1}N_{m}
-\sin\nu\cos{\cal K}_{m+1}B_{m}) .\\
\end{array}
\end{eqnarray*}
(2) First note that $\varphi(A\times B)=\frac{1}{2}\bigl[\varphi(A),
\varphi(B)\bigr]$, where $A, B\in\mathbb{R}^{3}$. We then have
\begin{eqnarray*}
\begin{array}{ll}
\varphi(B_{m+1}\times B_{m})&=\frac{1}{2}\bigl[\Phi_{m+1}E_{3}
\Phi_{m+1}^{-1}, \Phi_{m}E_{3}\Phi_{m}^{-1}\bigr]
=\frac{1}{2}\Phi_{m}\bigl[L_{m}E_{3}L_{m}^{-1}, E_{3}\bigr]
\Phi_{m}^{-1} \\
&=\frac{1}{2}\Phi_{m}\left(
\bigl[\mp\sin\nu\cos(\frac{\hat{w}_{m+1}}{2})E_{1}
\mp\sin\nu\sin(\frac{\hat{w}_{m+1}}{2})E_{2}+\cos\nu E_{3}, E_{3}
\bigr]\right)\Phi_{m}^{-1}\\
&=\Phi_{m}\left(\pm\sin\nu\cos(\frac{\hat{w}_{m+1}}{2})E_{2}
\mp\sin\nu\sin(\frac{\hat{w}_{m+1}}{2})E_{1}\right)\Phi_{m}^{-1}
=\sin\nu \varphi(T_{m}) .\\
\end{array}
\end{eqnarray*}
Therefore, we have $B_{m+1}\times B_{m}=\sin\nu T_{m}$.
\qed

Alternatively, when ${\cal K}_{m+1}=-\frac{1}{2}\hat{w}_{m+2}+
\frac{1}{2}\hat{w}_{m}$, we define the frame as follows.
\begin{eqnarray}\label{frenetframe2}
\left\{
\begin{array}{ll}
\varphi(T_{m})&=\pm\hat{\Phi}_{m}\left(\sin(\frac{\hat{w}_{m+1}}{2})E_{1}
+\cos(\frac{\hat{w}_{m+1}}{2})E_{2}\right)\hat{\Phi}_{m}^{-1} ,\\
\varphi(N_{m})&=\pm\hat{\Phi}_{m}\left(-\cos(\frac{\hat{w}_{m+1}}{2})E_{1}
+\sin(\frac{\hat{w}_{m+1}}{2})E_{2}\right)\hat{\Phi}_{m}^{-1} ,\\
\varphi(B_{m})&=\hat{\Phi}_{m}E_{3}\hat{\Phi}_{m}^{-1}.\\
\end{array}
\right.
\end{eqnarray}
We write $\hat{\Phi}_{m+1}=\hat{\Phi}_{m}\hat{L}_{m}$. We then have the
following.
\begin{Lemma}\label{frame2}
For ${\cal K}_{m+1}=-\frac{1}{2}\hat{w}_{m+2}+\frac{1}{2}\hat{w}_{m}$,
we have the following.
\item{\rm (1)} 
\begin{equation*}
\hat{L}^{\pm}_{m}=\left(\begin{matrix}\cos(\frac{\nu}{2})e^{\frac{i}{4}(\hat{w}_{m+1}+\hat{w}_{m})}&
\pm\sin(\frac{\nu}{2})e^{\frac{i}{4}(\hat{w}_{m+1}-\hat{w}_{m})}\\
\mp\sin(\frac{\nu}{2})e^{-\frac{i}{4}(\hat{w}_{m+1}-\hat{w}_{m})}&
\cos(\frac{\nu}{2})e^{-\frac{i}{4}(\hat{w}_{m+1}+\hat{w}_{m})}\\
\end{matrix}\right),
\end{equation*}
where the compound sign is in the same order as the compound sign 
in~{\rm (\ref{frenetframe2})}.

\item{\rm (2)} $B_{m+1}\times B_{m}=\sin\nu T_{m}$.
\end{Lemma}

\begin{Theorem}\label{fundamentaltheorem}
After choosing $f_{m}, g_{m}$
which satisfies $\exp(i\dfrac{\hat{w}_{m}}{2})
=\left.-i\dfrac{g_{m}}{f^{*}_{m}}
\right|_{\lambda=\lambda_{0},z=i\lambda_{0}}$
for a solution of the semi-discrete sine-Gordon equation,
under the conditions {\rm (3.1)}$\sim${\rm (3.4)}, 
we obtain the following.
\begin{eqnarray}\label{gammabm}
\hskip 1cm
\Gamma_{m}
&=\Biggl. \left(\begin{matrix}
\frac{H_{m}+H^{*}_{m}}{F_{m}} \\
\\
\frac{1}{i}\frac{H_{m}-H^{*}_{m}}{F_{m}}\\
\\
iR_{m}-\dfrac{1}{2}\frac{d}{dz}\log F_{m}\\
\end{matrix}\right)\Biggr{|}_{z=0} ,\quad
B_{m}=\Biggl. \dfrac{1}{F_{m}}\left(\begin{matrix}
f^{*}_{m}g_{m}+f_{m}g^{*}_{m}\\
\\
\frac{1}{i}\left(f^{*}_{m}g_{m}-f_{m}g^{*}_{m}\right)\\
\\
f_{m}f^{*}_{m}-g_{m}g^{*}_{m}\\
\end{matrix}\right)\Biggr{|}_{z=0} .
\end{eqnarray}
Moreover, we have  
$\Gamma_{m+1}-\Gamma_{m}=\varepsilon B_{m+1}\times B_{m}$.
\end{Theorem}

{\it Proof}. The form of $\Gamma_{m}$ 
follows from (3.1), (3.2), first two equations 
in (3.3) and (3.4). The form of $B_{m}$ 
easily follows from its definition (3.5)(or (3.7)).
We calculate $B_{m+1}\times B_{m}$, then
we find the following using the last four 
equations of (~\ref{fgFH}).
\begin{eqnarray*}
\begin{array}{ll}
\varepsilon B_{m+1}\times B_{m}
&=\dfrac{\varepsilon}{F_{m+1}F_{m}}
\left(\begin{matrix}
\frac{1}{i}\left(\Psi^{FH}_{m}-
\left(\Psi^{FH}_{m}\right)^{*}\right)\\
-\left(\Psi^{FH}_{m}+
\left(\Psi^{FH}_{m}\right)^{*}\right)\\
\frac{2}{i}\Psi^{FR}_{m}\\
\end{matrix}\right)\\
&\\
&=\left(\begin{matrix}
\frac{H_{m+1}+H_{m+1}^{*}}{F_{m+1}}
-\frac{H_{m}+H_{m}^{*}}{F_{m}}\\
\frac{1}{i}\frac{H_{m+1}-H_{m+1}^{*}}{
F_{m+1}}-\frac{1}{i}\frac{H_{m}
-H_{m}^{*}}{F_{m}}\\
\left.iR_{m+1}-\frac{1}{2}\frac{d}{dz}\right|_{z=0}
\log F_{m+1}-\left.iR_{m}+
\frac{1}{2}\frac{d}{dz}\right|_{z=0}
\log F_{m}\\
\end{matrix}\right)
=\Gamma_{m+1}-\Gamma_{m} .\\
\end{array}
\end{eqnarray*}

It follows from (2) of Lemma~\ref{frame} that
we may define $T_{m}$ and $N_{m}$ by
\begin{equation*}
T_{m}:=\dfrac{\Gamma_{m+1}-\Gamma_{m}}{
|\Gamma_{m+1}-\Gamma_{m}|}, \quad
N_{m}:=B_{m}\times T_{m} .
\end{equation*}

\qed

We give here the example of $\alpha_{m}$.
We consider the following Weierstrass ${\frak p}$-function
stated in~(\ref{WP}) :
\begin{equation*}
{\frak p}(w)=\left({\rm dn}(2iw+iK^{\prime})-ik{\rm sn}(2iw+
iK^{\prime})\right)^{2}+e_{1},
\end{equation*}
where $e_{1}=\dfrac{2}{3}(2k^{2}-1)$ and
$K^{\prime}=K(k^{\prime})$. The two periods of ${\frak p}(w)$
are given by $2\omega=2K^{\prime}$ and $2\omega^{\prime}=
iK+K^{\prime}$. In particular, we see that ${\frak p}(\frac{1}{2}K^{\prime})=1+\dfrac{2}{3}(2k^{2}-1)$. Let $\zeta_{\it W}$ be the
Weierstrass $\zeta$-function. We define $\alpha_{m}$ by
\begin{equation}\label{alpha}
\alpha_{m}
=\exp\left(\frac{i}{2}\eta_{m}
\displaystyle{\int_{w(\lambda_{0})}^{w(\lambda)}
\left({\frak p}(w)+\dfrac{\zeta_{\it W}(\omega)}{\omega}\right)dw}\right),
\end{equation}
where $w(\lambda)$ is a function of $\lambda$ and
$\eta_{m}$ is some real-valued function which depends only $m$ and $t$,
which is determined explicitly later.
Moreover, we choose $\lambda$ so that
the integration part in the $\alpha_{m}$
is real-valued.

\section{Semi-discrete surfaces corresponding to 
$\lq\lq{\rm dn}$-solution\rq\rq of semi-discrete SG-equation}

We construct a semi-discrete space curve which corresponds
to the solution $\hat{w}_{m}=2\arccos(\text{dn}(\psi_{m}))$
of the semi-discrete SG-equation.
Set
\begin{eqnarray*}
v_{m}&=\dfrac{\psi_{m}-K}{2iK^{\prime}}, \quad
v^{+}_{m}(\lambda,z)=v_{m}+\dfrac{\lambda}{kK^{\prime}}+
\dfrac{i}{kK^{\prime}}z,\quad 
v^{-}_{m}(\lambda,z)=v_{m}-\dfrac{\lambda}{kK^{\prime}}+
\dfrac{i}{kK^{\prime}}z ,
\end{eqnarray*}
where $\psi_{m}$ is a real-valued function of some parameters,
$K=K(k)$ is the complete elliptic integral of 1st kind with modulus $k$ and $K^{\prime}=K(k^{\prime})$ for 
$k^{\prime}=\sqrt{1-k^{2}}$.
We then see that $\left(v^{+}_{m}(\lambda,z)\right)^{*}
=-v^{-}_{m}(\lambda,z)$ 
for real parameters $\lambda$ and $z$.  
In this case, we may verify the following.
\begin{eqnarray}
\left\{
\begin{array}{ll}
\left(\theta_{3}(v^{+}_{m}(\lambda,z) |2\tau^{\prime})\right)^{*}
&=\theta_{3}(-v^{-}_{m}(\lambda,z) | 2\tau^{\prime})=
\theta_{3}(v^{-}_{m}(\lambda,z) |2\tau^{\prime}),\\
&\\
\left(\theta_{2}(v^{+}_{m}(\lambda,z) |2\tau^{\prime})\right)^{*}
&=\theta_{2}(-v^{-}_{m}(\lambda,z) | 2\tau^{\prime})=
\theta_{2}(v^{-}_{m}(\lambda,z) |2\tau^{\prime}),\\
\end{array}
\right.
\end{eqnarray}
where $\tau=i\dfrac{K^{\prime}}{K}$ and $\tau^{\prime}=i\dfrac{K}{K^{\prime}}$.

We choose $\lambda_{0}=\dfrac{kK^{\prime}}{2}$ in this section.
Set 
\begin{eqnarray*}
v^{+}_{m}(z)&=v_{m}+\dfrac{1}{2}+\dfrac{i}{kK^{\prime}}z,\quad
v^{-}_{m}(z)=v_{m}-\dfrac{1}{2}+\dfrac{i}{kK^{\prime}}z,\quad
v_{m}(z)=v_{m}+\dfrac{i}{kK^{\prime}}z.
\end{eqnarray*}

\begin{Lemma} Set $\varphi_{m}=m\alpha+\beta kt$
and $\psi_{m}=m\gamma+\beta t$.
We choose $\tau$-functions as follows.
\begin{eqnarray*}
\left\{
\begin{array}{ll}
f_{m}&=\exp(-\frac{i}{2}\varphi_{m})
\left(\theta_{3}(v^{-}_{m}(\lambda,z) | 2\tau^{\prime})+i 
\theta_{2}(v^{-}_{m}(\lambda,z) | 2\tau^{\prime})\right) ,\\
&\\
g_{m}&=\exp(\frac{i}{2}\varphi_{m})
\left(\theta_{3}(v^{+}_{m}(\lambda,z) | 2\tau^{\prime})+i 
\theta_{2}(v^{+}_{m}(\lambda,z) | 2\tau^{\prime})\right) ,\\
&\\
F_{m}&=\left.\left(f_{m}f^{*}_{m}+g_{m}g^{*}_{m}
\right)\right|_{\lambda=\frac{kK^{\prime}}{2}}
=2\theta_{3}(v_{m}(z) |\tau^{\prime})
\theta_{0}(0 | \tau^{\prime}) ,\\
&\\
H_{m}&=\frac{1}{kK^{\prime}}\exp(i\varphi_{m})\left(
\theta_{2}(v^{+}_{m}(z) |2\tau^{\prime}) \theta^{\prime}_{3}
(v^{+}_{m}(z) |2\tau^{\prime}) -
\theta_{3}(v^{+}_{m}(z) |2\tau^{\prime}) \theta^{\prime}_{2}
(v^{+}_{m}(z) |2\tau^{\prime})\right),\\
&\\
\eta_{m}&=\displaystyle{
\psi_{m}-\frac{\pi}{2E^{\prime}}-\frac{mk^{2}K^{\prime}}{E^{\prime}}
\int_{0}^{\gamma}{\rm sn}^{2}\psi d\psi}.\\
\end{array}
\right.
\end{eqnarray*}
We then have $\exp(i\dfrac{\hat{w}_{m}}{2})
=\left.-i\dfrac{g_{m}}{f_{m}^{*}}
\right|_{\lambda=\frac{kK^{\prime}}{2},
z=i\frac{kK^{\prime}}{2}}
={\rm dn}(\psi_{m})-ik{\rm sn}(\psi_{m})$, hence
$\hat{w}_{m}$ is a \lq\lq{\rm dn}-solution\rq\rq. 
We may verify that these functions satisfy 
{\rm (3.3)} and {\rm (3.4)} 
for $\varepsilon=+1$. 
\end{Lemma}

{\it Proof}. Let $\hat{w}_{m}$ be a \lq\lq{\rm dn}-
solution\rq\rq.   
For $\xi_{m}=m\Omega+\xi_{0}+At$, we take
$\Omega=\dfrac{\gamma}{4K}, \xi_{0}=\dfrac{1}{2}$
and $A=\dfrac{\beta}{4K}$,which yield
$4K\xi_{m}=\psi_{m}+2K$. Since $\sin\left(\dfrac{\hat{w}_{m}}{2}
\right)=k{\rm sn}(4K\xi_{m})=-k{\rm sn}(\psi_{m})$
and $\cos\left(\dfrac{\hat{w}_{m}}{2}\right)
={\rm dn}(4K\xi_{m})={\rm dn}(\psi_{m})$,
we have $\exp (i\dfrac{\hat{w}_{m}}{2})=
{\rm dn}(\psi_{m})-ik{\rm sn}(\psi_{m})$.
On the other hand, 
\begin{eqnarray*}
\begin{array}{ll}
\left.-i\dfrac{g_{m}}{f_{m}^{*}}
\right|_{\lambda=\frac{kK^{\prime}}{2},
z=i\frac{kK^{\prime}}{2}}&=
-i\dfrac{\theta_{3}(v_{m}|2\tau^{\prime})+
i\theta_{2}(v_{m}|2\tau^{\prime})}{
\theta_{3}(v_{m}|2\tau^{\prime})-
i\theta_{2}(v_{m}|2\tau^{\prime})}\\
&=-i\dfrac{\theta_{3}^{2}(v_{m}|2\tau^{\prime})-
\theta_{2}^{2}(v_{m}|2\tau^{\prime})+2i
\theta_{2}(v_{m}|2\tau^{\prime})
\theta_{3}(v_{m}|2\tau^{\prime})}{
\theta_{3}^{2}(v_{m}|2\tau^{\prime})+
\theta_{2}^{2}(v_{m}|2\tau^{\prime})}\\
&=-i\dfrac{\theta_{0}(v_{m}|\tau^{\prime})
\theta_{0}(0|\tau^{\prime})+i
\theta_{2}(v_{m}|\tau^{\prime})
\theta_{2}(0|\tau^{\prime})}{
\theta_{3}(v_{m}|\tau^{\prime})
\theta_{3}(0|\tau^{\prime})}\\
&=-ik{\rm sn}(\psi_{m})+{\rm dn}(\psi_{m}),\\
\end{array}
\end{eqnarray*}
where we have used Lemma~\ref{App2}, Lemma~\ref{App3} and $k=\left(\dfrac{\theta_{0}(0 |\tau^{\prime})}{
\theta_{3}(0 |\tau^{\prime})}\right)^{2}$.
For $F_{m}$,we use the fourth equation of 
Lemma~\ref{App2}
for $\lambda=\frac{kK^{\prime}}{2}$.
\begin{eqnarray*}
\begin{array}{ll}
F_{m}&=(\theta_{3}(v^{-}_{m}(\lambda,z)|2\tau^{\prime})
+i\theta_{2}(v^{-}_{m}(\lambda,z)|2\tau^{\prime})
(\theta_{3}(v^{+}_{m}(\lambda,z)|2\tau^{\prime})
-i\theta_{2}(v^{+}_{m}(\lambda,z)|2\tau^{\prime})\\
&\quad \left.+(\theta_{3}(v^{+}_{m}(\lambda,z)|2\tau^{\prime})
+i\theta_{2}(v^{+}_{m}(\lambda,z)|2\tau^{\prime})
(\theta_{3}(v^{-}_{m}(\lambda,z)|2\tau^{\prime})
-i\theta_{2}(v^{-}_{m}(\lambda,z)|2\tau^{\prime})
\right|_{\lambda=\frac{kK^{\prime}}{2}}\\
&\\
&=\left.2\left(\theta_{3}(v^{+}_{m}(\lambda,z)|2\tau^{\prime})
\theta_{3}(v^{-}_{m}(\lambda,z)|2\tau^{\prime})+
\theta_{2}(v^{+}_{m}(\lambda,z)|2\tau^{\prime})
\theta_{2}(v^{-}_{m}(\lambda,z)|2\tau^{\prime})\right)
\right|_{\lambda=\frac{kK^{\prime}}{2}} \\
&=2\theta_{3}(v_{m}(z) |\tau^{\prime})
\theta_{3}(\frac{1}{2} |\tau^{\prime})
=2\theta_{3}(v_{m}(z) |\tau^{\prime})
\theta_{0}(0 |\tau^{\prime}).\\
\end{array}
\end{eqnarray*}
For $H_{m}$, considering (3.4) we calculate as follows.
\begin{eqnarray*}
\begin{array}{ll}
H_{m}&=\left.\dfrac{i}{2}D_{\lambda}g_{m}\cdot f^{*}_{m}
\right|_{\lambda=\frac{kK^{\prime}}{2}}\\
&=\dfrac{i}{2kK^{\prime}}\exp(i\varphi_{m})\left\{
\left(\theta_{3}^{\prime}(v^{+}_{m}(z)|2\tau^{\prime})
+i\theta_{2}^{\prime}(v^{+}_{m}(z)|2\tau^{\prime})
\right)\left(
\theta_{3}(v^{+}_{m}(z)|2\tau^{\prime})
-i\theta_{2}(v^{+}_{m}(z)|2\tau^{\prime})
\right)\right. \\
&\quad\left. -
\left(\theta_{3}(v^{+}_{m}(z)|2\tau^{\prime})
+i\theta_{2}(v^{+}_{m}(z)|2\tau^{\prime})
\right)\left(
\theta_{3}^{\prime}(v^{+}_{m}(z)|2\tau^{\prime})
-i\theta_{2}^{\prime}(v^{+}_{m}(z)|2\tau^{\prime})
\right)
\right\}\\
&=\dfrac{1}{kK^{\prime}}\exp(i\varphi_{m})
\left(\theta_{2}(v^{+}_{m}(z)|2\tau^{\prime})
\theta_{3}^{\prime}(v^{+}_{m}(z)|2\tau^{\prime})
-
\theta_{3}(v^{+}_{m}(z)|2\tau^{\prime})
\theta_{2}^{\prime}(v^{+}_{m}(z)|2\tau^{\prime})
\right) .\\
\end{array}
\end{eqnarray*}
The rest two equations in (4.3) are prepared
to ensure that $\Gamma_{m+1}-\Gamma_{m}=\varepsilon
B_{m+1}\times B_{m}$. Therefore, we prove them
in the proof of the next Theorem.
\qed

\begin{Theorem}\label{dnsolution}
\begin{eqnarray*}
\hskip 1cm
\Gamma_{m}&=\left(\begin{matrix} 
\frac{1}{k}\cos(\varphi_{m}) 
{\rm dn}(\psi_{m})\\
&\\
\frac{1}{k}\sin(\varphi_{m}) {\rm dn}(\psi_{m}) \\
&\\
\displaystyle{-k \int_{0}^{\psi_{m}}{\rm sn}^{2}\psi d\psi 
+km\int_{0}^{\gamma}{\rm sn}^{2}\psi d\psi}\\
\end{matrix}\right), 
\quad
B_{m}&=\left(\begin{matrix} \cos(\varphi_{m})
{\rm sn}(\psi_{m})\\
\sin(\varphi_{m}){\rm sn}(\psi_{m}) \\
-{\rm cn}(\psi_{m})\\
\end{matrix}\right),
\end{eqnarray*}
where $\varphi_{m}=m\alpha +\beta k t, \psi_{m}=
m\gamma +\beta t$ and
$\cos(\alpha)={\rm dn}(\gamma), \sin(\alpha)=
k {\rm sn}(\gamma)$. 
We have also the following properties.
\item{\rm (i)} $<B_{m}, B_{m+1}>={\rm cn}(\gamma)$, thus
$\cos\nu={\rm cn}(\gamma)$ and $\sin\nu=\pm{\rm sn}(\gamma)>0$.
The compound sign is in the same order as the compound sign 
in~{\rm (\ref{frenetframe1})}.
\item{\rm (ii)} $\cos(\frac{\hat{w}_{m}}{2})={\rm dn}(\psi_{m})$, 
$\sin(\frac{\hat{w}_{m}}{2})=-k {\rm sn}(\psi_{m})$ and 
${\cal K}_{m+1}=\frac{1}{2}\hat{w}_{m+2}-\frac{1}{2}\hat{w}_{m}$. 
\item{\rm (iii)} $\Gamma_{m}(t)$ is an isoperimetric
deformation of the discrete space curve, that is,
the following holds.
$\dfrac{d\Gamma_{m}}{dt}=\pm\beta\left(
\cos(w_{m}) T_{m}+\sin(w_{m}) N_{m}\right)$, where
$w_{m}=\frac{1}{2}\hat{w}_{m}-\frac{1}{2}\hat{w}_{m+1}$.
\end{Theorem}

{\it Proof}. First, we rewrite $H_{m}$.
\begin{eqnarray*}
\begin{array}{ll}
&\quad\theta_{2}(v^{+}_{m}(z) |2\tau^{\prime}) \theta^{\prime}_{3}
(v^{+}_{m}(z) |2\tau^{\prime}) -
\theta_{3}(v^{+}_{m}(z) |2\tau^{\prime}) \theta^{\prime}_{2}
(v^{+}_{m}(z) |2\tau^{\prime})\\
&=\theta_{2}^{2}(v^{+}_{m}|2\tau^{\prime})
\left(\dfrac{\theta_{3}(v^{+}_{m}(z)|2\tau^{\prime})}{
\theta_{2}(v^{+}_{m}(z)|2\tau^{\prime})}\right)^{\prime}\\
&=\dfrac{1}{2}\left(
\theta_{3}(v^{+}_{m}(z)|\tau^{\prime})
\theta_{3}(0 |\tau^{\prime})-
\theta_{0}(v^{+}_{m}(z)|\tau^{\prime})
\theta_{0}(0|\tau^{\prime})\right)
\left(\dfrac{\theta_{3}(v^{+}_{m}(z)|2\tau^{\prime})}{
\theta_{2}(v^{+}_{m}(z)|2\tau^{\prime})}
\right)^{\prime} ,\\
\end{array}
\end{eqnarray*}
where we have used the first and second equations of
Lemma~\ref{App2} for $\lambda=\frac{kK^{\prime}}{2}$.
We here calculate the reciprocal of the expression 
in parentheses.
\begin{eqnarray*}
\begin{array}{ll}
&\quad\dfrac{\theta_{2}(v^{+}_{m}(z)|2\tau^{\prime})}{
\theta_{3}(v^{+}_{m}(z)|2\tau^{\prime})}
=\dfrac{2\theta_{2}^{2}(v^{+}_{m}(z)|2\tau^{\prime})}{
2\theta_{2}(v^{+}_{m}(z)|2\tau^{\prime})
\theta_{3}(v^{+}_{m}(z)|2\tau^{\prime})}\\
&\\
&=\dfrac{\theta_{3}(v^{+}_{m}(z)|\tau^{\prime})
\theta_{3}(0 |\tau^{\prime})-
\theta_{0}(v^{+}_{m}(z)|\tau^{\prime})
\theta_{0}(0|\tau^{\prime})}{
\theta_{2}(v^{+}_{m}(z)|\tau^{\prime})
\theta_{2}(0|\tau^{\prime})}\\
&\\
&=-\dfrac{\theta_{0}(v_{m}(z)|\tau^{\prime})
\theta_{3}(0 |\tau^{\prime})}{
\theta_{1}(v_{m}(z)|\tau^{\prime})
\theta_{2}(0|\tau^{\prime})}
+
\dfrac{\theta_{3}(v_{m}(z)|\tau^{\prime})
\theta_{0}(0 |\tau^{\prime})}{
\theta_{1}(v_{m}(z)|\tau^{\prime})
\theta_{2}(0|\tau^{\prime})}\\
&\\
&\displaystyle{\mathop{\longrightarrow}^{z=0}} 
-\dfrac{{\rm sn}(\psi_{m})}{i{\rm cn}(\psi_{m})}
+\dfrac{1}{i{\rm cn}(\psi_{m})}
=\dfrac{1-{\rm sn}(\psi_{m})}{i{\rm cn}(\psi_{m})},\\
\end{array}
\end{eqnarray*}
where we have used $\theta_{0}(v+\frac{1}{2})=
\theta_{3}(v), \theta_{2}(v+\frac{1}{2})=-
\theta_{1}(v), \theta_{3}(v+\frac{1}{2})=\theta_{0}(v)$
and Lemma~\ref{App3}. 
Therefore, we obtain the following.
\begin{equation*}
\left.
\left(\dfrac{\theta_{3}(v^{+}_{m}(z)|2\tau^{\prime})}{
\theta_{2}(v^{+}_{m}(z)|2\tau^{\prime})}
\right)^{\prime}\right|_{z=0}
=2iK^{\prime}\dfrac{d}{d\psi_{m}}
\left(\dfrac{i{\rm cn}(\psi_{m})}{1-{\rm sn}(\psi_{m})}
\right)
=-2K^{\prime}\dfrac{{\rm dn}(\psi_{m})}{
1-{\rm sn}(\psi_{m})},
\end{equation*}
and
\begin{equation}
\begin{array}{ll}
\left.
\dfrac{H_{m}}{F_{m}}\right|_{z=0}
&=-\dfrac{2}{k}\exp(i\varphi_{m})
\left(
\dfrac{{\rm dn}(\psi_{m})}{1-{\rm sn}(\psi_{m})}\right)
\dfrac{1}{2}\cdot
\dfrac{\theta_{0}(v_{m}|\tau^{\prime})
\theta_{3}(0 |\tau^{\prime})-
\theta_{3}(v_{m}|\tau^{\prime})
\theta_{0}(0 |\tau^{\prime})}{
2\theta_{3}(v_{m}|\tau^{\prime})
\theta_{0}(0 |\tau^{\prime})}\\
&\\
&=-\dfrac{2}{k}\exp(i\varphi_{m})
\dfrac{{\rm dn}(\psi_{m})}{1-{\rm sn}(\psi_{m})}
\dfrac{1}{4}({\rm sn}(\psi_{m})-1)
=\dfrac{1}{2k}\exp(i\varphi_{m}){\rm dn}(\psi_{m}) .\\
\end{array}
\end{equation}
We denote by $(\Gamma_{m})_{i}$ the $i$-th component
of $\Gamma_{m}$ for $i=1,2,3$.
Now, we see that
\begin{eqnarray*}
\begin{array}{ll}
(\Gamma_{m})_{1}&=
\left.\dfrac{H_{m}+H_{m}^{*}}{F_{m}}\right|_{z=0}
=\dfrac{1}{k}\cos(\varphi_{m}){\rm dn}(\psi_{m}),\\
&\\
(\Gamma_{m})_{2}&=
\left.\dfrac{1}{i}
\dfrac{H_{m}-H_{m}^{*}}{F_{m}}\right|_{z=0}
=\dfrac{1}{k}\sin(\varphi_{m}){\rm dn}(\psi_{m}).\\
\end{array}
\end{eqnarray*}
Lastly, we calculate the third component 
$(\Gamma_{m})_{3}$.
\begin{equation*}
\left.\dfrac{d}{dz}\right|_{z=0}\log F_{m}
=\dfrac{i}{kK^{\prime}}\dfrac{\theta_{3}^{\prime}
(v_{m}|\tau^{\prime})}{
\theta_{3}(v_{m}|\tau^{\prime})}
=\dfrac{i}{kK^{\prime}}\cdot 2iK^{\prime}
\dfrac{d}{d\psi_{m}}\log \theta_{3}(v_{m}|\tau^{\prime})
=-\dfrac{2}{k}
\dfrac{d}{d\psi_{m}}\log \theta_{3}(v_{m}|\tau^{\prime}).
\end{equation*} 
Here we point out the following is known.
\begin{equation*}
\theta_{3}(\frac{v}{\tau}|\tau^{\prime})=
\exp\left(\pi i(\frac{v^{2}}{\tau}-\frac{1}{4})\right)
\tau^{\frac{1}{2}}\theta_{3}(v|\tau).
\end{equation*}
We apply this formulae to our case by taking
$v=\tau v_{m}$. We then see that
\begin{eqnarray*}
\begin{array}{ll}
\theta_{3}(v_{m}|\tau^{\prime})
&=\exp\left(\pi i(\tau v_{m}^{2}-\frac{1}{4}\right)
\tau^{\frac{1}{2}}\theta_{3}(\frac{\psi_{m}}{2K}
-\frac{1}{2} |\tau)\\
&\\
&=\exp\left(\pi i(\tau v_{m}^{2}-\frac{1}{4}\right)
\tau^{\frac{1}{2}}\theta_{0}(\frac{\psi_{m}}{2K}|
\tau).\\
\end{array}
\end{eqnarray*}
We then have
\begin{eqnarray*}
\begin{array}{ll}
&\dfrac{d}{d\psi_{m}}\log \theta_{3}(v_{m} |\tau^{\prime})
=\dfrac{d}{d\psi_{m}}\left(\pi i\tau v^{2}_{m}\right)
+\dfrac{d}{d\psi_{m}}\theta_{0}(\frac{\psi_{m}}{2K}
|\tau) 
=\dfrac{\pi i}{K}v_{m}+
\displaystyle{\int_{0}^{\psi_{m}}{\rm dn}^{2}\psi
d\psi-\dfrac{E}{K}\psi_{m}}\\
&\\
&=\displaystyle{\int_{0}^{\psi_{m}}{\rm dn}^{2}\psi
d\psi}+\left(\dfrac{\pi}{2KK^{\prime}}-\dfrac{E}{K}
\right)\psi_{m}-\dfrac{\pi}{2K^{\prime}}\\
&\\
&=\displaystyle{\int_{0}^{\psi_{m}}{\rm dn}^{2}\psi
d\psi}+\left(\dfrac{E^{\prime}}{K^{\prime}}-1\right)
\psi_{m}
-\dfrac{\pi}{2K^{\prime}}
=-k^{2}\displaystyle{\int_{0}^{\psi_{m}}{\rm sn}^{2}\psi
d\psi}+\dfrac{E^{\prime}}{K^{\prime}}\psi_{m}
-\dfrac{\pi}{2K^{\prime}},\\
\end{array}
\end{eqnarray*}
where we used the Legendre equation $EK^{\prime}+
E^{\prime}K-KK^{\prime}=\dfrac{\pi}{2}$ in the 
second equality to the last.
Thus, we obtain the following.
\begin{equation}
\left.-\dfrac{1}{2}\dfrac{d}{dz}\right|_{z=0}
\log F_{m}
=\dfrac{1}{k}\dfrac{d}{d\psi_{m}}\log\theta_{3}(
v_{m}|\tau^{\prime})
=-k\displaystyle{\int_{0}^{\psi_{m}}{\rm sn}^{2}\psi
d\psi}+\dfrac{E^{\prime}}{kK^{\prime}}\psi_{m}
-\dfrac{\pi}{2kK^{\prime}}.
\end{equation}
On the other hand, when we take $w(\lambda)=\dfrac{\lambda}{k}$ in~(\ref{alpha}), 
we have from (ii) of Lemma~\ref{App2},
\begin{eqnarray}
\begin{array}{ll}
iR_{m}&=\left.i\alpha^{-1}_{m}\dfrac{\partial\alpha_{m}}{
\partial \lambda}\right|_{\lambda=\frac{kK^{\prime}}{2}}
=\left.
i\cdot\dfrac{i}{2}\eta_{m}\left({\frak p}(\frac{\lambda}{k})
+\dfrac{\zeta_{\rm W}(\omega)}{\omega}\right)
\dfrac{dw}{d\lambda}
\right|_{\lambda=\frac{kK^{\prime}}{2}}\\
&\\
&=-\dfrac{1}{2k}\eta_{m}
\left({\frak p}(\frac{\omega}{2})+\dfrac{
\zeta_{\rm W}(\omega)}{\omega}\right)
=-\dfrac{E^{\prime}}{kK^{\prime}}\eta_{m}\\
&=-\dfrac{E^{\prime}}{kK^{\prime}}\psi_{m}
+\dfrac{\pi}{2kK^{\prime}}
+mk\displaystyle{\int_{0}^{\gamma}{\rm sn}^{2}\psi d\psi}.\\
\end{array}
\end{eqnarray} 
It follows from (4.5) and (4.6) that
\begin{equation*}
(\Gamma_{m})_{3}=\left.
iR_{m}-\dfrac{1}{2}\dfrac{d}{dz}\right|_{z=0}
\log F_{m}
=\displaystyle{-k\int_{0}^{\psi_{m}}{\rm sn}^{2}\psi d\psi
+mk\int_{0}^{\gamma}{\rm sn}^{2}\psi d\psi} .
\end{equation*}
Next, we show $B_{m}$ using (3.6). 
Set $\lambda=\dfrac{kK^{\prime}}{2}$ and $z=0$.
We write $v^{\pm}_{m}=v_{m}\pm\dfrac{1}{2}$.
We first calculate
$f_{m}^{*}g_{m}$.
\begin{eqnarray*}
\begin{array}{ll}
f^{*}_{m}g_{m}&=\exp(i\varphi_{m})
\left(\theta_{3}(v^{+}_{m}|2\tau^{\prime})
-i\theta_{2}v^{+}_{m}|2\tau^{\prime})\right)
\left(\theta_{3}(v^{+}_{m}|2\tau^{\prime})
+i\theta_{2}v^{+}_{m}|2\tau^{\prime})\right)\\
&=\exp(i\varphi_{m})\left(
\theta_{3}^{2}(v^{+}_{m}|2\tau^{\prime})
+\theta_{2}^{2}(v^{+}_{m}|2\tau^{\prime})\right)\\
&=\exp(i\varphi_{m})\theta_{3}(v^{+}_{m}|\tau^{\prime})
\theta_{3}(0|\tau^{\prime})
=\exp(i\varphi_{m})\theta_{0}(v_{m}|\tau^{\prime})
\theta_{3}(0|\tau^{\prime}),
\end{array}
\end{eqnarray*}
which implies that
\begin{equation*}
\dfrac{f^{*}_{m}g_{m}}{F_{m}}=
\exp(i\varphi_{m})\dfrac{\theta_{0}(v_{m}|\tau^{\prime})
\theta_{3}(0|\tau^{\prime})}{
\theta_{3}(v_{m}|\tau^{\prime})
\theta_{0}(0|\tau^{\prime})}
=\dfrac{1}{2}\exp(i\varphi_{m}){\rm sn}(\psi_{m}).
\end{equation*}
Therefore, we see that
\begin{eqnarray*}
\begin{array}{ll}
(B_{m})_{1}&=\dfrac{f^{*}_{m}g_{m}+f_{m}g^{*}_{m}}{
F_{m}}=\cos(\varphi_{m}){\rm sn}(\psi_{m}),\\
(B_{m})_{2}&=\dfrac{1}{i}
\dfrac{f^{*}_{m}g_{m}-f_{m}g^{*}_{m}}{
F_{m}}=\sin(\varphi_{m}){\rm sn}(\psi_{m}).\\
\end{array}
\end{eqnarray*}
Finally, to obtain the third component of $B_{m}$
we calculate $f_{m}f^{*}_{m}-g_{m}g^{*}_{m}$.
\begin{eqnarray*}
\begin{array}{ll}
f_{m}f^{*}_{m}-g_{m}g^{*}_{m}&=
\left(\theta_{3}(v_{m}^{-}|2\tau^{\prime})
+i\theta_{2}(v_{m}^{-}|2\tau^{\prime})\right)
\left(\theta_{3}(v_{m}^{+}|2\tau^{\prime})
-i\theta_{2}(v_{m}^{+}|2\tau^{\prime})\right)\\
\quad &-
\left(\theta_{3}(v_{m}^{+}|2\tau^{\prime})
+i\theta_{2}(v_{m}^{+}|2\tau^{\prime})\right)
\left(\theta_{3}(v_{m}^{-}|2\tau^{\prime})
-i\theta_{2}(v_{m}^{-}|2\tau^{\prime})\right)\\
&=2i\left(\theta_{3}(v^{+}_{m}|2\tau^{\prime})
\theta_{2}(v^{-}_{m}|2\tau^{\prime})
-\theta_{2}(v^{+}_{m}|2\tau^{\prime})
\theta_{3}(v^{-}_{m}|2\tau^{\prime})\right)\\
&=2i\theta_{1}(v_{m}|\tau^{\prime})
\theta_{1}(\frac{1}{2}|\tau^{\prime})
=2i\theta_{1}(v_{m}|\tau^{\prime})
\theta_{2}(0|\tau^{\prime}) ,\\
\end{array}
\end{eqnarray*}
which, together with (ii) of Lemma~\ref{App3}, yields
\begin{equation*}
(B_{m})_{3}=\dfrac{f_{m}f^{*}_{m}-g_{m}g^{*}_{m}}{
F_{m}}=i\dfrac{\theta_{1}(v_{m}|\tau^{\prime})
\theta_{2}(0|\tau^{\prime})}{
\theta_{3}(v_{m}|\tau^{\prime})
\theta_{0}(0|\tau^{\prime})}
=-{\rm cn}(\psi_{m}) .
\end{equation*}
Next, we prove that $\Gamma_{m+1}-\Gamma_{m}=B_{m+1}\times B_{m}$.
For this, we calculate 
$\left.\dfrac{\Psi^{FH}_{m}}{F_{m+1}F_{m}}\right|_{z=0}$ in the same way and
obtain the following.
\begin{equation}
\left.\dfrac{1}{i}\dfrac{\Psi^{FH}_{m}}{F_{m+1}F_{m}}\right|_{z=0}
=\dfrac{i}{2}\exp(i\varphi_{m+1}){\rm sn}(\psi_{m+1}){\rm cn}(\psi_{m})
-\dfrac{i}{2}\exp(i\varphi_{m}){\rm sn}(\psi_{m}){\rm cn}(\psi_{m+1}).
\end{equation}
On the other hand, it follows from (4.4) that
\begin{equation}
\left.\left(\dfrac{H_{m+1}}{F_{m+1}}-\dfrac{H_{m}}{F_{m}}\right)\right|_{z=0}
=\dfrac{1}{2k}\left(\exp(i\varphi_{m+1}){\rm dn}(\psi_{m+1})
-\exp(i\varphi_{m}){\rm dn}(\psi_{m})\right) .
\end{equation}
Using the assumption $\exp(i\alpha)=\cos(\alpha)+i\sin(\alpha)
={\rm dn}(\gamma)+ik{\rm sn}(\gamma)$ and comparing (4.7) and (4.8),
we see from (iv) and (v) of Lemma~\ref{App7} that
\begin{equation*}
\left.\left(\dfrac{H_{m+1}}{F_{m+1}}-\dfrac{H_{m}}{F_{m}}\right)\right|_{z=0}
=\left.\dfrac{1}{i}\dfrac{\Psi^{FH}_{m}}{F_{m+1}F_{m}}\right|_{z=0},
\end{equation*}
which implies that $(\Gamma_{m+1})_{a}-(\Gamma_{m})_{a}=\left(
B_{m+1}\times B_{m}\right)_{a}$ for $a=1,2$.
The third component of $B_{m+1}\times B_{m}$ is given by
\begin{equation*}
\left(B_{m+1}\times B_{m}\right)_{3}=
-k{\rm sn}(\gamma){\rm sn}(\psi_{m}){\rm sn}(\psi_{m+1}).
\end{equation*}
On the other hand, $\left(\Gamma_{m+1}\right)_{3}-\left(\Gamma_{m}\right)_{3}$
is given by
\begin{equation*}
\left(\Gamma_{m+1}\right)_{3}-\left(\Gamma_{m}\right)_{3}
=-k\int_{\psi_{m}}^{\psi_{m+1}}{\rm sn}^{2}\psi d\psi +k\int_{0}^{\gamma}
{\rm sn}^{2}\psi d\psi .
\end{equation*}
To show the equality $(\Gamma_{m+1})_{3}-(\Gamma_{m})_{3}=\left(
B_{m+1}\times B_{m}\right)_{3}$, first we show that the derivations by $t$
of both equations coincide. We must have the following equation.
\begin{equation*}
{\rm sn}^{2}(\psi_{m+1})-{\rm sn}^{2}(\psi_{m})
={\rm sn}(\gamma){\rm cn}(\psi_{m}){\rm dn}(\psi_{m}){\rm sn}(\psi_{m+1})
+{\rm sn}(\gamma){\rm sn}(\psi_{m}){\rm cn}(\psi_{m+1}){\rm dn}(\psi_{m+1}),
\end{equation*}
which, together with the equation multiplying (iv) of Lemma~\ref{App7} by
${\rm sn}(\psi_{m}){\rm cn}(\psi_{m+1})$, yields
\begin{eqnarray*}
\begin{array}{ll}
&{\rm dn}(\gamma){\rm sn}(\psi_{m}){\rm cn}(\psi_{m}){\rm sn}(\psi_{m+1})
{\rm cn}(\psi_{m+1})+{\rm sn}^{2}(\psi_{m}){\rm sn}^{2}(\psi_{m+1})\\
&={\rm sn}^{2}(\psi_{m})+{\rm sn}(\gamma){\rm sn}(\psi_{m}){\rm cn}(\psi_{m+1})
{\rm dn}(\psi_{m+1}) \\
&={\rm sn}^{2}(\psi_{m+1})-{\rm sn}(\gamma){\rm cn}(\psi_{m}){\rm dn}(\psi_{m})
{\rm sn}(\psi_{m+1}) .\\
\end{array}
\end{eqnarray*}
The last equation is equivalent to (vi) of 
Lemma~\ref{App7}.
Therefore, we see that there is a function $f(m,\gamma)$ such that
\begin{equation*}
(\Gamma_{m+1})_{3}-(\Gamma_{m})_{3}+f(m,\gamma)=\left(
B_{m+1}\times B_{m}\right)_{3}.
\end{equation*}
Setting $t=0$ in the equation above, 
we easily see that $f(m,\gamma)$ dose not depend on $m$.
Therefore, we see that $f(m,\gamma)=f(0,\gamma)=0$, hence
we are done.
Now, we show the rest properties of $\Gamma_{m}, B_{m}$ and
$T_{m}$.
$<B_{m}, B_{m+1}>={\rm cn}(\gamma)$ follows from (i) of Lemma~\ref{App7}.
Using $\cos(\alpha)={\rm dn}(\gamma), \sin(\alpha)=k{\rm sn}(\gamma)$
and (iv) of Lemma~\ref{App7}, we obtain
\begin{eqnarray*}
\begin{array}{ll}
B_{m+1}\times B_{m}&=\left(\begin{matrix}
-\sin(\varphi_{m+1}){\rm sn}(\psi_{m+1}){\rm cn}(\psi_{m})
+\sin(\varphi_{m}){\rm sn}(\psi_{m}){\rm cn}(\psi_{m+1})\\
\cos(\varphi_{m+1}){\rm sn}(\psi_{m+1}){\rm cn}(\psi_{m})
-\cos(\varphi_{m}){\rm sn}(\psi_{m}){\rm cn}(\psi_{m+1})\\
-k{\rm sn}(\gamma){\rm sn}(\psi_{m}){\rm sn}(\psi_{m+1})\\
\end{matrix}\right)\\
&\\
&=\left(\begin{matrix}
-{\rm sn}(\gamma)\sin(\varphi_{m}){\rm dn}(\psi_{m+1})-k{\rm sn}(\gamma)
\cos(\varphi_{m}){\rm sn}(\psi_{m+1}){\rm cn}(\psi_{m})\\
{\rm sn}(\gamma)\cos(\varphi_{m}){\rm dn}(\psi_{m+1})-k{\rm sn}(\gamma)
\sin(\varphi_{m}){\rm sn}(\psi_{m+1}){\rm cn}(\psi_{m})\\
-k{\rm sn}(\gamma){\rm sn}(\psi_{m}){\rm sn}(\psi_{m+1})\\
\end{matrix}\right).
\end{array}
\end{eqnarray*} 
We consider the case where the sign in the definition~(\ref{frenetframe1})
is \lq\lq$+$\rq\rq. In this case, we choose $\gamma$ so that 
${\rm sn}(\gamma)>0$.
Since $\Gamma_{m+1}-\Gamma_{m}=B_{m+1}\times B_{m}$, divide it by ${\rm sn}(\gamma)$
we have $T_{m}$ and $N_{m}=B_{m}\times T_{m}$ as follows.
\begin{eqnarray*}
\begin{array}{ll}
T_{m}&=\left(\begin{matrix}
-\sin(\varphi_{m}){\rm dn}(\psi_{m+1})
-k\cos(\varphi_{m}){\rm sn}(\psi_{m+1}){\rm cn}(\psi_{m})\\
\cos(\varphi_{m}){\rm dn}(\psi_{m+1})
-k\sin(\varphi_{m}){\rm sn}(\psi_{m+1}){\rm cn}(\psi_{m})\\
-k{\rm sn}(\psi_{m}){\rm sn}(\psi_{m+1})\\
\end{matrix}\right) ,\\
&\\
N_{m}&=\left(\begin{matrix}
-k\sin(\varphi_{m}){\rm sn}(\psi_{m+1})
+\cos(\varphi_{m}){\rm cn}(\psi_{m}){\rm dn}(\psi_{m+1})\\
k\cos(\varphi_{m}){\rm sn}(\psi_{m+1})
+\sin(\varphi_{m}){\rm cn}(\psi_{m}){\rm dn}(\psi_{m+1})\\
{\rm sn}(\psi_{m}){\rm dn}(\psi_{m+1})\\
\end{matrix}\right) ,\\
\end{array}
\end{eqnarray*}
which are compatible with the definition~(\ref{frenetframe1}).
For this, we use the following.
\begin{eqnarray*}
\begin{array}{ll}
\left.
\Phi_{m}E_{1}\Phi_{m}^{-1}\right|_{z=0,\lambda=\frac{kK^{\prime}}{2}}
&=\dfrac{1}{F_{m}}\left(\begin{matrix}{\rm Re}((f^{*}_{m})^{2}
-(g_{m})^{2})\\
{\rm Im}((f^{*}_{m})^{2}
-(g_{m})^{2})\\
-(f_{m}g_{m}+f^{*}_{m}g^{*}_{m})\\
\end{matrix}\right)
=\left(\begin{matrix}-\cos(\varphi_{m}){\rm cn}(\psi_{m})\\
-\sin(\varphi_{m}){\rm cn}(\psi_{m})\\
-{\rm sn}(\psi_{m})\\
\end{matrix}\right), \\
&\\
\left.
\Phi_{m}E_{2}\Phi_{m}^{-1}\right|_{z=0,\lambda=\frac{kK^{\prime}}{2}}
&=\dfrac{1}{F_{m}}\left(\begin{matrix}{\rm Re}(i(f^{*}_{m})^{2}
+i(g_{m})^{2})\\
{\rm Im}(i(f^{*}_{m})^{2}
+i(g_{m})^{2})\\
i(f_{m}g_{m}-f^{*}_{m}g^{*}_{m})\\
\end{matrix}\right)
=\left(\begin{matrix}
-\sin(\varphi_{m})\\
\cos(\varphi_{m})\\
0\\
\end{matrix}\right) .\\
\end{array}
\end{eqnarray*}
We then see that $\sin\nu=<B_{m+1}, N_{m}>
=<B_{m+1}\times B_{m}, T_{m}>
={\rm sn}(\gamma)>0$, 
which also follows from (iii) of Lemma~\ref{App7} directly.
For the case where the sign in the definition~(\ref{frenetframe1}) is 
\lq\lq$-$\rq\rq,
we choose $\gamma$ so that ${\rm sn}(\gamma)<0$ and define
$T_{m}$ by $T_{m}=-\frac{1}{{\rm sn}(\gamma)}B_{m+1}\times B_{m}$.
We then have $\sin\nu=-{\rm sn}(\gamma)>0$. 
Under the condition $\cos(\frac{\hat{w}}{2})={\rm dn}(
\psi_{m}), \sin(\frac{\hat{w}}{2})=-k{\rm sn}(\psi_{m})$, 
to calculate ${\cal K}_{m+1}$ we use (v) and (iv) of Lemma~\ref{App7}.
\begin{eqnarray*}
\begin{array}{ll}
&-\sin{\cal K}_{m+1}=<N_{m+1}, T_{m}>\\
&=k{\rm dn}(\gamma){\rm sn}(\psi_{m+2}){\rm dn}(\psi_{m+1})+k{\rm sn}(\gamma){\rm cn}(\psi_{m+1}){\rm dn}(\psi_{m+1}){\rm dn}(\psi_{m+2})\\
&\quad +k^{3}{\rm sn}(\gamma){\rm sn}(\psi_{m+1}){\rm sn}(\psi_{m+2})
{\rm cn}(\psi_{m})-k{\rm dn}(\gamma){\rm sn}(\psi_{m+1}){\rm cn}(\psi_{m}){\rm cn}(\psi_{m+1}){\rm dn}(\psi_{m+2})\\
&\qquad -k{\rm sn}(\psi_{m}){\rm sn}^{2}(\psi_{m+1})
{\rm dn}(\psi_{m+2})\\
&=k{\rm sn}(\psi_{m+2}){\rm dn}(\psi_{m})-k{\rm sn}(\psi_{m})
{\rm cn}^{2}(\psi_{m+1}){\rm dn}(\psi_{m+2})-
k{\rm sn}(\psi_{m})
{\rm sn}^{2}(\psi_{m+1}){\rm dn}(\psi_{m+2})\\
&=k{\rm sn}(\psi_{m+2}){\rm dn}(\psi_{m})
-k{\rm sn}(\psi_{m}){\rm dn}(\psi_{m+2}).\\
\end{array}
\end{eqnarray*}
Therefore, we see that ${\cal K}_{m+1}=\frac{1}{2}\hat{w}_{m+2}-
\frac{1}{2}\hat{w}_{m}$.
Finally, calculating the differential of $\Gamma_{m}(t)$ by $t$, we have
\begin{equation*}
\dfrac{1}{\beta}\dfrac{d\Gamma_{m}}{dt}=
\left(\begin{matrix}
-\sin(\varphi_{m}){\rm dn}(\psi_{m})
-k\cos(\varphi_{m}){\rm sn}(\psi_{m}){\rm cn}(\psi_{m})\\
\cos(\varphi_{m}){\rm dn}(\psi_{m})
-k\sin(\varphi_{m}){\rm sn}(\psi_{m}){\rm cn}(\psi_{m})\\
-k{\rm sn}^{2}(\psi_{m})\\
\end{matrix}\right),
\end{equation*}
which implies that $\dfrac{d\Gamma_{m}}{dt}\perp B_{m}$ and
\begin{eqnarray*}
\begin{array}{ll}
\pm<\dfrac{1}{\beta}\dfrac{d\Gamma_{m}}{dt},T_{m}>&=
{\rm dn}(\psi_{m}){\rm dn}(\psi_{m+1})+k^{2}{\rm sn}
(\psi_{m}){\rm sn}(\psi_{m+1})
=\cos(\frac{\hat{w}_{m}}{2}-\frac{\hat{w}_{m+1}}{2})=\cos(w_{m}),\\
\pm<\dfrac{1}{\beta}\dfrac{d\Gamma_{m}}{dt},N_{m}>&=
k{\rm sn}(\psi_{m+1}){\rm dn}(\psi_{m})-k{\rm sn}(\psi_{m}){\rm dn}(\psi_{m+1})
=\sin(\frac{\hat{w}_{m}}{2}-\frac{\hat{w}_{m+1}}{2})=\sin(w_{m}).\\
\end{array}
\end{eqnarray*}
\qed

Similarly, in Lemma 4.2, we change $\varphi_{m}$ by
$\varphi_{m}-m\pi$.
\begin{Lemma} We choose $\tau$-functions as follows.
\begin{eqnarray*}
\left\{
\begin{array}{ll}
f_{m}&=i^{m}\exp(-\frac{i}{2}\varphi_{m})
\left(\theta_{3}(v^{-}_{m}(\lambda,z) | 2\tau^{\prime})+i 
\theta_{2}(v^{-}_{m}(\lambda,z) | 2\tau^{\prime})\right) ,\\
&\\
g_{m}&=(-i)^{m}\exp(\frac{i}{2}\varphi_{m})
\left(\theta_{3}(v^{+}_{m}(\lambda,z) | 2\tau^{\prime})+i 
\theta_{2}(v^{+}_{m}(\lambda,z) | 2\tau^{\prime})\right) ,\\
&\\
F_{m}&=f_{m}f^{*}_{m}+g_{m}g^{*}_{m}=
2\theta_{3}(v_{m}(z) |\tau^{\prime})
\theta_{0}(0 | \tau^{\prime}) ,\\
&\\
H_{m}&=\frac{(-1)^{m}}{kK^{\prime}}\exp(i\varphi_{m})\left(
\theta_{2}(v^{+}_{m}(z) |2\tau^{\prime}) \theta^{\prime}_{3}
(v^{+}_{m}(z) |2\tau^{\prime}) -
\theta_{3}(v^{+}_{m}(z) |2\tau^{\prime}) \theta^{\prime}_{2}
(v^{+}_{m}(z) |2\tau^{\prime})\right),\\
&\\
\eta_{m}&=\displaystyle{
\psi_{m}-\frac{\pi}{2E^{\prime}}-\frac{mk^{2}K^{\prime}}{E^{\prime}}
\int_{0}^{\gamma}{\rm sn}^{2}\psi d\psi} .\\
\end{array}
\right.
\end{eqnarray*}
\end{Lemma}

These satisfy (3.3) and (3.4) for $\varepsilon=-1$. In this case,
we have $\Gamma_{m+1}-\Gamma_{m}=-B_{m+1}\times B_{m}$.
We then have the following.
\begin{Corollary}\label{dncorollary}
\begin{eqnarray*}
\hskip 1cm
\Gamma_{m}&=\left(\begin{matrix} 
\frac{(-1)^{m}}{k}\cos(\varphi_{m}) 
{\rm dn}(\psi_{m})\\
&\\
\frac{(-1)^{m}}{k}\sin(\varphi_{m}) 
{\rm dn}(\psi_{m}) \\
&\\
\displaystyle{-k \int_{0}^{\psi_{m}}{\rm sn}^{2}\psi d\psi 
+km\int_{0}^{\gamma}{\rm sn}^{2}\psi d\psi}\\
\end{matrix}\right), \quad
B_{m}&=\left(\begin{matrix} (-1)^{m}\cos(\varphi_{m})
{\rm sn}(\psi_{m})\\
(-1)^{m}\sin(\varphi_{m}){\rm sn}(\psi_{m}) \\
-{\rm cn}(\psi_{m})\\
\end{matrix}\right),
\end{eqnarray*}
where $\varphi_{m}=m\alpha +\beta k t, \psi_{m}=
m\gamma +\beta t$ and
$\cos(\alpha)=-{\rm dn}(\gamma), \sin(\alpha)=
k {\rm sn}(\gamma)$. 
We have also the following properties.
\item{\rm (i)} $<B_{m}, B_{m+1}>={\rm cn}(\gamma)$,
thus $\cos\nu={\rm cn}(\gamma)$ and
$\sin\nu=\pm{\rm sn}(\gamma)<0$.
The compound sign is in the same order as the compound sign 
in~{\rm (\ref{frenetframe2})}.
\item{\rm (ii)} $\cos(\frac{\hat{w}_{m}}{2})
={\rm dn}(\psi_{m})$, 
$\sin(\frac{\hat{w}_{m}}{2})
=-k{\rm sn}(\psi_{m})$ and
${\cal K}_{m+1}=-\frac{1}{2}\hat{w}_{m+2}+\frac{1}{2}\hat{w}_{m}$.
\item{\rm (iii)} $\Gamma_{m}(t)$ is an
isoperimetric deformation of the discrete
space curve, that is, the following holds.
$\dfrac{d\Gamma_{m}}{dt}=\pm\beta\left(
\cos(w_{m}) T_{m}+\sin(w_{m}) N_{m}\right)$, where
$w_{m}=\frac{1}{2}\hat{w}_{m}+\frac{1}{2}\hat{w}_{m+1}$.
\end{Corollary}

{\it Proof}. We just show the properties.
$<B_{m+1},B_{m}>={\rm cn}(\gamma)$ follows from (i) of Lemma~\ref{App7}.
We consider the case where the sign in the definition~(\ref{frenetframe2})
is \lq\lq$+$\rq\rq. In this case, we choose $\gamma$ so that
${\rm sn}(\gamma)<0$. 
Since $\Gamma_{m+1}-\Gamma_{m}=-B_{m+1}\times B_{m}$, we define $T_{m}$
and $N_{m}$ by 
$T_{m}=\dfrac{1}{{\rm sn}(\gamma)}B_{m+1}\times B_{m}$ and 
$N_{m}=B_{m}\times T_{m}$. We then have
\begin{eqnarray*}
\begin{array}{ll}
T_{m}&=
\left(\begin{matrix}
(-1)^{m}\left(-\sin(\varphi_{m}){\rm dn}(\psi_{m+1})+k\cos(\varphi_{m})
{\rm sn}(\psi_{m+1}){\rm cn}(\psi_{m})\right)\\
(-1)^{m}\left(\cos(\varphi_{m}){\rm dn}(\psi_{m+1})+k\sin(\varphi_{m})
{\rm sn}(\psi_{m+1}){\rm cn}(\psi_{m})\right)\\
k{\rm sn}(\psi_{m}){\rm sn}(\psi_{m+1})\\
\end{matrix}\right),\\
N_{m}&=
\left(\begin{matrix}
(-1)^{m}\left(k\sin(\varphi_{m}){\rm sn}(\psi_{m+1})+\cos(\varphi_{m})
{\rm cn}(\psi_{m}){\rm dn}(\psi_{m+1})\right)\\
(-1)^{m}\left(-k\cos(\varphi_{m}){\rm sn}(\psi_{m+1})+\sin(\varphi_{m})
{\rm cn}(\psi_{m}){\rm dn}(\psi_{m+1})\right)\\
{\rm sn}(\psi_{m}){\rm dn}(\psi_{m+1})\\
\end{matrix}\right) ,\\
\end{array}
\end{eqnarray*}
which are compatible with the definition~(\ref{frenetframe2}). 
We may calculate $<N_{m+1}, T_{m}>$ as in the proof of Theorem~\ref{dnsolution} and find that 
${\cal K}_{m+1}=-\frac{1}{2}\hat{w}_{m+2}+\frac{1}{2}\hat{w}_{m}$.
Moreover, we then see that 
$\sin\nu=<B_{m+1}, N_{m}>=$ 
$<B_{m+1}\times B_{m}, T_{m}>
={\rm sn}(\gamma)<0$. 
For the case where the sign in~(\ref{frenetframe2}) is \lq\lq$-$\rq\rq,
we choose $\gamma$ so that ${\rm sn}(\gamma)>0$ and define 
$T_{m}$ by $T_{m}=-\frac{1}{{\rm sn}(\gamma)}B_{m+1}\times B_{m}$.
We then see that $\sin\nu=-{\rm sn}(\gamma)<0$. 
Finally, 
$\dfrac{1}{\beta}
\dfrac{d\Gamma_{m}}{dt}$ is given by
\begin{equation*}
\dfrac{1}{\beta}\dfrac{d\Gamma_{m}}{dt}=
\left(\begin{matrix}
(-1)^{m}\left(-\sin(\varphi_{m}){\rm dn}(\psi_{m})
-k\cos(\varphi_{m}){\rm sn}(\psi_{m}){\rm cn}(\psi_{m})\right)\\
(-1)^{m}\left(\cos(\varphi_{m}){\rm dn}(\psi_{m})
-k\sin(\varphi_{m}){\rm sn}(\psi_{m}){\rm cn}(\psi_{m})\right)\\
-k{\rm sn}^{2}(\psi_{m})\\
\end{matrix}\right) ,
\end{equation*}
which implies that 
\begin{eqnarray*}
\begin{array}{ll}
\pm<\dfrac{1}{\beta}\dfrac{d\Gamma_{m}}{dt},T_{m}>
&={\rm dn}(\psi_{m}){\rm dn}(\psi_{m+1})-k^{2}{\rm sn}(\psi_{m}){\rm sn}(\psi_{m+1}) \\
&=\cos(\frac{\hat{w}_{m}}{2}+\frac{\hat{w}_{m+1}}{2})=\cos(w_{m}),\\
\pm<\dfrac{1}{\beta}\dfrac{d\Gamma_{m}}{dt},N_{m}>
&=-k{\rm sn}(\psi_{m+1}){\rm dn}(\psi_{m})
-k{\rm sn}(\psi_{m}){\rm dn}(\psi_{m+1})\\
&=\sin(\frac{\hat{w}_{m}}{2}+\frac{\hat{w}_{m+1}}{2})=\sin(w_{m}).\\
\end{array}
\end{eqnarray*}
\qed

{\bf [Periodicity]} 
Assume that $m=2n$ and $\gamma=\pm K$ in $\Gamma_{m}$.
The torsion $\nu$ is given by $\cos\nu={\rm cn}(K)=0$.
In this case, we have $\cos\alpha=\pm{\rm dn}(K)=\pm k^{\prime}$.
Since ${\rm sn}(\psi+2K)=-{\rm sn}(\psi)$, we see that
\begin{eqnarray*}
\begin{array}{ll}
\displaystyle{
\int_{0}^{2nK+\beta t}{\rm sn}^{2}(\psi) d\psi}
&=\displaystyle{\int_{0}^{2nK}{\rm sn}^{2}(\psi) d\psi+
\int_{0}^{\beta t}{\rm sn}^{2}(\psi) d\psi}\\
&=\displaystyle{2n\int_{0}^{K}{\rm sn}^{2}(\psi) d\psi+
\int_{0}^{\beta t}{\rm sn}^{2}(\psi) d\psi}\\
\end{array}
\end{eqnarray*}
Therefore, the third element of $\Gamma_{m}$ 
is independent of $m=2n$. It follows from ${\rm dn}(2nK+\beta t)
={\rm dn}(\beta t)$ that ${\rm dn}(\psi_{m})$ is also
independent of $m=2n$. For $\cos(\varphi_{m})$ and
$\sin(\varphi_{m})$, 
these are periodic if $\alpha=\frac{\pi}{n}$ or
$\alpha=\dfrac{n-1}{n}\pi$ according as $\cos\alpha=k^{\prime}$
or $-k^{\prime}$. 
This is possible when we choose $k^{\prime}$ so that
$\frac{\pi}{n}=\arccos(k^{\prime})$. 
Thus, for arbitrary $n\geqq 3$, there is a $k^{\prime}$
such that $\Gamma_{m+2n}=\Gamma_{m}$ holds for any $m$.

\section{Semi-discrete surfaces corresponding 
to \lq\lq${\rm cn}$-solution\rq\rq of semi-discrete SG-equation}

We construct a semi-discrete space curve which corresponds
to the solution $\hat{w}_{m}=2\arccos({\rm cn}\psi_{m}))$
of the semi-discrete SG-equation.
Set
\begin{eqnarray}
\left\{
\begin{array}{ll}
&v_{m}=\dfrac{\psi_{m}-K}{2iK^{\prime}}, \quad
v_{m}(z)=v_{m}+\dfrac{i}{K^{\prime}}z,\\
&\\
&v^{+}_{m}(\lambda,z)=
v_{m}+\dfrac{1}{2}+\tau^{\prime}+
\dfrac{\lambda}{K^{\prime}}+\dfrac{i}{K^{\prime}}z ,\\
&\\
&v^{-}_{m}(\lambda,z)=
v_{m}+\dfrac{1}{2}+\tau^{\prime}-
\dfrac{\lambda}{K^{\prime}}+\dfrac{i}{K^{\prime}}z ,\\
&\\
&v^{+}_{m}(z)=v_{m}+1+\tau^{\prime}+\dfrac{i}{K^{\prime}}z,\quad
v^{-}_{m}(z)=v_{m}+\tau^{\prime}+\dfrac{i}{K^{\prime}}z .\\
\end{array}
\right.
\end{eqnarray}

We then see that 
$\left(v^{+}_{m}(\lambda,z)\right)^{*}
=-v^{-}_{m}(\lambda,z)+1$ and
$\left(v^{+}_{m}(z)\right)^{*}
=-v^{-}_{m}(z)+1$.
In this case, we may verify the following.
\begin{eqnarray}
\left\{
\begin{array}{ll}
\left(\theta_{3}(v^{+}_{m}(\lambda,z) |2\tau{\prime})\right)^{*}
&=\theta_{3}(-v^{-}_{m}(\lambda,z)+1 | 2\tau^{\prime})=
\theta_{3}(v^{-}_{m}(\lambda,z) |2\tau^{\prime}),\\
&\\
\left(\theta_{2}(v^{+}_{m}(\lambda,z) |2\tau{\prime})\right)^{*}
&=\theta_{2}(-v^{-}_{m}(\lambda,z)+1 | 2\tau^{\prime})=
-\theta_{2}(v^{-}_{m}(\lambda,z) |2\tau^{\prime}).\\
\end{array}
\right.
\end{eqnarray}
For $v^{\pm}_{m}(z)$, the same equations hold as above.

\begin{Lemma} Set $\varphi_{m}=m\alpha+\beta t$ and $\psi_{m}=
m\gamma +\beta t$. We choose $\tau$-functions as follows.
\begin{eqnarray*}
\left\{
\begin{array}{ll}
f_{m}&=\exp(-\frac{i}{2}\varphi_{m})
\left(\theta_{3}(v^{-}_{m}(\lambda,z) | 2\tau^{\prime})+ 
\theta_{2}(v^{-}_{m}(\lambda,z) | 2\tau^{\prime})\right) ,\\
&\\
g_{m}&=\exp(\frac{i}{2}\varphi_{m})
\left(\theta_{3}(v^{+}_{m}(\lambda,z) | 2\tau^{\prime})+ 
\theta_{2}(v^{+}_{m}(\lambda,z) | 2\tau^{\prime})\right) ,\\
&\\
F_{m}&=f_{m}f^{*}_{m}+g_{m}g^{*}_{m}=
2\beta_{m}\theta_{3}(v_{m}(z) |\tau^{\prime})
\theta_{3}(0 | \tau^{\prime}) ,\\
&{\rm with}\hskip .1cm
\beta_{m}=\exp(-2\pi iv_{m}+\frac{2\pi}{K^{\prime}}z
+\frac{K}{K^{\prime}}\pi) ,\\
&\\
H_{m}&=\dfrac{-i}{K^{\prime}}\exp(i\varphi_{m})
\left(
\theta_{2}(v^{+}_{m}(z) |2\tau^{\prime}) \theta^{\prime}_{3}
(v^{+}_{m}(z) |2\tau^{\prime}) -
\theta_{3}(v^{+}_{m}(z) |2\tau^{\prime}) \theta^{\prime}_{2}
(v^{+}_{m}(z) |2\tau^{\prime})\right),\\
&\\
\eta_{m}&=\displaystyle{
\psi_{m}-\frac{3\pi}{2E^{\prime}}-\frac{mk^{2}K^{\prime}}{E^{\prime}}
\int_{0}^{\gamma}{\rm sn}^{2}\psi d\psi} .\\
\end{array}
\right.
\end{eqnarray*}
\end{Lemma}
We then have $\exp(i\dfrac{\hat{w}_{m}}{2})
=\left.-i\dfrac{g_{m}}{f_{m}^{*}}
\right|_{\lambda=\frac{K^{\prime}}{2},
z=i\frac{K^{\prime}}{2}}
={\rm cn}(\psi_{m})+i{\rm sn}(\psi_{m})$, hence
$\hat{w}_{m}$ is a \lq\lq{\rm cn}-solution\rq\rq. 
We may verify that these functions satisfy (3.3) and (3.4) for
$\varepsilon=+1$. 

{\it Proof}. Let $\hat{w}_{m}$ be a \lq\lq{\rm cn}-solution\rq\rq.
Then, for the same choice of $\xi_{m}$ as that in the proof of Lemma 4.2
except for the chice of $\xi_{0}=0$,
we have $\sin\left(\frac{\hat{w}_{m}}{2}\right)={\rm sn}(\psi_{m})$
and $\cos\left(\frac{\hat{w}_{m}}{2}\right)={\rm cn}(\psi_{m})$.
On the other hand,
\begin{eqnarray*}
\begin{array}{ll}
\left.-i\dfrac{g_{m}}{f_{m}^{*}}
\right|_{\lambda=\frac{K^{\prime}}{2},
z=i\frac{K^{\prime}}{2}}&=
-i\dfrac{\theta_{3}(v_{m}+\frac{1}{2}+\tau^{\prime}|2\tau^{\prime})+
\theta_{2}(v_{m}+\frac{1}{2}+\tau^{\prime}|2\tau^{\prime})}{
\theta_{3}(v_{m}+\frac{1}{2}+\tau^{\prime}|2\tau^{\prime})-
\theta_{2}(v_{m}+\frac{1}{2}+\tau^{\prime}|2\tau^{\prime})}\\
&=-i\dfrac{\theta_{3}(v_{m}+\frac{1}{2}+\tau^{\prime}|\tau^{\prime})
\theta_{3}(0|\tau^{\prime})+
\theta_{2}(v_{m}+\frac{1}{2}+\tau^{\prime}|\tau^{\prime})
\theta_{2}(0|\tau^{\prime})}{
\theta_{0}(v_{m}+\frac{1}{2}+\tau^{\prime}|\tau^{\prime})
\theta_{0}(0|\tau^{\prime})}\\
&=-i\dfrac{-\theta_{0}(v_{m}|\tau^{\prime})
\theta_{3}(0|\tau^{\prime})+
\theta_{1}(v_{m}|\tau^{\prime})
\theta_{2}(0|\tau^{\prime})}{
\theta_{3}(v_{m}|\tau^{\prime})
\theta_{0}(0|\tau^{\prime})}\\
&=i{\rm sn}(\psi_{m})+{\rm cn}(\psi_{m}),\\
\end{array}
\end{eqnarray*}
where we have used Lemma~\ref{App2} and 
Lemma~\ref{App3}.
For $F_{m}$,we use the fourth equation of 
Lemma~\ref{App2}
for $\lambda=\frac{K^{\prime}}{2}$.
\begin{eqnarray*}
\begin{array}{ll}
F_{m}&=(\theta_{3}(v^{-}_{m}(\lambda,z)|2\tau^{\prime})
+\theta_{2}(v^{-}_{m}(\lambda,z)|2\tau^{\prime})
(\theta_{3}(v^{+}_{m}(\lambda,z)|2\tau^{\prime})
-\theta_{2}(v^{+}_{m}(\lambda,z)|2\tau^{\prime})\\
&\quad \left.+(\theta_{3}(v^{+}_{m}(\lambda,z)|2\tau^{\prime})
+\theta_{2}(v^{+}_{m}(\lambda,z)|2\tau^{\prime})
(\theta_{3}(v^{-}_{m}(\lambda,z)|2\tau^{\prime})
-\theta_{2}(v^{-}_{m}(\lambda,z)|2\tau^{\prime})
\right|_{\lambda=\frac{K^{\prime}}{2}}\\
&\\
&=2\left(\theta_{3}(v^{+}_{m}(z)|2\tau^{\prime})
\theta_{3}(v^{-}_{m}(z)|2\tau^{\prime})-
\theta_{2}(v^{+}_{m}(z)|2\tau^{\prime})
\theta_{2}(v^{-}_{m}(z)|2\tau^{\prime})\right)
\\
&\\
&=2\theta_{0}(v_{m}(z)+\frac{1}{2}+\tau^{\prime} |\tau^{\prime})
\theta_{0}(\frac{1}{2} |\tau^{\prime})
=2\beta_{m}\theta_{3}(v_{m}(z) |\tau^{\prime})
\theta_{3}(0 |\tau^{\prime}),\\
\end{array}
\end{eqnarray*}
where $\beta_{m}=\exp(-\pi i(2v_{m}(z)+\tau^{\prime}))
=\exp(-2\pi iv_{m}+\frac{2\pi}{K^{\prime}}z+
\frac{K}{K^{\prime}}\pi)$.
For $H_{m}$, considering (3.4) we calculate as follows.
\begin{eqnarray*}
\begin{array}{ll}
H_{m}&=\left.\dfrac{i}{2}D_{\lambda}g_{m}\cdot f^{*}_{m}
\right|_{\lambda=\frac{K^{\prime}}{2}}\\
&=\dfrac{i}{2K^{\prime}}\exp(i\varphi_{m})\left\{
\left(\theta_{3}^{\prime}(v^{+}_{m}(z)|2\tau^{\prime})
+\theta_{2}^{\prime}(v^{+}_{m}(z)|2\tau^{\prime})
\right)\left(
\theta_{3}(v^{+}_{m}(z)|2\tau^{\prime})
-\theta_{2}(v^{+}_{m}(z)|2\tau^{\prime})
\right)\right. \\
&\qquad\qquad\qquad\left. -
\left(\theta_{3}(v^{+}_{m}(z)|2\tau^{\prime})
+\theta_{2}(v^{+}_{m}(z)|2\tau^{\prime})
\right)\left(
\theta_{3}^{\prime}(v^{+}_{m}(z)|2\tau^{\prime})
-\theta_{2}^{\prime}(v^{+}_{m}(z)|2\tau^{\prime})
\right)
\right\}\\
&\\
&=\dfrac{-i}{K^{\prime}}\exp(i\varphi_{m})
\left(\theta_{2}(v^{+}_{m}(z)|2\tau^{\prime})
\theta_{3}^{\prime}(v^{+}_{m}(z)|2\tau^{\prime})
-
\theta_{3}(v^{+}_{m}(z)|2\tau^{\prime})
\theta_{2}^{\prime}(v^{+}_{m}(z)|2\tau^{\prime})
\right) .\\
\end{array}
\end{eqnarray*}
The rest two equations in (3.3) are prepared
to ensure that $\Gamma_{m+1}-\Gamma_{m}=
B_{m+1}\times B_{m}$. Therefore, we prove them
in the proof of the next Theorem.
\qed

We then have the following.

\begin{Theorem}\label{cnsolution}
\begin{eqnarray*}
\hskip 1cm
\Gamma_{m}&=\left(\begin{matrix} k \cos(\varphi_{m}) 
{\rm cn}(\psi_{m})\\
k \sin(\varphi_{m}) {\rm cn}(\psi_{m}) \\
\displaystyle{-k^{2}\int_{0}^{\psi_{m}}{\rm sn}^{2}\psi d\psi 
+k^{2}m\int_{0}^{\gamma}{\rm sn}^{2}\psi d\psi}\\
\end{matrix}\right), \quad
B_{m}&=\left(\begin{matrix} -k \cos(\varphi_{m})
{\rm sn}(\psi_{m})\\
-k \sin(\varphi_{m}){\rm sn}(\psi_{m}) \\
{\rm dn}(\psi_{m})\\
\end{matrix}\right),
\end{eqnarray*}
where $\varphi_{m}=m\alpha +\beta t, \psi_{m}=
m\gamma +\beta t$ and
$\cos(\alpha)={\rm cn}(\gamma), \sin(\alpha)=
{\rm sn}(\gamma)$. 
We have also the following properties.
\item{\rm (i)} $<B_{m}, B_{m+1}>={\rm dn}(\gamma)$,
thus $\cos\nu={\rm dn}(\gamma)$ and
$\sin\nu=\pm k{\rm sn}(\gamma)>0$.
The compound sign is in the same order 
as the compound sign in~{\rm (\ref{frenetframe1})}.
\item{\rm (ii)} 
$\cos(\frac{\hat{w}_{m}}{2})={\rm cn}(\psi_{m})$, 
$\sin(\frac{\hat{w}_{m}}{2})={\rm sn}(\psi_{m})$ and 
${\cal K}_{m+1}=\frac{1}{2}\hat{w}_{m+2}-\frac{1}{2}\hat{w}_{m}$.
\item{\rm (iii)} 
$\Gamma_{m}(t)$ is an 
isoperimetric deformation of the discrete 
space curve, that is, the following
holds.
$\dfrac{d\Gamma_{m}}{dt}=
\pm\beta k\left(
\cos(w_{m}) T_{m}+\sin(w_{m}) N_{m}\right)$, where
$w_{m}=-\frac{1}{2}\hat{w}_{m}+\frac{1}{2}\hat{w}_{m+1}$. 
\end{Theorem}

{\it Proof}. Sinec the essential part of the proof is almost same as 
that of Theorem~\ref{dnsolution}, we here present
the outline.
\begin{eqnarray*}
\begin{array}{ll}
&\quad\dfrac{\theta_{2}(v^{+}_{m}(z)|2\tau^{\prime})}{
\theta_{3}(v^{+}_{m}(z)|2\tau^{\prime})}
=\dfrac{2\theta_{2}^{2}(v^{+}_{m}(z)|2\tau^{\prime})}{
2\theta_{2}(v^{+}_{m}(z)|2\tau^{\prime})
\theta_{3}(v^{+}_{m}(z)|2\tau^{\prime})}\\
&\\
&=\dfrac{\theta_{3}(v^{+}_{m}(z)|\tau^{\prime})
\theta_{3}(0 |\tau^{\prime})-
\theta_{0}(v^{+}_{m}(z)|\tau^{\prime})
\theta_{0}(0|\tau^{\prime})}{
\theta_{2}(v^{+}_{m}(z)|\tau^{\prime})
\theta_{2}(0|\tau^{\prime})}\\
&\\
&=-\dfrac{\theta_{3}(v_{m}(z)|\tau^{\prime})
\theta_{3}(0 |\tau^{\prime})}{
\theta_{2}(v_{m}(z)|\tau^{\prime})
\theta_{2}(0|\tau^{\prime})}
-
\dfrac{\theta_{0}(v_{m}(z)|\tau^{\prime})
\theta_{0}(0 |\tau^{\prime})}{
\theta_{2}(v_{m}(z)|\tau^{\prime})
\theta_{2}(0|\tau^{\prime})}\\
&\\
&\displaystyle{\mathop{\longrightarrow}^{z=0}} 
-\dfrac{1}{{\rm dn}(\psi_{m})}
-\dfrac{k{\rm sn}(\psi_{m})}{{\rm dn}(\psi_{m})}
=-\dfrac{1+k{\rm sn}(\psi_{m})}{{\rm dn}(\psi_{m})}.\\
\end{array}
\end{eqnarray*}
Therefore, we obtain the following.
\begin{equation*}
\left.
\left(\dfrac{\theta_{3}(v^{+}_{m}(z)|2\tau^{\prime})}{
\theta_{2}(v^{+}_{m}(z)|2\tau^{\prime})}
\right)^{\prime}\right|_{z=0}
=2iK^{\prime}\dfrac{d}{d\psi_{m}}
\left(-\dfrac{{\rm dn}(\psi_{m})}{
1+k{\rm sn}(\psi_{m})}
\right)
=\dfrac{2ikK^{\prime}{\rm cn}(\psi_{m})}{
1+k{\rm sn}(\psi_{m})},
\end{equation*}
and
\begin{equation}
\begin{array}{ll}
\left.
\dfrac{H_{m}}{F_{m}}\right|_{z=0}
&=-\dfrac{i}{K^{\prime}}\exp(i\varphi_{m})
\left(
\dfrac{2ikK^{\prime}{\rm cn}(\psi_{m})}{1+k{\rm sn}(\psi_{m})}\right)
\dfrac{1}{4}\cdot
\left\{1+
\dfrac{\theta_{0}(v_{m}|\tau^{\prime})
\theta_{0}(0 |\tau^{\prime})}{
\theta_{3}(v_{m}|\tau^{\prime})
\theta_{3}(0 |\tau^{\prime})}\right\}\\
&\\
&=\dfrac{k}{2}\exp(i\varphi_{m})
{\rm cn}(\psi_{m}) .\\
\end{array}
\end{equation}
We then see that 
$(\Gamma_{m})_{1}=k\cos(\varphi_{m}){\rm cn}(\psi_{m})$ and
$(\Gamma_{m})_{2}=k\sin(\varphi_{m}){\rm cn}
(\psi_{m})$. For $(\Gamma_{m})_{3}$, 
the calculation is alomost same as that of Theorem 4.3. 
The only difference is to take $w(\lambda)=\lambda, \lambda
=\frac{K^{\prime}}{2}$.
Next, we prove that $\Gamma_{m+1}-\Gamma_{m}=B_{m+1}\times B_{m}$.
For this, we calculate 
$\left.\dfrac{\Psi^{FH}_{m}}{F_{m+1}F_{m}}\right|_{z=0}$ 
in the same way and obtain the following.
\begin{equation}
\left.\dfrac{1}{i}\dfrac{\Psi^{FH}_{m}}{F_{m+1}F_{m}}\right|_{z=0}
=\dfrac{k}{2}i\exp(i\varphi_{m+1}){\rm sn}(\psi_{m+1}){\rm dn}(\psi_{m})
-\dfrac{k}{2}i\exp(i\varphi_{m}){\rm sn}(\psi_{m}){\rm dn}(\psi_{m+1}).
\end{equation}
On the other hand, it follows from (5.5) that
\begin{equation}
\left.\left(\dfrac{H_{m+1}}{F_{m+1}}-\dfrac{H_{m}}{F_{m}}\right)\right|_{z=0}
=\dfrac{k}{2}\left(\exp(i\varphi_{m+1}){\rm cn}(\psi_{m+1})
-\exp(i\varphi_{m}){\rm cn}(\psi_{m})\right) .
\end{equation}
Using the assumption $\exp(i\alpha)=\cos(\alpha)+i\sin(\alpha)
={\rm cn}(\gamma)+i{\rm sn}(\gamma)$ and comparing (5.6) and (5.7),
we see from (vii) and (viii) of Lemma~\ref{App7} that
\begin{equation*}
\left.\left(\dfrac{H_{m+1}}{F_{m+1}}-\dfrac{H_{m}}{F_{m}}\right)\right|_{z=0}
=\left.\dfrac{1}{i}\dfrac{\Psi^{FH}_{m}}{F_{m+1}F_{m}}\right|_{z=0},
\end{equation*}
which implies that $(\Gamma_{m+1})_{a}-(\Gamma_{m})_{a}=\left(
B_{m+1}\times B_{m}\right)_{a}$ for $a=1,2$.
The proof for $(\Gamma_{m+1})_{3}-(\Gamma_{m})_{3}=
\left(B_{m+1}\times B_{m}\right)_{3}$ is same as that of Theorem 4.3.
The rest properties of $\Gamma_{m}, B_{m}$ and
$T_{m}$ are easily verified by direct calculations. 
We consider the case where the sign in 
the definition~(\ref{frenetframe1}) is 
\lq\lq$+$\rq\rq. 
In this case, we choose $\gamma$ so that
$k{\rm sn}(\gamma)>0$. Since
$\Gamma_{m+1}-\Gamma_{m}=B_{m+1}\times B_{m}$,
divide it by $k{\rm sn}(\gamma)$ we have
$T_{m}$ and $N_{m}=B_{m}\times T_{m}$ as follows.
\begin{eqnarray*}
\begin{array}{ll}
T_{m}&=
\left(\begin{matrix}
-\sin(\varphi_{m}){\rm cn}(\psi_{m+1})-\cos(\varphi_{m})
{\rm sn}(\psi_{m+1}){\rm dn}(\psi_{m})\\
\cos(\varphi_{m}){\rm cn}(\psi_{m+1})-\sin(\varphi_{m})
{\rm sn}(\psi_{m+1}){\rm dn}(\psi_{m})\\
-k{\rm sn}(\psi_{m}){\rm sn}(\psi_{m+1})\\
\end{matrix}\right) ,\\
N_{m}&=
\left(\begin{matrix}
\sin(\varphi_{m}){\rm sn}(\psi_{m+1})-\cos(\varphi_{m})
{\rm cn}(\psi_{m+1}){\rm dn}(\psi_{m})\\
-\cos(\varphi_{m}){\rm sn}(\psi_{m+1})-\sin(\varphi_{m})
{\rm cn}(\psi_{m+1}){\rm dn}(\psi_{m})\\
-k{\rm sn}(\psi_{m}){\rm cn}(\psi_{m+1})\\
\end{matrix}\right) ,\\
\dfrac{1}{\beta k}\dfrac{d\Gamma_{m}}{dt}&=
\left(\begin{matrix}
-\sin(\varphi_{m}){\rm cn}(\psi_{m})-\cos(\varphi_{m})
{\rm sn}(\psi_{m}){\rm dn}(\psi_{m})\\
\cos(\varphi_{m}){\rm cn}(\psi_{m})-\sin(\varphi_{m})
{\rm sn}(\psi_{m}){\rm dn}(\psi_{m})\\
-k{\rm sn}^{2}(\psi_{m})\\
\end{matrix}\right) ,\\
\end{array}
\end{eqnarray*}
which are compatible with 
the definition~(\ref{frenetframe1}).
We then see that $\sin\nu=<B_{m+1}\times
B_{m}, T_{m}>=k{\rm sn}(\gamma)>0$. 
For the case where the sign in the 
definition~(\ref{frenetframe1}) is 
\lq\lq$-$\rq\rq, we choose $\gamma$ so 
that $k{\rm sn}(\gamma)<0$ and define
$T_{m}$ by 
$T_{m}=-\frac{1}{k{\rm sn}(\gamma)}
B_{m+1}\times B_{m}$. We then see that
$\sin\nu=-k{\rm sn}(\gamma)>0$.
For ${\cal K}_{m+1}$, we calculate $<N_{m+1}, T_{m}>$ using
(v) and (iv) of Lemma~\ref{App7}.
\begin{eqnarray*}
\begin{array}{ll}
&-\sin{\cal K}_{m+1}=<N_{m+1}, T_{m}>\\
&=-{\rm cn}(\gamma){\rm sn}(\psi_{m+2}){\rm cn}(\psi_{m+1})
-{\rm sn}(\gamma){\rm sn}(\psi_{m+1}){\rm sn}(\psi_{m+2}){\rm dn}(\psi_{m})\\
&\quad -{\rm sn}(\gamma){\rm cn}(\psi_{m+1}){\rm cn}(\psi_{m+2})
{\rm dn}(\psi_{m+1})+{\rm cn}(\gamma){\rm sn}(\psi_{m+1}){\rm cn}(\psi_{m+2}){\rm dn}(\psi_{m}){\rm dn}(\psi_{m+1})\\
&\qquad +k^{2}{\rm sn}(\psi_{m}){\rm sn}^{2}(\psi_{m+1})
{\rm cn}(\psi_{m+2})\\
&=-{\rm cn}(\psi_{m}){\rm sn}(\psi_{m+2})+{\rm sn}(\psi_{m})
{\rm cn}(\psi_{m+2}){\rm dn}^{2}(\psi_{m+1})+
k^{2}{\rm sn}(\psi_{m})
{\rm sn}^{2}(\psi_{m+1}){\rm cn}(\psi_{m+2})\\
&=-{\rm cn}(\psi_{m}){\rm sn}(\psi_{m+2})
+{\rm sn}(\psi_{m}){\rm cn}(\psi_{m+2}).\\
\end{array}
\end{eqnarray*}
Threfore, we see that ${\cal K}_{m+1}=
\frac{1}{2}\hat{w}_{m+2}-\frac{1}{2}\hat{w}_{m}$.
\qed

Similarly, in Lemma 5.3, we change $\varphi_{m}$ by
$\varphi_{m}-m\pi$.

\begin{Lemma} We choose $\tau$-functions as follows.
\begin{eqnarray*}
\left\{
\begin{array}{ll}
f_{m}&=i^{m}\exp(-\frac{i}{2}\varphi_{m})
\left(\theta_{3}(v^{-}_{m}(\lambda,z) | 2\tau^{\prime})+ 
\theta_{2}(v^{-}_{m}(\lambda,z) | 2\tau^{\prime})\right) ,\\
&\\
g_{m}&=(-i)^{m}\exp(\frac{i}{2}\varphi_{m})
\left(\theta_{3}(v^{+}_{m}(\lambda,z) | 2\tau^{\prime})+ 
\theta_{2}(v^{+}_{m}(\lambda,z) | 2\tau^{\prime})\right) ,\\
&\\
F_{m}&=f_{m}f^{*}_{m}+g_{m}g^{*}_{m}=
2\beta_{m}\theta_{3}(v_{m}(z) |\tau^{\prime})
\theta_{3}(0 | \tau^{\prime}) ,\\
&{\rm with}\hskip .1cm
\beta_{m}=\exp(-2\pi iv_{m}+\frac{2\pi}{K^{\prime}}z
+\frac{K}{K^{\prime}}\pi) ,\\
&\\
H_{m}&=(-1)^{m+1}\dfrac{i}{K^{\prime}}\exp(i\varphi_{m})\left(
\theta_{2}(v^{+}_{m}(z) |2\tau^{\prime}) \theta^{\prime}_{3}
(v^{+}_{m}(z) |2\tau^{\prime}) -
\theta_{3}(v^{+}_{m}(z) |2\tau^{\prime}) \theta^{\prime}_{2}
(v^{+}_{m}(z) |2\tau^{\prime})\right),\\
&\\
\eta_{m}&=\displaystyle{
\psi_{m}-\frac{3\pi}{2E^{\prime}}-\frac{mk^{2}K^{\prime}}{E^{\prime}}
\int_{0}^{\gamma}{\rm sn}^{2}\psi d\psi} .\\
\end{array}
\right.
\end{eqnarray*}
\end{Lemma}

We then have the following.

\begin{Corollary}\label{cncorollary}
\begin{eqnarray*}
\hskip 1cm
\Gamma_{m}&=\left(\begin{matrix} 
(-1)^{m}k \cos(\varphi_{m}) 
{\rm cn}(\psi_{m})\\
(-1)^{m}k \sin(\varphi_{m}) {\rm cn}(\psi_{m}) \\
\displaystyle{-k^{2}\int_{0}^{\psi_{m}}{\rm sn}^{2}\psi d\psi 
+k^{2}m\int_{0}^{\gamma}{\rm sn}^{2}\psi d\psi}\\
\end{matrix}\right), \quad
B_{m}&=\left(\begin{matrix} -(-1)^{m}k \cos(\varphi_{m})
{\rm sn}(\psi_{m})\\
-(-1)^{m}k \sin(\varphi_{m}){\rm sn}(\psi_{m}) \\
{\rm dn}(\psi_{m})\\
\end{matrix}\right),
\end{eqnarray*}
where $\varphi_{m}=m\alpha +\beta t, \psi_{m}=
m\gamma +\beta t$ and
$\cos(\alpha)=-{\rm cn}(\gamma), \sin(\alpha)=
{\rm sn}(\gamma)$. 
We have also the following properties.
\item{\rm (i)} $<B_{m}, B_{m+1}>={\rm dn}(\gamma)$,
thus $\cos\nu={\rm dn}(\gamma)$ and 
$\sin\nu=\pm k{\rm sn}(\gamma)<0$.
The compound sign is in the same order 
as the compound sign in~{\rm (\ref{frenetframe2})}.
\item{\rm (ii)} 
$\cos(\frac{\hat{w}_{m}}{2})={\rm cn}(\psi_{m})$,  
$\sin(\frac{\hat{w}_{m}}{2})={\rm sn}(\psi_{m})$ and 
${\cal K}_{m+1}=-\frac{1}{2}\hat{w}_{m+2}+\frac{1}{2}\hat{w}_{m}$.
\item{\rm (iii)} 
$\Gamma_{m}(t)$ is an 
isoperimetric deformation of the discrete 
space curve, that is, the following
holds.
$\dfrac{d\Gamma_{m}}{dt}=\pm\beta k\left(
\cos(w_{m}) T_{m}+\sin(w_{m}) N_{m}\right)$, where
$w_{m}=-\frac{1}{2}\hat{w}_{m}-\frac{1}{2}\hat{w}_{m+1}$. 
\end{Corollary}

{\it Proof}. The main part of the proof
is almost same as that of Theorem~\ref{cnsolution}.
We consider the case where the sign in the 
definition~\ref{frenetframe2} is \lq\lq$+$\rq\rq.
Since we have $\Gamma_{m+1}-\Gamma_{m}=
-B_{m+1}\times B_{m}$, we choose $\gamma$ so 
that $k{\rm sn}(\gamma)<0$.
We define $T_{m}$ and $N_{m}$ by
$T_{m}=\frac{1}{k{\rm sn}(\gamma)}B_{m+1}
\times B_{m}$ and $N_{m}=B_{m}\times T_{m}$.
We then have
\begin{eqnarray*}
\begin{array}{ll}
T_{m}&=
\left(\begin{matrix}
(-1)^{m}\left(-\sin(\varphi_{m}){\rm cn}(\psi_{m+1})+\cos(\varphi_{m})
{\rm sn}(\psi_{m+1}){\rm dn}(\psi_{m})\right)\\
(-1)^{m}\left(\cos(\varphi_{m}){\rm cn}(\psi_{m+1})+\sin(\varphi_{m})
{\rm sn}(\psi_{m+1}){\rm dn}(\psi_{m})\right)\\
k{\rm sn}(\psi_{m}){\rm sn}(\psi_{m+1})\\
\end{matrix}\right) ,\\
N_{m}&=
\left(\begin{matrix}
(-1)^{m}\left(-\sin(\varphi_{m}){\rm sn}(\psi_{m+1})-\cos(\varphi_{m})
{\rm cn}(\psi_{m+1}){\rm dn}(\psi_{m})\right)\\
(-1)^{m}\left(\cos(\varphi_{m}){\rm sn}(\psi_{m+1})-\sin(\varphi_{m})
{\rm cn}(\psi_{m+1}){\rm dn}(\psi_{m})\right)\\
-k{\rm sn}(\psi_{m}){\rm cn}(\psi_{m+1})\\
\end{matrix}\right) ,\\
\end{array}
\end{eqnarray*}
which are compatible with 
the definition~\ref{frenetframe2}.
We then see that $\sin\nu=k{\rm sn}(\gamma)<0$. 
For the case where the sign in the 
definition~\ref{frenetframe2} is 
\lq\lq$-$\rq\rq, we choose $\gamma$
so that $k{\rm sn}(\gamma)>0$ and
define $T_{m}$ by $T_{m}=
-\frac{1}{k{\rm sn}(\gamma)}B_{m+1}
\times B_{m}$. We then see that
$\sin\nu=-k{\rm sn}(\gamma)<0$.
Moreover, $\dfrac{1}{\beta k}\dfrac{d\Gamma_{m}}{dt}$ 
is given by
\begin{equation*}
\dfrac{1}{\beta k}\dfrac{d\Gamma_{m}}{dt}=
\left(\begin{matrix}
(-1)^{m}\left(-\sin(\varphi_{m}){\rm cn}(\psi_{m})-\cos(\varphi_{m})
{\rm sn}(\psi_{m}){\rm dn}(\psi_{m})\right)\\
(-1)^{m}\left(\cos(\varphi_{m}){\rm cn}(\psi_{m})-\sin(\varphi_{m})
{\rm sn}(\psi_{m}){\rm dn}(\psi_{m})\right)\\
-k{\rm sn}^{2}(\psi_{m})\\
\end{matrix}\right) .
\end{equation*}
In the same way as the proof of Theorem~\ref{cnsolution},
we may prove the rest properties.
\qed

{\bf [Periodicity]} 
Assume that $m=2n$ and $\gamma=K$ in $\Gamma_{m}$.
The third element of $\Gamma_{m}$ is independent of
$m=2n$ as in \S 4.
The torsion $\nu$ is given by $\cos\nu={\rm dn}(K)=k^{\prime}$.
In this case, we have $\cos\alpha=\pm{\rm cn}(K)=0,
\sin\alpha={\rm sn}(K)=1$,
that is, $\alpha=\dfrac{\pi}{2}$. Since we know that
\begin{eqnarray*}
\begin{array}{ll}
\cos(n\pi+\beta t){\rm cn}(2nK+\beta t)&=
\cos(\beta t){\rm cn}(\beta t),\\
\sin(n\pi+\beta t){\rm cn}(2nK+\beta t)&=
\sin(\beta t){\rm cn}(\beta t),\\
\end{array}
\end{eqnarray*}
we have $\Gamma_{m+2}=\Gamma_{m}$ holds for any $m$.

\section{Discrete $K$-surfaces corresponding to ${\rm dn}$-
solution or ${\rm cn}$-solution of the discrete sine-Gordon equation}

As an application of the results in $\S 4$ and $\S 5$, we may construct
a discrete $K$-surface in the sense of Bobenko-Pinkall(\cite{bobpink}).
In \cite{bobpink}, a general formulae in terms of the Riemann theta functions for the immersion of discrete
$K$-surface is obtained using the theory of the spectral curve
and the Baker-Akhiezer function.
However, when one try to produce a non-trivial concrete example of
the discrete $K$-surface using their formulae, it is difficult 
to use the formulae directly. Therefore, it is realistic to limit 
the genus of the spectral curve to 1. In this point of view, we may
use the results obtained in $\S 4$ and $\S 5$.
For a discrete $K$-surface $F_{m,n}$, first we construct
the unitary frame $\Phi_{m,n}$ associated to $F_{m,n}$.
We write
\begin{equation}\label{twoframes}
\Phi_{m+1,n}=\Phi_{m,n}L_{m,n} ,\qquad
\Phi_{m,n+1}=\Phi_{m,n}\hat{L}_{m,n}.
\end{equation}
Referring to the unitary frames $L^{\pm}_{m}$ and $\hat{L}^{\pm}_{m}$ 
in Lemma~\ref{frame} and Lemma~\ref{frame2}, respectively,
we choose $L^{\pm}_{m,n}$ and $\hat{L}^{\pm}_{m,n}$ as follows.
\begin{eqnarray}
\quad
\left\{
\begin{array}{ll}
L^{\pm}_{m,n}&=
\left(\begin{matrix}\cos(\frac{\nu^{1}}{2})
e^{-\frac{i}{4}(\hat{w}_{m+1,n}
-\hat{w}_{m,n})}&
\pm\sin(\frac{\nu^{1}}{2})e^{-\frac{i}{4}(\hat{w}_{m+1,n}
+\hat{w}_{m,n})}\\
\mp\sin(\frac{\nu^{1}}{2})e^{\frac{i}{4}(\hat{w}_{m+1,n}
+\hat{w}_{m,n})}&
\cos(\frac{\nu^{1}}{2})e^{\frac{i}{4}(\hat{w}_{m+1,n}
-\hat{w}_{m,n})}\\
\end{matrix}\right), \\
&\\
\hat{L}^{\pm}_{m,n}&=
\left(\begin{matrix}\cos(\frac{\nu^{2}}{2})
e^{\frac{i}{4}(\hat{w}_{m,n+1}
+\hat{w}_{m,n})}&
\pm\sin(\frac{\nu^{2}}{2})e^{\frac{i}{4}(\hat{w}_{m,n+1}
-\hat{w}_{m,n})}\\
\mp\sin(\frac{\nu^{2}}{2})e^{-\frac{i}{4}(\hat{w}_{m,n+1}
-\hat{w}_{m,n})}&
\cos(\frac{\nu^{2}}{2})e^{-\frac{i}{4}(\hat{w}_{m,n+1}
+\hat{w}_{m,n})}\\
\end{matrix}\right) . \\
\end{array}
\right.
\end{eqnarray}
However, these does not satisfy the compatibility condition
$L^{*}_{m,n}\hat{L}^{*}_{m+1,n}=\hat{L}^{*}_{m,n}L^{*}_{m,n+1}$ for 
the translation $\Phi_{m,n} \longrightarrow \Phi_{m+1,n+1}$.
As in \cite{bobpink}, we use the gauge transformation 
$\widetilde{\Phi}_{m,n}=\Phi_{m,n}\exp(\varepsilon_{m,n}E_{3})$.
We then see that $L^{\pm}_{m,n}$ and $\hat{L}^{\pm}_{m,n}$ are transformed
to ${\cal L}^{\pm}_{m,n}=\exp(-\varepsilon_{m,n}E_{3})
L^{\pm}_{m,n}\exp(\varepsilon_{m+1,n}E_{3})$ and
$\hat{\cal L}^{\pm}_{m,n}=\exp(-\varepsilon_{m,n}E_{3})
\hat{L}^{\pm}_{m,n}\exp(\varepsilon_{m,n+1}E_{3})$, respectively.
We choose $\varepsilon_{m,n}$'s as follows.
\begin{equation*}
\varepsilon_{m,n}=\frac{1}{4}\hat{w}_{m,n}, \quad
\varepsilon_{m+1,n}=\frac{1}{4}\hat{w}_{m+1,n},\quad
\varepsilon_{m,n+1}=\frac{1}{4}\hat{w}_{m,n}
+\frac{1}{4}\hat{w}_{m,n+1}.
\end{equation*}
We then have the following.
\begin{eqnarray}
\quad
\left\{
\begin{array}{ll}
{\cal L}^{\pm}_{m,n}&=
\left(\begin{matrix}\cos(\frac{\nu^{1}}{2})
e^{-\frac{i}{2}(\hat{w}_{m+1,n}
-\hat{w}_{m,n})}&
\pm\sin(\frac{\nu^{1}}{2})\\
\mp\sin(\frac{\nu^{1}}{2})&
\cos(\frac{\nu^{1}}{2})e^{\frac{i}{2}(\hat{w}_{m+1,n}
-\hat{w}_{m,n})}\\
\end{matrix}\right), \\
&\\
\hat{\cal L}^{\pm}_{m,n}&=
\left(\begin{matrix}\cos(\frac{\nu^{2}}{2})
&\pm\sin(\frac{\nu^{2}}{2})e^{\frac{i}{2}(\hat{w}_{m,n+1}
+\hat{w}_{m,n})}\\
\mp\sin(\frac{\nu^{2}}{2})e^{-\frac{i}{2}(\hat{w}_{m,n+1}
+\hat{w}_{m,n})}&
\cos(\frac{\nu^{2}}{2})\\
\end{matrix}\right) . \\
\end{array}
\right.
\end{eqnarray}
For the compatibility condition 
${\cal L}^{*}_{m,n}
\hat{\cal L}^{*}_{m+1,n}=\hat{\cal L}^{*}_{m,n}
{\cal L}^{*}_{m,n+1}$, we have the following.
\begin{Proposition}{\rm (cf. }\cite{bobpink}{\rm )}
Let $\hat{\gamma}$ be the one defined in the 
sine-Gordon equation~{\rm (\ref{sineGordon})}. 
\item{\rm (1)} 
${\cal L}^{\pm}_{m,n}
\hat{\cal L}^{\pm}_{m+1,n}=\hat{\cal L}^{\pm}_{m,n}
{\cal L}^{\pm}_{m,n+1}$
holds if and only if
$\hat{w}_{m,n}$ is a solution of the discrete sine-Gordon equation
with $\hat{\gamma}=-\tan(\frac{\nu^{1}}{2})\tan(\frac{\nu^{2}}{2})$.
\item{\rm (2)}
${\cal L}^{\pm}_{m,n}
\hat{\cal L}^{\mp}_{m+1,n}=\hat{\cal L}^{\mp}_{m,n}
{\cal L}^{\pm}_{m,n+1}$
holds if and only if
$\hat{w}_{m,n}$ is a solution of the discrete sine-Gordon equation
with $\hat{\gamma}=\tan(\frac{\nu^{1}}{2})\tan(\frac{\nu^{2}}{2})$.

\end{Proposition}

{\it Proof}. We prove the case (1) only. The case (2) is similarily
proved.
For simplicity of the notation, we use 
the following abbreviation.
\begin{equation*}
A=\hat{w}_{m+1,n+1},\quad B=\hat{w}_{m,n},\quad
C=\hat{w}_{m+1,n},\quad D=\hat{w}_{m,n+1} .
\end{equation*}
We then see that 
${\cal L}^{\pm}_{m,n}
\hat{\cal L}^{\pm}_{m+1,n}=\hat{\cal L}^{\pm}_{m,n}
{\cal L}^{\pm}_{m,n+1}$ holds if and only if
\begin{eqnarray*}
\cos(\frac{\nu^{1}}{2})\cos(\frac{\nu^{2}}{2})
\left(e^{-\frac{i}{2}(A-D)}-e^{-\frac{i}{2}(C-B)}\right)
=\sin(\frac{\nu^{1}}{2})\sin(\frac{\nu^{2}}{2})
\left(e^{\frac{i}{2}(D+B)}-e^{-\frac{i}{2}(A+C)}\right),
\end{eqnarray*}
which is equivalent to
\begin{equation*}
\cos(\frac{\nu^{1}}{2})\cos(\frac{\nu^{2}}{2})
\left(e^{-\frac{i}{2}(A+B)}-e^{-\frac{i}{2}(C+D)}\right)
=\sin(\frac{\nu^{2}}{2})\sin(\frac{\nu^{2}}{2})
\left(1-e^{-\frac{i}{2}(A+B+C+D)}\right),
\end{equation*}
which is also equivalent to
\begin{eqnarray*}
\begin{array}{ll}
&\cos(\frac{1}{2}(A+B))-\cos(\frac{1}{2}(C+D))
-i\left(\sin(\frac{1}{2}(A+B))-\sin(\frac{1}{2}(C+D))
\right)\\
&=\tan(\frac{\nu^{1}}{2})\tan(\frac{\nu^{2}}{2})
\left(
1-\cos(\frac{1}{2}(A+B+C+D))+i\sin(
\frac{1}{2}(A+B+C+D))\right) .
\end{array}
\end{eqnarray*}
For $U+V=\frac{1}{2}(A+B), U-V=\frac{1}{2}(C+D)$, we have
$U=\frac{1}{4}(A+B+C+D), V=\frac{1}{4}(A+B-C-D)$.
Thus, the last equation is equivalent to
\begin{equation*}
-2\sin(U)\sin(V)-2i\cos(U)\sin(V)=
\tan(\frac{\nu^{1}}{2})\tan(\frac{\nu^{2}}{2})
\left(2\sin^{2}(U)+2i\sin(U)\cos(U)\right),
\end{equation*}
which is again equivalent to 
$-\sin(V)=\tan(\frac{\nu^{1}}{2})\tan(\frac{\nu^{2}}{2})\sin(U)$,
which is the discrete sine-Gordon equation stated in \S2 for
$\hat{\gamma}=-\tan(\frac{\nu^{1}}{2})\tan(\frac{\nu^{2}}{2})$.
For example, for \lq\lq{\rm dn}-solution\rq\rq,
$\tan(\frac{\nu^{1}}{2})\tan(\frac{\nu^{2}}{2})=
\frac{{\rm sn}(\frac{\gamma}{2}){\rm dn}(\frac{\gamma}{2})}{
{\rm cn}(\frac{\gamma}{2})}
\frac{{\rm sn}(\frac{\delta}{2}){\rm dn}(\frac{\delta}{2})}{
{\rm cn}(\frac{\delta}{2})}$. For \lq\lq {\rm cn}-solution\rq\rq,
$\tan(\frac{\nu^{1}}{2})\tan(\frac{\nu^{2}}{2})=k^{2}
\frac{{\rm sn}(\frac{\gamma}{2}){\rm cn}(\frac{\gamma}{2})}{
{\rm dn}(\frac{\gamma}{2})}
\frac{{\rm sn}(\frac{\delta}{2}){\rm cn}(\frac{\delta}{2})}{
{\rm dn}(\frac{\delta}{2})}$. For this, just use 
$\tan(\frac{\nu}{2})=\dfrac{\sin\nu}{1+\cos\nu}$ 
and each values $\sin\nu$ and $\cos\nu$ 
in \lq\lq {\rm dn}-solution\rq\rq or
\lq\lq {\rm cn}-solution\rq\rq.
\qed

The following is an example of discrete $K$-surface corresponding
to the \lq\lq dn-solution\rq\rq of the discrete sine-Gordon equation. 

\begin{Theorem}\label{ksurface1}
 The following $F_{m,n}$ gives a discrete
$K$-surface.
\begin{equation}
F_{m,n}=\left(\begin{matrix}
\frac{1}{k}(-1)^{n}\cos(\varphi_{m,n}){\rm dn}(\psi_{m,n})\\
&\\
\frac{1}{k}(-1)^{n}\sin(\varphi_{m,n}){\rm dn}(\psi_{m,n})\\
\displaystyle{-k\int_{0}^{\psi_{m,n}}{\rm sn}^{2}\psi d\psi
+k\left(m\int_{0}^{\gamma}{\rm sn}^{2}\psi d\psi +
n\int_{0}^{\delta}{\rm sn}^{2}\psi d\psi\right)}\\
\end{matrix}\right) ,
\end{equation}
where $\varphi_{m,n}=\alpha m +\beta n, 
\psi_{m,n}=\gamma m +\delta n$ with
$\cos(\alpha)={\rm dn}(\gamma), \sin(\alpha)=k {\rm sn}(\gamma),
\cos(\beta)=-{\rm dn}(\delta), \sin(\beta)=k {\rm sn}(
\delta)$. Moreover, we have the following properties.
\item{\rm (i)} $F_{m+1,n}-F_{m,n}=N_{m+1,n}\times N_{m,n}$,
\item{\rm (ii)} $F_{m,n+1}-F_{m,n}=-N_{m,n+1}\times N_{m,n}$,
where
\begin{equation}
N_{m,n}=\left(\begin{matrix}
(-1)^{n}\cos(\varphi_{m,n}){\rm sn}(\psi_{m,n}) \\
(-1)^{n}\sin(\varphi_{m,n}){\rm sn}(\psi_{m,n}) \\
-{\rm cn}(\psi_{m,n})\\
\end{matrix}\right)
\end{equation}
\end{Theorem}

{\it Proof}. The form of $F_{m,n}$ 
follows from Theorem~\ref{dnsolution} and
Corollary~\ref{dncorollary}, where we write $B_{m}$'s as $N_{m,n}$.
Therefore, we see that the properties (i) and (ii) hold. 
We then see that both of $F_{m+1,n}-F_{m,n}$ and $F_{m,n+1}-F_{m,n}$
are perpendicular to $N_{m,n}$. Hence, there is a plane ${\cal P}_{m,n}$
such that three points $F_{m,n}, F_{m+1,n}, F_{m,n+1}$ lies on
${\cal P}_{m,n}$. However, since we also have
$F_{m-1,n}-F_{m,n}=N_{m-1,n}\times N_{m,n}\perp N_{m,n}$ and 
$F_{m,n-1}-F_{m,n}=N_{m,n-1}\times N_{m,n}\perp N_{m,n}$,
the rest two points $F_{m-1,n}$ and $F_{m,n-1}$ also lies on 
${\cal P}_{m,n}$. Thus, $F_{m,n}$ is a discrete $K$-surface.
\qed

For $F_{m,n}$ in Theorem~\ref{ksurface1}, while ${\cal L}^{\pm}_{m,n}$
corresponds to the curvature ${\cal K}_{m+1}
=\frac{1}{2}\hat{w}_{m+2}-\frac{1}{2}\hat{w}_{m}$, 
$\hat{\cal L}^{\pm}_{m,n}$ corresponds to the curvature
${\cal K}_{m+1}
=-\frac{1}{2}\hat{w}_{m+2}+\frac{1}{2}\hat{w}_{m}$. 

The next example of a discrete $K$-surface corresponds
to the \lq\lq cn-solution\rq\rq of the discrete sine-Gordon equation. 

We have the following theorem.

\begin{Theorem}\label{ksurface2} 
The following $\hat{F}_{m,n}$ gives a discrete
$K$-surface.
\begin{equation}
\hat{F}_{m,n}=\left(\begin{matrix}
(-1)^{n}k \cos(\varphi_{m,n}){\rm cn}(\psi_{m,n})\\
&\\
(-1)^{n}k \sin(\varphi_{m,n}){\rm cn}(\psi_{m,n})\\
\displaystyle{-k^{2}\int_{0}^{\psi_{m,n}}{\rm sn}^{2}\psi d\psi
+k^{2}\left(m\int_{0}^{\gamma}{\rm sn}^{2}\psi d\psi +
n\int_{0}^{\delta}{\rm sn}^{2}\psi d\psi\right)}\\
\end{matrix}\right) ,
\end{equation}
where $\varphi_{m,n}=\alpha m +\beta n, 
\psi_{m,n}=\gamma m +\delta n$ with
$\cos(\alpha)={\rm cn}(\gamma), \sin(\alpha)={\rm sn}(\gamma),
\cos(\beta)=-{\rm cn}(\delta), \sin(\beta)={\rm sn}(
\delta)$. Moreover, we have the following properties.
\item{\rm (i)} $\hat{F}_{m+1,n}-\hat{F}_{m,n}=N_{m+1,n}\times N_{m,n}$,
\item{\rm (ii)} $\hat{F}_{m,n+1}-\hat{F}_{m,n}=-N_{m,n+1}\times N_{m,n}$,
where
\begin{equation}
N_{m,n}=\left(\begin{matrix}
(-1)^{n}k \cos(\varphi_{m,n}){\rm sn}(\psi_{m,n}) \\
(-1)^{n}k \sin(\varphi_{m,n}){\rm sn}(\psi_{m,n}) \\
-{\rm dn}(\psi_{m,n})\\
\end{matrix}\right)
\end{equation}
\end{Theorem}

{\it Proof}. This follows from Theorem~\ref{cnsolution} and
Corollary~\ref{cncorollary}.
\qed
\vspace{15pt}

{\bf [Periodicity]} 
\begin{enumerate}
\item For $F_{m,n}$, we may utilize the arguments of the periodicity in \S4.
\begin{enumerate}
\item The case of $\gamma=\delta=K$. Assume that $m=2p$ and
$n=2q$, where $p, q\in\mathbb{N}$. Since $\cos(\alpha)={\rm dn}(K)=
k^{\prime}$ and $\cos(\beta)=-{\rm dn}(K)=-k^{\prime}$,
we must have $\alpha=\frac{\pi}{p}$ and $\beta=\frac{q-1}{q}\pi$,
thus we have $p=q$. Choosing $k^{\prime}$ so that $\cos(\frac{\pi}{p})
=k^{\prime}$ for any $p$ with $p> 2$. We then see that
$F_{m+2p, n+2p}=F_{m,n}$ for any $m, n\in\mathbb{N}$ and any 
$p> 2$.
\item The case of $\gamma=K, \delta=2K$. In this case, 
$\cos(\beta)=-{\rm dn}(2K)=-1$, hence $\beta=\pi$. Let $p> 2$ be the 
same one as in (i). We see that $F_{m+2p, n+2}=F_{m,n}$ for any
$m,n\in\mathbb{N}$.
\item The case of $\gamma=2K$ means $F_{m+1,n}=F_{m,n}$.
\end{enumerate}
\item 
For $\hat{F}_{m,n}$, we may utilize the arguments of the periodicity in \S5.
\begin{enumerate}
\item The case of $\gamma=\delta=K$. 
In this case, we have $\alpha=\beta=\frac{\pi}{2}$.
We then have $\varphi_{m,n}=\frac{\pi}{2}(m+n)$ and
$\psi_{m,n}=K(m+n)$. We see that
$\hat{F}_{m+4,n}=\hat{F}_{m,n}, 
\hat{F}_{m,n+4}=\hat{F}_{m,n}$ and $\hat{F}_{m+2,n+2}
=\hat{F}_{m,n}$
for any $m,n\in\mathbb{N}$. 
\vspace{10pt}
\item The case of $\gamma=K, \delta=2K$.
We then have $\alpha=\frac{\pi}{2}$ and $\beta=2\pi$. 
We see that $\hat{F}_{m+4,n}=\hat{F}_{m,n}, 
\hat{F}_{m,n+2}=\hat{F}_{m,n}$ and
$\hat{F}_{m+4,n+2}=\hat{F}_{m,n}$ for any $m,n\in\mathbb{N}$.
\item The case of $\gamma=2K, \delta=K$.
We then have $\alpha=\pi$ and $\beta=\frac{\pi}{2}$.
We see that $\hat{F}_{m+2,n}=\hat{F}_{m,n}, 
\hat{F}_{m,n+4}=\hat{F}_{m,n}$ and 
$\hat{F}_{m+1,n+2}=\hat{F}_{m,n}$ for any $m,n\in\mathbb{N}$.
\end{enumerate}
\end{enumerate}

\vspace{20pt}

\section{Appendix}
The following formulae are well known.
\begin{Lemma}\label{App1}
\begin{eqnarray*}
{\rm (i)}\left\{
\begin{array}{ll}
\theta_{3}(x\pm y)\theta_{0}(x\mp y)\theta_{3}(0)\theta_{0}(0)
&=\theta_{3}(x)\theta_{0}(x)\theta_{3}(y)\theta_{0}(y)\mp
\theta_{1}(x)\theta_{2}(x)\theta_{1}(y)\theta_{2}(y), \\
\theta_{1}(x\pm y)\theta_{2}(x\mp y)\theta_{3}(0)\theta_{0}(0)
&=\theta_{1}(x)\theta_{2}(x)\theta_{3}(y)\theta_{0}(y)\pm
\theta_{3}(x)\theta_{0}(x)\theta_{1}(y)\theta_{2}(y),\\
\theta_{1}(x\pm y)\theta_{3}(x\mp y)\theta_{2}(0)\theta_{0}(0)
&=\theta_{1}(x)\theta_{3}(x)\theta_{2}(y)\theta_{0}(y)\pm
\theta_{2}(x)\theta_{0}(x)\theta_{1}(y)\theta_{3}(y),\\
\theta_{2}(x\pm y)\theta_{0}(x\mp y)\theta_{2}(0)\theta_{0}(0)
&=\theta_{2}(x)\theta_{0}(x)\theta_{2}(y)\theta_{0}(y)\mp
\theta_{1}(x)\theta_{3}(x)\theta_{1}(y)\theta_{3}(y) .\\
\end{array}
\right.
\end{eqnarray*}
\begin{eqnarray*}
\hskip 1cm
{\rm (ii)}\left\{
\begin{array}{ll}
\theta_{3}(x+y |\tau^{\prime}) \theta_{3}(x-y |\tau^{\prime})
&=\theta_{3}(2x |2\tau^{\prime}) \theta_{3}(2y |2\tau^{\prime})
+\theta_{2}(2x |2\tau^{\prime}) \theta_{2}(2y |2\tau^{\prime}),\\
\theta_{0}(x+y |\tau^{\prime}) \theta_{0}(x-y |\tau^{\prime})
&=\theta_{3}(2x |2\tau^{\prime}) \theta_{3}(2y |2\tau^{\prime})
-\theta_{2}(2x |2\tau^{\prime}) \theta_{2}(2y |2\tau^{\prime}),\\
\theta_{2}(x+y |\tau^{\prime}) \theta_{2}(x-y |\tau^{\prime})
&=\theta_{2}(2x |2\tau^{\prime}) \theta_{3}(2y |2\tau^{\prime})
+\theta_{3}(2x |2\tau^{\prime}) \theta_{2}(2y |2\tau^{\prime}),\\
\theta_{1}(x+y |\tau^{\prime}) \theta_{1}(x-y |\tau^{\prime})
&=\theta_{3}(2x |2\tau^{\prime}) \theta_{2}(2y |2\tau^{\prime})
-\theta_{2}(2x |2\tau^{\prime}) \theta_{3}(2y |2\tau^{\prime}).\\
\end{array}
\right.
\end{eqnarray*}
\end{Lemma}

Using Lemma~\ref{App1} and considering the cases of 
$2x=2y=v^{\pm}_{m}(\lambda,z)$
or $2x=v^{+}_{m}(\lambda,z), 2y=v^{-}_{m}(\lambda,z)$, we obtain the following.

\begin{Lemma}\label{App2}
For $v^{\pm}_{m}(\lambda,z)=v_{m}\pm\frac{\lambda}{kK^{\prime}}
+\frac{i}{kK^{\prime}}z, v_{m}(z)=v_{m}+\frac{i}{kK^{\prime}}z$,
we have
\begin{eqnarray*}
\left\{
\begin{array}{ll}
\theta_{3}(v^{\pm}_{m}(\lambda,z)|\tau^{\prime}) \theta_{3}(0|\tau^{\prime})
&=\theta_{3}^{2}(v^{\pm}_{m}(\lambda,z) |2\tau^{\prime}) +
\theta_{2}^{2}(v^{\pm}_{m}(\lambda,z)|2\tau^{\prime}),\\
\theta_{0}(v^{\pm}_{m}(\lambda,z)|\tau^{\prime}) \theta_{0}(0|\tau^{\prime})
&=\theta_{3}^{2}(v^{\pm}_{m}(\lambda,z) |2\tau^{\prime}) -
\theta_{2}^{2}(v^{\pm}_{m}(\lambda,z)|2\tau^{\prime}),\\
\theta_{2}(v^{\pm}_{m}(\lambda,z)|\tau^{\prime}) \theta_{2}(0|\tau^{\prime})
&=2\theta_{2}(v^{\pm}_{m}(\lambda,z) |2\tau^{\prime})
\theta_{3}(v^{\pm}_{m}(\lambda,z)|2\tau^{\prime}),\\
\theta_{3}(v_{m}(z) |\tau^{\prime}) \theta_{3}(\frac{\lambda}{kK^{\prime}}|\tau^{\prime})
&=\theta_{3}(v^{+}_{m}(\lambda,z)|2\tau^{\prime}) 
\theta_{3}(v^{-}_{m}(\lambda,z)|2\tau^{\prime})
+\theta_{2}(v^{+}_{m}(\lambda,z)|2\tau^{\prime})
\theta_{2}(v^{-}_{m}(\lambda,z)|2\tau^{\prime}),\\
\theta_{0}(v_{m}(z) |\tau^{\prime}) \theta_{0}(\frac{\lambda}{kK^{\prime}}|\tau^{\prime})
&=\theta_{3}(v^{+}_{m}(\lambda,z)|2\tau^{\prime}) 
\theta_{3}(v^{-}_{m}(\lambda,z)|2\tau^{\prime})
-\theta_{2}(v^{+}_{m}(\lambda,z)|2\tau^{\prime})
\theta_{2}(v^{-}_{m}(\lambda,z)|2\tau^{\prime}),\\
\theta_{2}(v_{m}(z) |\tau^{\prime}) \theta_{2}(\frac{\lambda}{kK^{\prime}}|\tau^{\prime})
&=\theta_{2}(v^{+}_{m}(\lambda,z)|2\tau^{\prime}) 
\theta_{3}(v^{-}_{m}(\lambda,z)|2\tau^{\prime})
+\theta_{3}(v^{+}_{m}(\lambda,z)|2\tau^{\prime})
\theta_{2}(v^{-}_{m}(\lambda,z)|2\tau^{\prime}),\\
\theta_{1}(v_{m}(z) |\tau^{\prime}) \theta_{1}(\frac{\lambda}{kK^{\prime}}|\tau^{\prime})
&=\theta_{3}(v^{+}_{m}(\lambda,z)|2\tau^{\prime}) 
\theta_{2}(v^{-}_{m}(\lambda,z)|2\tau^{\prime})
-\theta_{2}(v^{+}_{m}(\lambda,z)|2\tau^{\prime})
\theta_{3}(v^{-}_{m}(\lambda,z)|2\tau^{\prime}).\\
\end{array}
\right.
\end{eqnarray*}
\end{Lemma}

\begin{Lemma}\label{App3} 
Let $\tau^{\prime}=i\frac{K}{K^{\prime}},
\psi_{m}=m\gamma+\beta t$ and 
$v_{m}=\frac{\psi_{m}-K}{2iK^{\prime}}$.
\item{\rm (i)} $\dfrac{\theta_{0}(v_{m}|\tau^{\prime}) \theta_{3}(0|\tau^{\prime})}{\theta_{3}(v_{m}|\tau^{\prime}) \theta_{0}(0|\tau^{\prime})}={\rm sn}(\psi_{m})$,\quad
{\rm (ii)} $\dfrac{\theta_{1}(v_{m}|\tau^{\prime}) \theta_{2}(0|\tau^{\prime})}{\theta_{3}(v_{m}|\tau^{\prime}) \theta_{0}(0|\tau^{\prime})}=i{\rm cn}(\psi_{m})$,
\item{\rm (iii)} 
$\dfrac{\theta_{2}(v_{m}|\tau^{\prime}) \theta_{2}(0|\tau^{\prime})}{\theta_{3}(v_{m}|\tau^{\prime}) \theta_{3}(0|\tau^{\prime})}={\rm dn}(\psi_{m})$.
\end{Lemma}

The Weierstrass ${\frak p}$-function defined in (0.4) has two periods $2\omega=2K^{\prime}$ and $2\omega^{\prime}=
iK+K^{\prime}$. Remark that $\dfrac{\omega^{\prime}}{\omega}
=\dfrac{1}{2}+\dfrac{\tau^{\prime}}{2}$.

\begin{Lemma}\label{App4}
Let $\zeta_{\rm W}$ be the Weierstrass 
$\zeta$-function satisfying 
$\dfrac{d\zeta_{\rm W}}{dw}=-{\frak p}(w)$.
We then have
\item{\rm (i)} $\zeta_{\rm W}(\frac{\omega}{2})
-\frac{1}{2}\zeta_{\rm W}(\omega)=k$,\quad
{\rm (ii)} ${\frak p}(\frac{\omega}{2})+\dfrac{\zeta_{\rm W}(\omega)}{\omega}=\dfrac{2E^{\prime}}{K^{\prime}}$.
\end{Lemma}

{\it Proof}. (i) Use the following formulae.
\begin{equation*}
\zeta_{\rm W}(u+v)-\zeta_{\rm W}(u)-\zeta_{\rm W}(v)=
\dfrac{1}{2}\dfrac{{\frak p}^{\prime}(u)-{\frak p}^{\prime}(v)}
{{\frak p}(u)-{\frak p}(v)}.
\end{equation*}
It follows from (0.4) that 
${\frak p}^{\prime}(w)=4k{\rm cn}(2iw+iK^{\prime})\left(
{\frak p}(w)-e_{1}\right)$, where $e_{1}=\frac{2}{3}(2k^{2}-1)$.
Therefore, we see that
\begin{eqnarray*}
\begin{array}{ll}
\dfrac{1}{2}\displaystyle{\lim_{u\rightarrow \frac{\omega}{2}}
\dfrac{{\frak p}^{\prime}(u)-{\frak p}^{\prime}(\omega)}{
{\frak p}(u)-{\frak p}(\omega)}}
&=\dfrac{1}{2}\displaystyle{\lim_{u\rightarrow \frac{\omega}{2}}4k{\rm cn}(2iu+iK^{\prime})}\\
&=\dfrac{1}{2}\cdot 4k{\rm cn}(2iK^{\prime})=-2k,\\
\end{array}
\end{eqnarray*}
which, together with the formulae above for $u=\frac{\omega}{2},
v=-\omega$, implies the result (i). For (ii), first note that
$-2i\omega^{\prime}+iK^{\prime}=K$. We calculate the integral
\begin{eqnarray*}
\begin{array}{ll}
\displaystyle{\int_{-\frac{\omega}{2}}^{-\omega^{\prime}}
{\frak p}(w) dw}&=\displaystyle{\dfrac{1}{2i}\int_{0}^{K}
\left(({\rm dn}(u)-ik{\rm sn}(u))^{2}+e_{1}\right) du}\\
&=\displaystyle{\dfrac{1}{2i}\int_{0}^{K}\left(
2{\rm dn}^{2}(u)+(e_{1}-1)-2ik{\rm sn}(u){\rm dn}(u)\right)du}\\
&=\displaystyle{\dfrac{1}{2i}\left\{
2\int_{0}^{K}{\rm dn}^{2}(u) du +(e_{1}-1)K+2ik\left[{\rm cn}(u)\right]_{0}^{K}\right\}}\\
&=-iE-\dfrac{iK}{2}(e_{1}-1)-k .\\
\end{array}
\end{eqnarray*} 
On the other hand, we have 
$\displaystyle{\int_{-\frac{\omega}{2}}^{-\omega^{\prime}}
{\frak p}(w) dw=\left[-\zeta_{\rm W}(w)\right]_{-\frac{\omega}{2}}^{-\omega^{\prime}}=\zeta_{\rm W}(\omega^{\prime})-
\zeta_{\rm W}(\frac{\omega}{2})}$,
which together with the Legendre identity $\zeta_{\rm W}(\omega)\dfrac{\omega^{\prime}}{\omega}-\dfrac{\pi i}{2\omega}
=\zeta_{\rm W}(\omega^{\prime})$, $\dfrac{\omega^{\prime}}{\omega}=\dfrac{1}{2}+\dfrac{\tau^{\prime}}{2}$ and the above,
implies that 
\begin{equation}
\zeta_{\rm W}(\frac{\omega}{2})-\dfrac{1}{2}\zeta_{\rm W}(\omega)=k+iE+\dfrac{iK}{2}(e_{1}-1)+\dfrac{\tau^{\prime}}{2}
\zeta_{\rm W}(\omega)-\dfrac{\pi i}{2K^{\prime}} ,
\end{equation}
which ,together with result (i), yields that
\begin{equation*}
\zeta_{\rm W}(\omega)=\dfrac{\pi}{K}-\dfrac{2K^{\prime}E}{K}
+K^{\prime}(1-e_{1}).
\end{equation*}
Therefore, from the Legendre formulae
$KE^{\prime}+EK^{\prime}-KK^{\prime}=\dfrac{\pi}{2}$ 
we have
\begin{equation}
\dfrac{\zeta_{\rm W}(\omega)}{\omega}
=\dfrac{\pi}{KK^{\prime}}-\dfrac{2E}{K}+(1-e_{1})
=\dfrac{2E^{\prime}}{K^{\prime}}-(e_{1}+1) .
\end{equation}
On the other hand, we see that
\begin{equation*}
{\frak p}(\frac{\omega}{2})=({\rm dn}(2iK^{\prime})-ik
{\rm sn}(2iK^{\prime}))^{2}+e_{1}=e_{1}+1.
\end{equation*}
Therefore, we have the result (ii).
\qed

\begin{Lemma}\label{App7}
For $\psi_{m}=m\gamma+\beta t$, 
the following equations hold identically.
\item{\rm (i)} 
${\rm dn}(\gamma){\rm sn}(\psi_{m}){\rm sn}(\psi_{m+1})+{\rm cn}(\psi_{m}){\rm cn}(\psi_{m+1})={\rm cn}(\gamma)$,
\item{\rm (ii)}
$k^{2}{\rm cn}(\gamma){\rm sn}(\psi_{m}){\rm sn}(\psi_{m+1})+
{\rm dn}(\psi_{m}){\rm dn}(\psi_{m+1})={\rm dn}(\gamma)$ ,
\item{\rm (iii)}
$k^{2}{\rm sn}(\gamma){\rm sn}^{2}(\psi_{m+1})+
{\rm dn}(\gamma){\rm sn}(\psi_{m+1}){\rm cn}(\psi_{m}){\rm dn}(\psi_{m+1})
-{\rm sn}(\psi_{m}){\rm cn}(\psi_{m+1}){\rm dn}(\psi_{m+1})=
{\rm sn}(\gamma)$ ,
\item{\rm (iv)} 
${\rm sn}(\gamma){\rm dn}(\psi_{m+1})+{\rm sn}(\psi_{m}){\rm cn}(\psi_{m+1})
={\rm dn}(\gamma){\rm sn}(\psi_{m+1}){\rm cn}(\psi_{m})$ ,
\item{\rm (v)}
${\rm dn}(\gamma){\rm dn}(\psi_{m+1})+k^{2}{\rm sn}(\gamma)
{\rm sn}(\psi_{m+1}){\rm cn}(\psi_{m})={\rm dn}(\psi_{m})$ ,
\item{\rm (vi)}
\begin{equation*}
{\rm dn}(\gamma){\rm sn}(\psi_{m}){\rm cn}(\psi_{m}){\rm sn}(\psi_{m+1})
{\rm cn}(\psi_{m+1})+{\rm sn}(\gamma){\rm cn}(\psi_{m}){\rm dn}(\psi_{m})
{\rm sn}(\psi_{m+1})
={\rm cn}^{2}(\psi_{m}){\rm sn}^{2}(\psi_{m+1}) ,
\end{equation*}
\item{\rm (vii)} 
${\rm cn}(\gamma){\rm cn}(\psi_{m+1})+{\rm sn}(\gamma)
{\rm sn}(\psi_{m+1}){\rm dn}(\psi_{m})={\rm cn}(\psi_{m})$,
\item{\rm (viii)}
${\rm sn}(\gamma){\rm cn}(\psi_{m+1})+{\rm sn}(\psi_{m})
{\rm dn}(\psi_{m+1})={\rm cn}(\gamma){\rm sn}(\psi_{m+1})
{\rm dn}(\psi_{m})$,
\item{\rm (ix)}
\begin{equation*}
{\rm dn}^{2}(\gamma){\rm sn}(\psi_{m-1}){\rm sn}(\psi_{m+1})+
{\rm cn}(\psi_{m-1}){\rm cn}(\psi_{m+1})
+{\rm sn}^{2}(\gamma){\rm dn}(\psi_{m-1}){\rm dn}(\psi_{m+1})
={\rm cn}^{2}(\gamma) .
\end{equation*}
\end{Lemma}

{\it Proof}. 
These equations hold because of $\psi_{m-1}=\psi_{m}-\gamma,
\psi_{m+1}=\psi_{m}+\gamma$ and the addition formulae for
${\rm sn}, {\rm cn}$ and ${\rm dn}$-functions of the following.
\begin{eqnarray*}
\begin{array}{ll}
{\rm sn}(\psi_{m+1})&=\dfrac{{\rm cn}(\gamma){\rm dn}(\gamma){\rm sn}(\psi_{m})
+{\rm sn}(\gamma){\rm cn}(\psi_{m}){\rm dn}(\psi_{m})}{1-k^{2}{\rm sn}^{2}(\gamma){\rm sn}^{2}(\psi_{m})},\\
&\\
{\rm cn}(\psi_{m+1})&=\dfrac{{\rm cn}(\gamma){\rm cn}(\psi_{m})
-{\rm sn}(\gamma){\rm dn}(\gamma){\rm sn}(\psi_{m}){\rm dn}(\psi_{m})}{1-k^{2}{\rm sn}^{2}(\gamma){\rm sn}^{2}(\psi_{m})},\\
&\\
{\rm dn}(\psi_{m+1})&=\dfrac{{\rm dn}(\gamma){\rm dn}(\psi_{m})
-k^{2}{\rm sn}(\gamma){\rm cn}(\gamma){\rm sn}(\psi_{m})
{\rm cn}(\psi_{m})}{1-k^{2}{\rm sn}^{2}(\gamma){\rm sn}^{2}(\psi_{m})} .\\
\end{array}
\end{eqnarray*}
For $\psi_{m-1}$, replace $\gamma$ by $-\gamma$ in the equations above.
\qed

\begin{figure}[htbp]
\begin{center}
\includegraphics[height=168mm,width=140mm]{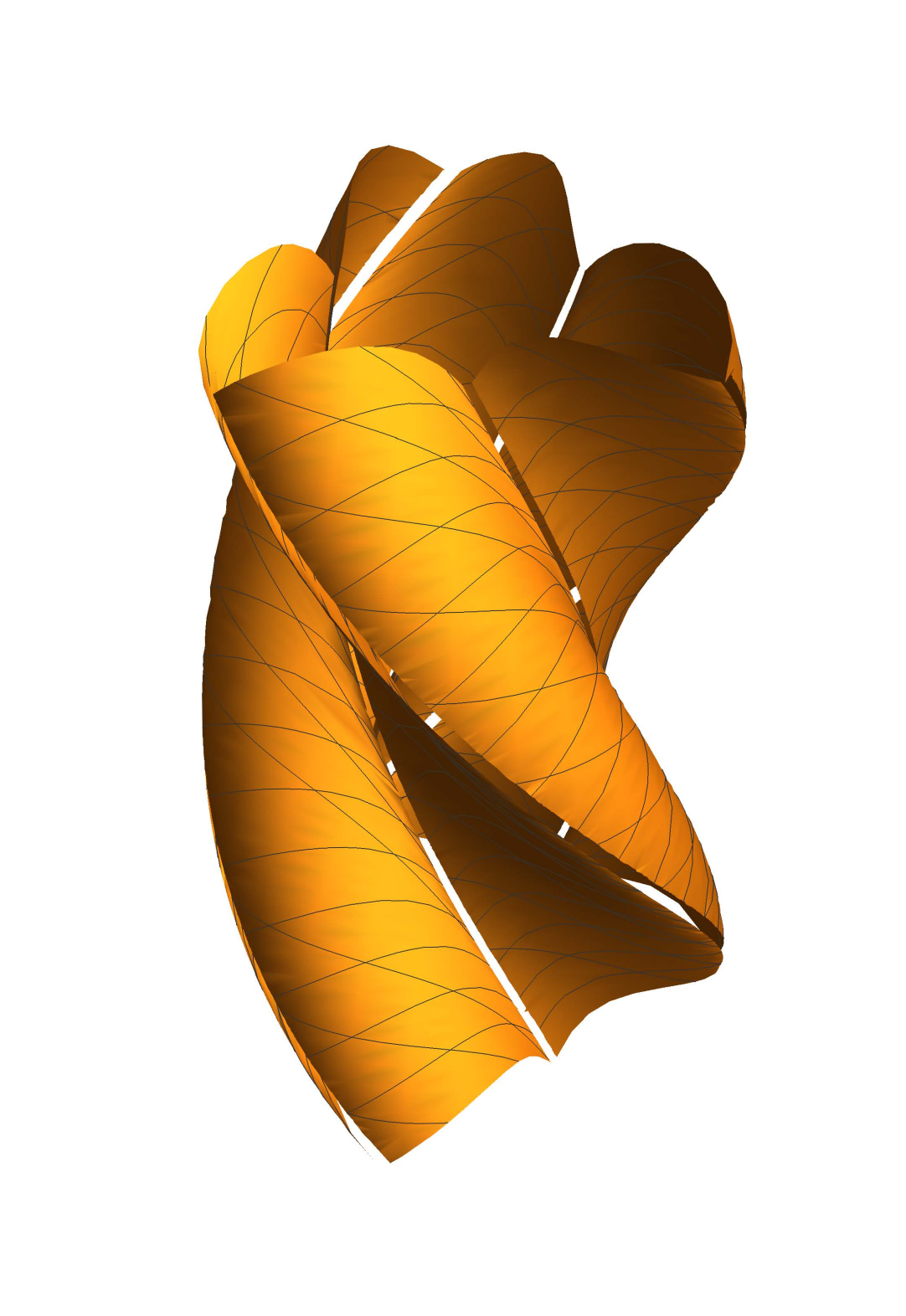}
\caption{\lq\lq {\rm dn}-solution semi-discrete surface\rq\rq; $\gamma=K$, n=3, m=6}
\end{center}
\end{figure}

\begin{figure}[htbp]
\begin{center}
\includegraphics[height=168mm,width=140mm]{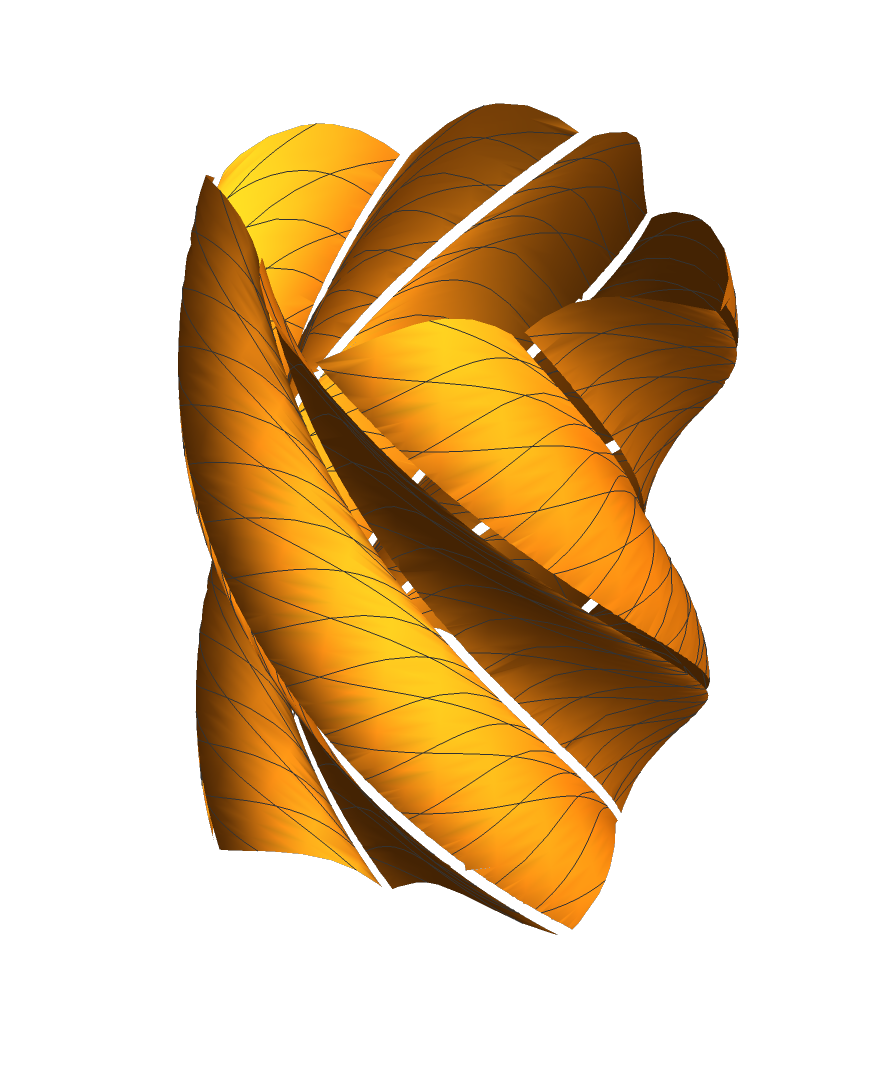}
\caption{\lq\lq {\rm dn}-solution semi-discrete surface\rq\rq; $\gamma=K$, n=4, m=8}
\end{center}
\end{figure}

\begin{figure}[htbp]
\begin{center}
\includegraphics[height=168mm,width=140mm]{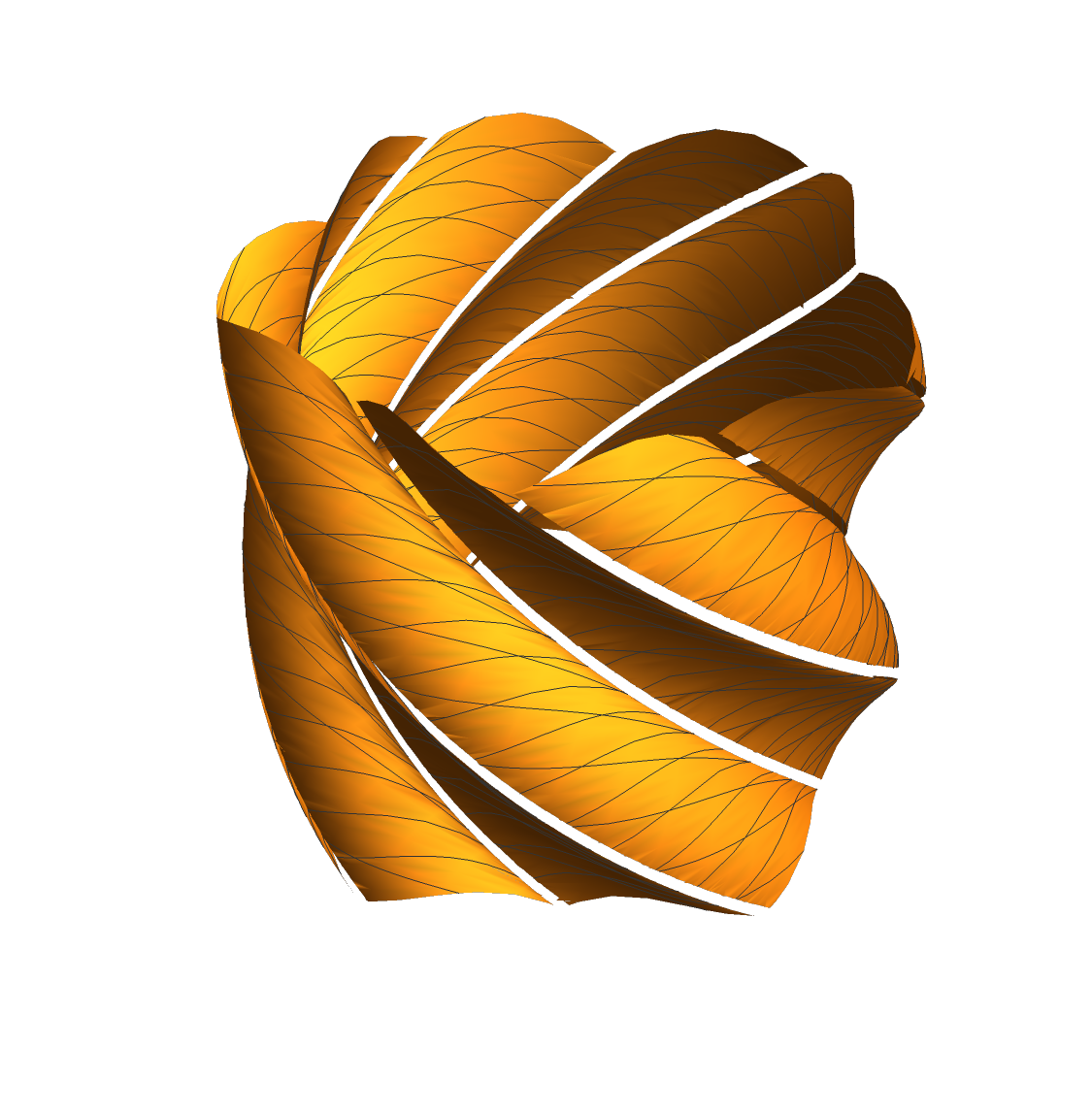}
\caption{\lq\lq {\rm dn}-solution semi-discrete surface\rq\rq; $\gamma=K$, n=5, m=10}
\end{center}
\end{figure}

\begin{figure}[htbp]
\begin{center}
\includegraphics[height=168mm,width=140mm]{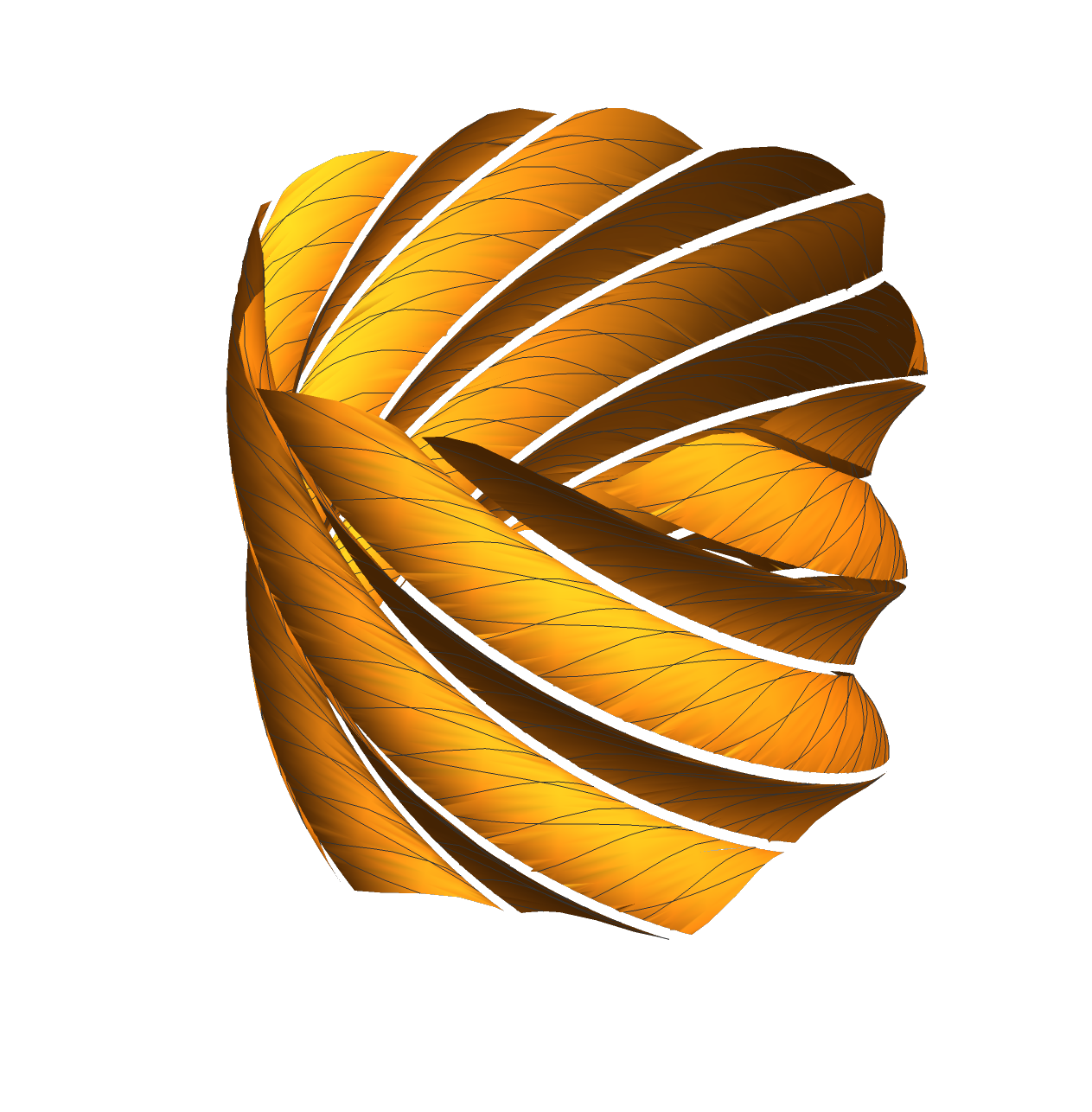}
\caption{\lq\lq {\rm dn}-solution semi-discrete surface\rq\rq; $\gamma=K$, n=6, m=12}
\end{center}
\end{figure}

\begin{figure}[htbp]
\begin{center}
\includegraphics[height=168mm,width=140mm]{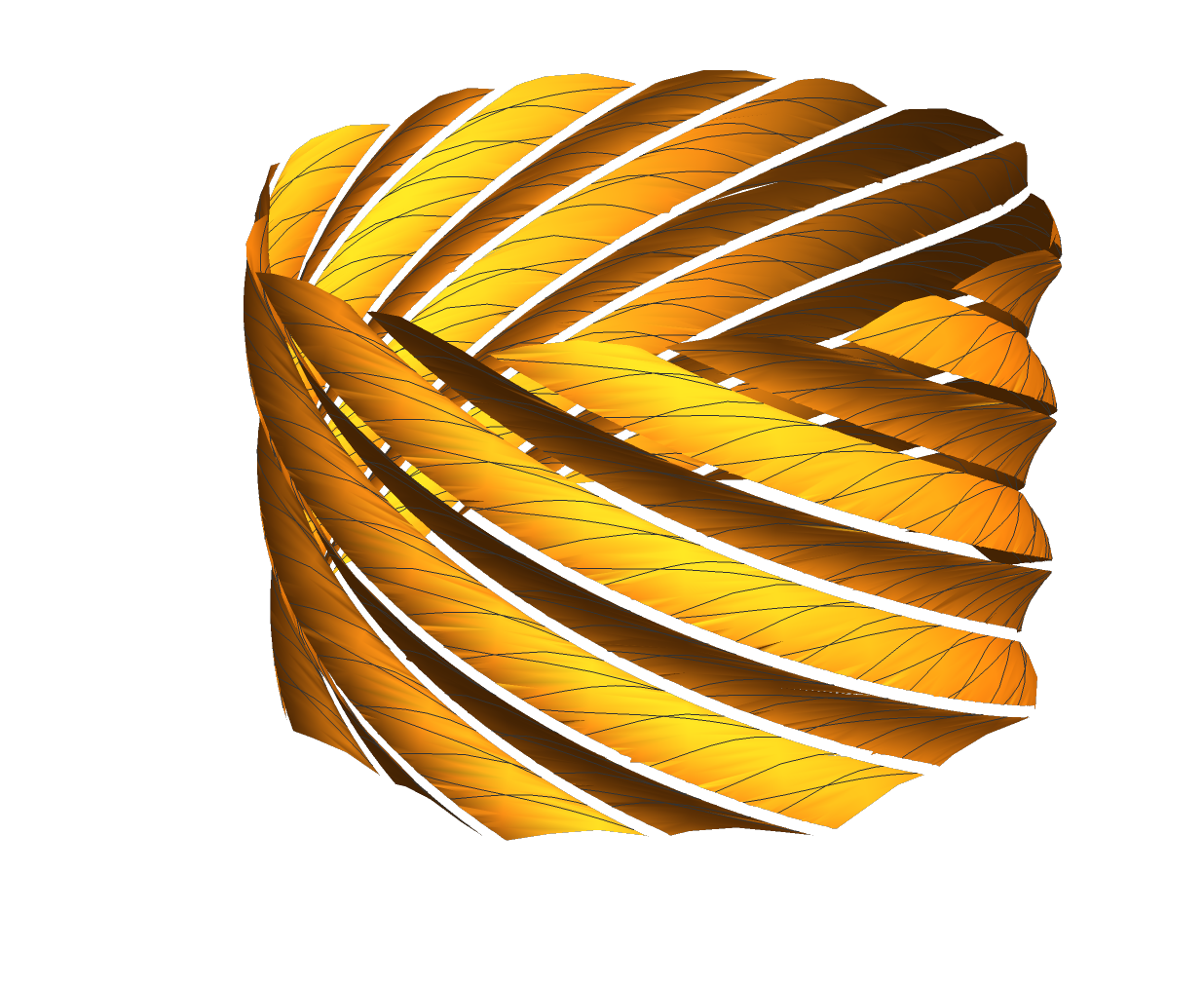}
\caption{\lq\lq {\rm dn}-solution semi-discrete surface\rq\rq; $\gamma=K$, n=8, m=16}
\end{center}
\end{figure}

\begin{figure}[htbp]
\begin{center}
\includegraphics[height=168mm,width=140mm]{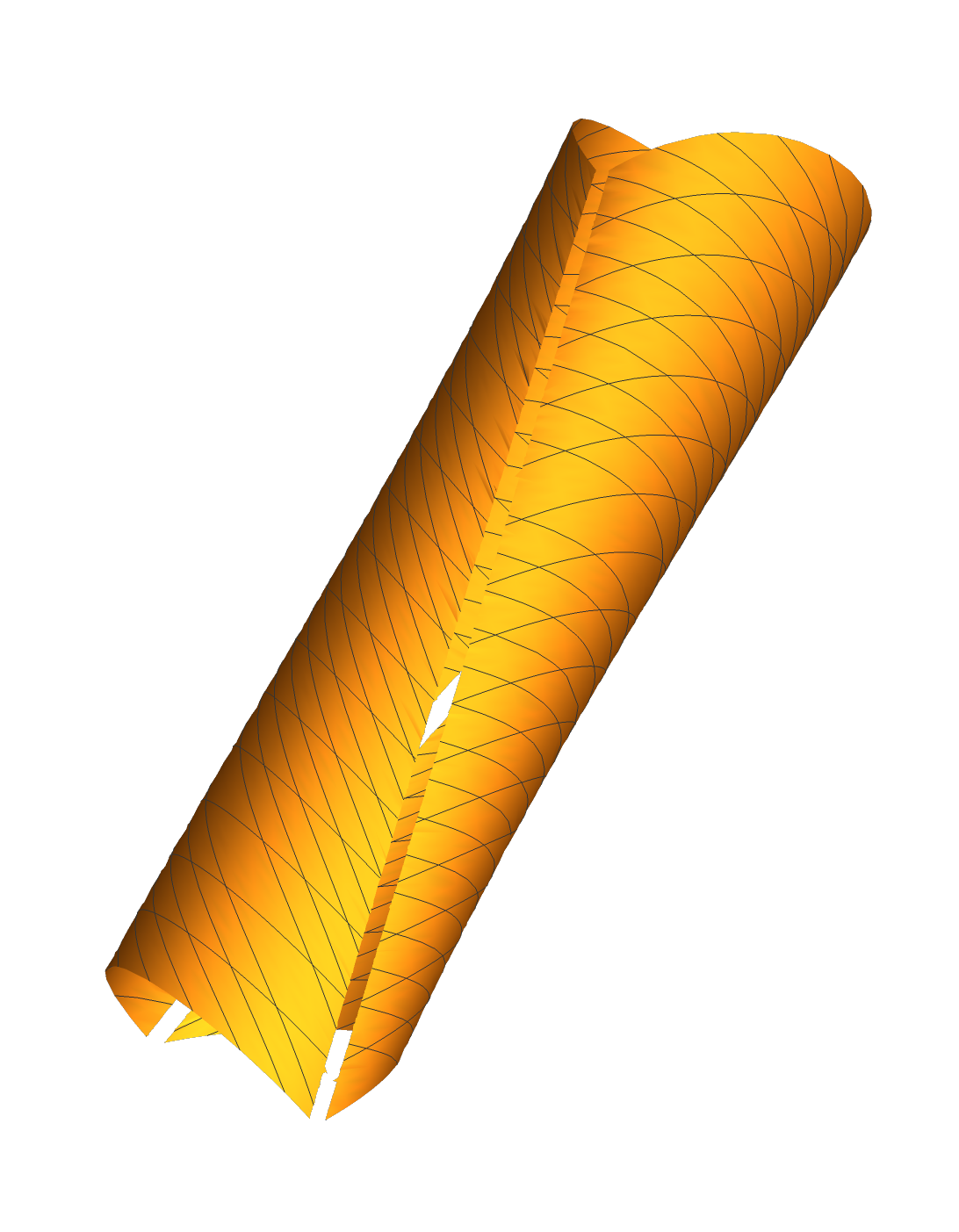}
\caption{\lq\lq {\rm cn}-solution semi-discrete surface\rq\rq; $\sin\alpha={\rm sn}(K/n)$, 
n=1, m=2}
\end{center}
\end{figure}

\begin{figure}[htbp]
\begin{center}
\includegraphics[height=168mm,width=140mm]{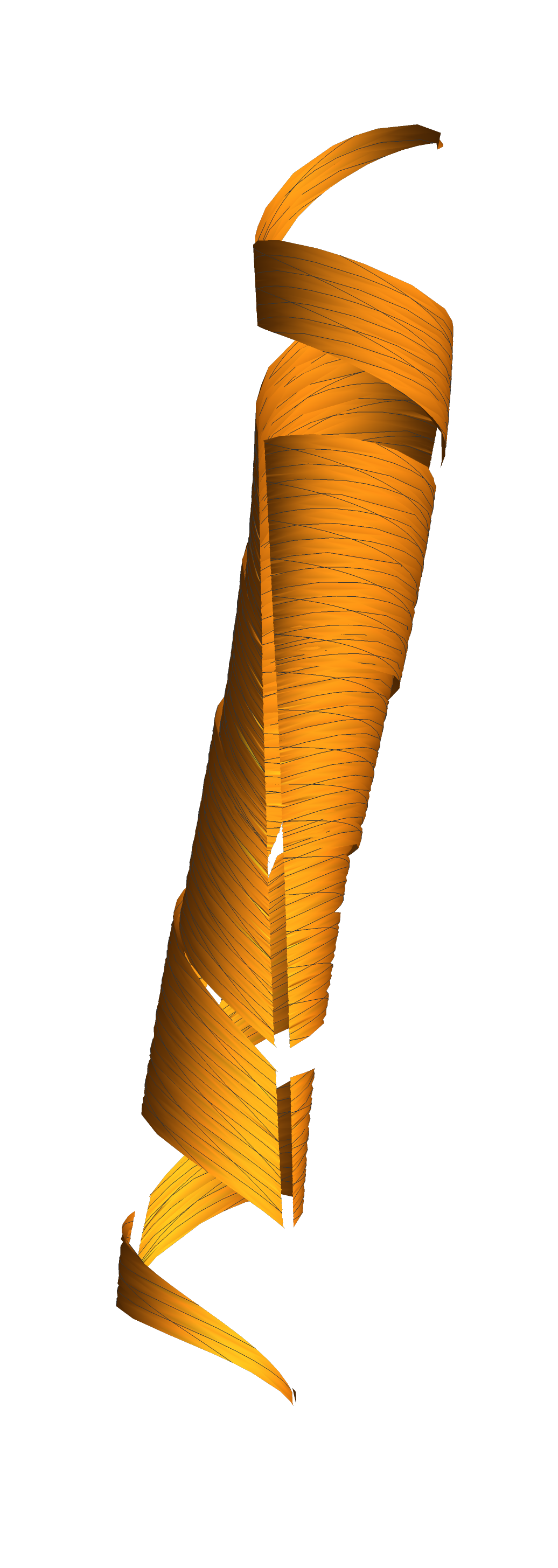}
\caption{\lq\lq {\rm cn}-solution semi-discrete surface\rq\rq; $\sin\alpha={\rm sn}(K/n)$, 
n=2, m=16}
\end{center}
\end{figure}

\begin{figure}[htbp]
\begin{center}
\includegraphics[height=168mm,width=140mm]{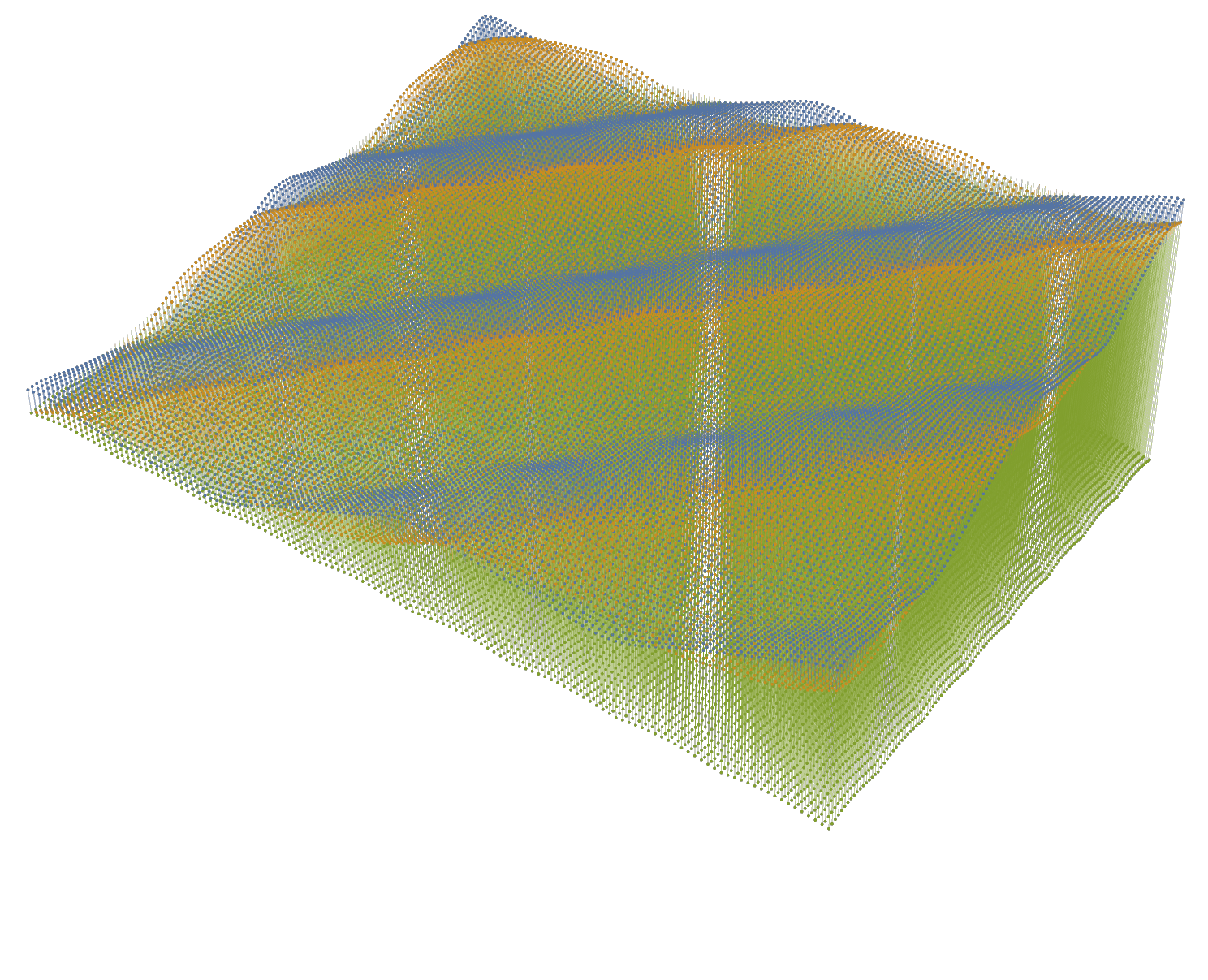}
\caption{\lq\lq Discrete {\rm K}-surface\rq\rq; $\sin\alpha=0.8\times{\rm sn}(K/16)$, 
m=n=128}
\end{center}
\end{figure}


\vspace{15pt}

\begin{thebibliography}{99}
\bibitem{bobpink} A. Bobenko and U. Pinkall,
{\it Discrete surfaces with constant negative 
Gaussian curvature and the Hirota equation}, 
Journal of Differential Geometry {\bf 43}(1996),
no.3, 527-611.
\bibitem{hirota} R. Hirota, {\it Nonlinear partial difference equations III
discrete sine-Gordon equation}, J. Phys. Soc. Japan {\bf 43},1977, 2079-2086.
\bibitem{ikmo} J. Inoguchi, K. Kajiwara, N. Matsuura and Y. Ohta,
{\it Discrete mKdV and discrete sine-Gordon flows on discrete space curves},
J. Physics A: Math. Theor. {\bf 47}, 2014, 235202.
\bibitem{kkp} S. Kaji, K. Kajiwara and H. Park,
{\it Linkage mechanisms governed by integrable deformations 
of discrete space curves},
Nonlinear Systems and Their Remarkable Mathematical Structures: 
Volume 2 (CRC Press,2019),356-381.
\bibitem{kks} S. Kaji, K. Kajiwara and S. Shigetomi,
{\it An explicit construction of Kaleidocycles},
a preprint.
\item K. Sogo, {\it Variational discretization of Euler's elastica problem},
J. Phys. Soc. Japan {\bf 75}, 2006, 064007.
\bibitem{udagawa} S. Udagawa, K. Kajiwara and J. Inoguchi.
{\it Solutions of sine-Gordon equation and its 
discretization}(written in Japanese), Bull. Liberal Arts and Sciences,
Nihon Univ. School of Medicine, no. {\bf 50},
2022, 6-26.
\item NIST, Digital Library of Mathematical FUnctions, 
{\it https://dlmf.nist.gov/}
\end{thebibliography}
\end{document}